\documentclass[a4paper,11pt,reqno]{article}
\usepackage{geometry}
\usepackage[normalem]{ulem}

\usepackage{amsthm} 

\usepackage{amsmath,amstext,amscd,amssymb,euscript,mathrsfs}
\usepackage{calrsfs}
\usepackage{epsfig}
\usepackage{mathtools}
\usepackage[colorlinks=false]{hyperref} 
\usepackage[english]{babel}
\usepackage{verbatim} 
\usepackage{enumerate}
\usepackage{bbm}
\usepackage{color}
\usepackage{wrapfig} 

\newcommand{\SortNoop}[1]{}

\usepackage{enumitem}
\setenumerate{label={\normalfont(\alph*)},topsep=4pt,itemsep=0pt} 

\theoremstyle{plain}
\newtheorem{theorem}{Theorem}[section]
\newtheorem{lemma}[theorem]{Lemma}
\newtheorem{proposition}[theorem]{Proposition}
\newtheorem{corollary}[theorem]{Corollary}

\theoremstyle{definition}
\newtheorem{definition}[theorem]{Definition}

\newtheorem{remark}[theorem]{Remark}

\newtheorem{obs}[theorem]{\!\!}

\newcommand{\eps}{\epsilon}
\newcommand{\e}{{\rm e}}

\newcommand{\arroweps}{\xrightarrow{\epsilon \downarrow 0}} 

\usepackage{times}
\usepackage{stmaryrd}

\newcommand{\1}{\mathbbm{1}}

\newcommand{\dd}{\, \mathrm{d}}

\newcommand{\T}{\mathbb T}
\newcommand{\E}{\mathbb E}

\newcommand{\Q}{\mathbb Q}
\newcommand{\Z}{\mathbb Z}
\newcommand{\R}{\mathbb R}
\newcommand{\N}{\mathbb N}

\newcommand{\cC}{\mathcal{C}}
\newcommand{\cD}{\mathcal{D}}

\newcommand{\cF}{\mathcal{F}}

\newcommand{\cS}{\mathcal{S}}
\newcommand{\cT}{\mathcal{T}}

\newcommand{\cW}{\mathcal{W}}

\newcommand{\fd}{\mathfrak{d}}

\newcommand{\fn}{\mathfrak{n}}

\newcommand{\fq}{\mathfrak{q}}

\newcommand{\fC}{\mathfrak{C}}   
\newcommand{\fD}{\mathfrak{D}}

\newcommand{\fL}{\mathfrak{L}}

\newcommand{\fR}{\mathfrak{R}}

\newcommand{\fU}{\mathfrak{U}}

\newcommand{\fX}{\mathfrak{X}}

\newcommand{\sH}{\mathscr{H}}

\renewcommand{\P}{\mathbb{P}}

\newcommand{\supp}{\operatorname{supp}} 

\usepackage{mathtools} 

\newcommand{\rD}{\mathrm{D}}

\newcommand{\para}{\varolessthan}

\newcommand{\reso}{\varodot}

\newcommand{\bxi}{{\boldsymbol{\xi}}}

\newcommand{\btheta}{{\boldsymbol{\theta}}}
\newcommand{\blambda}{{\boldsymbol{\lambda}}}

\newcommand{\I}{\mathbb{I}}

\newcommand{\limn}{\lim_{n\rightarrow \infty}}

\renewcommand{\epsilon}{\varepsilon}

\usepackage[disable]{todonotes} 

\usepackage[normalem]{ulem} 

\newenvironment{colorwil}
    {\color[rgb]{0,0.5,0}    }

\newcommand{\bw}{\bgroup\color[rgb]{0,0.5,0}}
\newcommand{\ew}{\egroup}

\newcommand{\bwil}{\begin{colorwil}}
\newcommand{\ewil}{\end{colorwil}}

\newenvironment{calc}    {\color{orange}     }    {     }

\newcommand{\cnewline}{\\}
\newcommand{\cand}{&}

\usepackage{environ}

\usepackage{xcolor}

\newcommand{\br}{\bgroup\color{red}}
\newcommand{\er}{\egroup}

\newcommand{\bb}{\bgroup\color{blue}}
\newcommand{\eb}{\egroup}

\RenewEnviron{calc}{}  
\renewcommand{\cnewline}{}  
\renewcommand{\cand}{}  

\RenewEnviron{hide}{}  

\usepackage{tocloft}
\begin{document}

\renewcommand{\thefootnote}{\Roman{footnote}}

\title{Longtime asymptotics of the two-dimensional parabolic Anderson model with white-noise potential}

\author{
\renewcommand{\thefootnote}{\Roman{footnote}}
Wolfgang K\"onig
\footnotemark[1]
\\
\renewcommand{\thefootnote}{\Roman{footnote}}
Nicolas Perkowski
\footnotemark[2]
\\
\renewcommand{\thefootnote}{\Roman{footnote}}
Willem van Zuijlen
\footnotemark[3]
}

\footnotetext[1]{TU Berlin and WIAS Berlin, Mohrenstra{\ss}e 39, 10117 Berlin, Germany, {\tt koenig@wias-berlin.de}}
\footnotetext[2]{FU Berlin, Arnimallee 7, 14195 Berlin, Germany {\tt perkowski@math.fu-berlin.de}
}
\footnotetext[3]{WIAS Berlin, Mohrenstra{\ss}e 39, 10117 Berlin, Germany, {\tt vanzuijlen@wias-berlin.de}}

\date{Oktober 15, 2021}

\maketitle 

\renewcommand{\thefootnote}{\arabic{footnote}} 


\begin{abstract}
We consider the parabolic Anderson model (PAM) $\partial_t u = \frac12 \Delta u + \xi u$ in $\R^2$ with a Gaussian (space) white-noise potential $\xi$. 
We prove that the almost-sure large-time asymptotic behaviour of the total mass at time $t$, written $U(t)$, is given by $\log U(t)\sim  \chi t \log t$ for $t \to \infty$, 
with the deterministic constant $\chi$ identified in terms of a variational formula. 
In earlier work of one of the authors this constant was used to describe the asymptotic behaviour $\blambda_1(Q_t)\sim\chi\log t$ of the principal eigenvalue $\blambda_1(Q_t)$ of the Anderson operator with Dirichlet boundary conditions on the box $Q_t= [-\frac{t}{2},\frac{t}{2}]^2$. 

\bigskip

\textbf{R\'esum\'e} Nous considérons le modèle parabolique d'Anderson (PAM) $\partial_t u = \frac12 \Delta u + \xi u$ dans $\R^2$ avec un potentiel de bruit blanc $\xi$ en espace. 
Nous prouvons que le comportement asymptotique presque sûr de la masse totale au temps $t$, écrite $U(t)$, est donné par $\log U(t)\sim \chi t \log t$ pour $t \to \infty$, 
avec une constante déterministe $\chi$ 
que nous identifions à l'aide 
d'une formule variationnelle. 
Cette constante a déjà été utilisée,
dans un travail antérieur de l'un des auteurs, 
pour décrire le comportement asymptotique 
$\blambda_1(Q_t)\sim\chi\log t$
de la valeur propre principale $\blambda_1(Q_t)$ de l'opérateur d'Anderson muni de conditions aux limites de Dirichlet sur la boîte $Q_t= [-\frac{t}{2},\frac{t}{2}]^2$.

\bigskip




%
\end{abstract}




\setcounter{tocdepth}{1} 
\tableofcontents

\section{Introduction and main results}\label{sec-Intro}

In this paper, we continue the programme of proving intermittency properties of the parabolic Anderson model in two dimensions with Gaussian white-noise potential. In Section \ref{sec-PAM} we introduce the model and its solution, in Section \ref{sec-purpose} we describe the difficulties and our strategy to resolve them, in Section \ref{sec-results} we formulate our main results, and in Section \ref{sec-discussion} we relate our findings to analogous earlier results that consider more regular random potentials. 

\subsection{The parabolic Anderson model and intermittency}\label{sec-PAM}

\noindent We consider the solution to the (Cauchy problem for the) heat equation with random potential, formally defined by
\begin{equation}\label{PAM}
\begin{aligned}
\partial_t u(t,x)&=\tfrac 12 \Delta u(t,x)+\xi(x) u(t,x),\qquad (t,x) \in (0,\infty) \times \R^2, \\ 
u(0,\cdot)&=\delta_0,
\end{aligned}
\end{equation}
This equation is called the {\em parabolic Anderson model (PAM)}, the random Schr\"odinger operator $\frac 12\Delta +\xi$ on the right-hand side is called the {\em Anderson Hamiltonian}. In our case, the potential $\xi $ is a Gaussian white noise in two dimensions, i.e., a distribution rather than a function. Here, the solution $u(t,\cdot)$ needs to be constructed via a renormalization procedure. Indeed, $\xi$ is replaced by a mollified version $\xi_\eps$ minus a correction $c_\eps \approx \frac 1{2\pi} \log \eps  $ and it is proved that the solution $u_\eps$ of the PAM with potential $\xi_\eps-c_\eps$ has a limit as $\eps\to 0$. This is the solution $u\colon [0,\infty)\times\R^2\to[0,\infty)$ that we will consider here. It was first constructed -- on the torus $\T^2$ instead of $\R^2$ -- by Hairer~\cite{Ha14} and by Gubinelli, Imkeller and Perkowski~\cite{Ha14, GuImPe15}, using the framework of regularity structures and the one of paracontrolled distributions, respectively. 
A construction on the entire $\R^2$ is due to Hairer and Labb\'e \cite{HaLa15} who realized that with a partial Cole-Hopf transform one can avoid using paracontrolled distributions or regularity structures.

As has been proved for certain choices of random potentials, the PAM displays an interesting {\em intermittency} effect. 
This means that it admits a highly pronounced concentration property on large space-time scales, which distinguishes it clearly from models that show a diffusive behaviour in the vicinity of the (functional) central limit theorem. Indeed, earlier investigations of the PAM on $\Z^d$ with i.i.d.\ potential, and  on $\R^d$ with regular potential, have revealed that the function $u(t,\cdot)$ is highly concentrated on few small islands 
 that are far from each other and carry most of the {\em total mass} of the solution,
\begin{equation}
U(t)=\int_{\R^2} u(t,x)\,{\rm d} x.
\end{equation}
The main source of conjectures and proofs has been -- under the assumption of sufficient regularity -- the Fourier expansion in terms of the eigenvalues $\lambda_1>\lambda_2\geq \lambda_3\geq \dots$ and corresponding $L^2$-orthonormal basis of eigenfunctions ${\rm e}_1,{\rm e}_2,{\rm e}_3 \dots$ of $\frac 12\Delta+\xi$,
\begin{equation}\label{eigenvalueexpansion}
u(t,x)=\sum_n {\rm e}^{t \lambda_n}{\rm e}_n(x){\rm e}_n(0).
\end{equation}
Indeed, according to the phenomenon of Anderson localization, the leading eigenfunctions ${\rm e}_1,{\rm e}_2,\dots$ are supposed to be concentrated in such islands.  A proper proof of this concentration and the one of $u(t,\cdot)$ in the same islands has been given in terms of a kind of spatial extreme-value picture in large centred boxes in some few cases, for example for $\xi$ an i.i.d.\ potential on $\Z^d$ with double-exponential tails by Biskup, K\"onig and dos Santos \cite{BiKodS18}. See K\"onig \cite{Ko16} and Astrauskas \cite{As16} for two extensive surveys of the mathematical treatment of the PAM until 2016.

\subsection{Our main purpose and strategy}\label{sec-purpose}

The present paper is a contribution to the development of 
methods for proper formulation and proofs of intermittency for the PAM in the case of a Gaussian white-noise potential in two dimensions. This problem makes sense for dimensions $d=1,2,3$. For $d\ge 4$ the PAM with white-noise potential is scaling-critical respectively scaling-supercritical in the sense of Hairer~\cite{Ha14} and there is no known solution theory.  
The cases $d=2,3$ are conceptually similar to each other and most arguments developed for dimension $2$ are expected to extend to dimension $3$, but to require more technicalities. On the other hand, the one-dimensional case is simpler for several reasons, and there already exist very good localization results for the leading eigenfunctions of the Anderson Hamiltonian by Dumaz and Labb\'e~\cite{DuLa17}.

The main purpose of the present paper is to identify the {\em almost-sure large-$t$ asymptotics} of the total mass, $U(t)$, of the solution to \eqref{PAM}.
Based on earlier treatment of this question for more regular random fields, it is easy to guess what the answer should be. 
Indeed, 
in case $\xi$ is a more regular potential, then the total mass of the solution to the PAM admits a Feynman--Kac representation, informally written as 
\begin{equation}\label{FKform}
U(t)=
 \E^{0}\big[\e^{\int_0^t \xi(B_s)\,{\rm d}s}  \big], 
\end{equation}
for a standard Brownian motion $B$ in $\R^2$ starting from $0$.  
From this one reads that it is probabilistically very costly for the Brownian motion starting from the origin,  to reach a remote site of distance  $t$  
from the starting site. 
This should imply the asymptotic behaviour 
 $\log U(t)\sim \log U_{t}(t)$,  
where  $U_{t}(t)$  
is the total mass of the solution 
 $u_{t}(t,\cdot)$  
to \eqref{PAM} in the centred box 
 $Q_{t}= [ \frac{t}{2}, \frac{t}{2}]^2$  of diameter  $t$  with Dirichlet boundary condition. Using the expansion  \eqref{eigenvalueexpansion} in that box, we guess that 
 $\log U_{t}(t)\sim t \blambda_1(Q_{t}) $, 
where  $\blambda_1(Q_{t})$  is the corresponding principal eigenvalue of the Anderson Hamiltonian $\frac 12\Delta +\xi$ (after renormalization): 
See Allez and Chouk \cite{AlCh15} and also Gubinelli, Ugurcan and Zachhuber~\cite{GuUgZa20} and Labb\'e~\cite{La19} for the derivation of the spectrum  and Rayleigh-Ritz representations for the eigenvalues of the Anderson Hamiltonian on boxes. 
The asymptotics of the eigenvalues as the size of the box diverges have been described by Chouk and van Zuijlen~\cite{ChvZ21}, who show that almost surely  $\blambda_1(Q_{t})\sim \chi \log t $ as $t \to \infty$, where $\chi > 0$ is a deterministic constant given by the variation formula \eqref{eqn:chi_var_expr}. 
Summarizing, we expect that the large-$t$  asymptotics of $U(t)$ are given as $\log U(t)\sim \chi t \log t$. 
Such a line of arguments has been made rigorous for various types of regular or i.i.d.\ random potentials in $\R^d$ and in $\Z^d$ as the ambient space; see  \cite{Ko16} and \cite{As16}. 

However, in the present case of $u$ being a renormalized solution for a Gaussian white-noise potential, it presents a formidable task to carry through this programme. 
The difficulty lies in the fact that the Feynman--Kac representation~\eqref{FKform}  is only formal notation and meaningless. We could interpret it rigorously with the help of the random polymer measure in the work by Cannizzaro and Chouk~\cite{CaCh18}, but that does not seem very helpful because the construction of the polymer measure actually involves the total mass. Also,  under the polymer measure the coordinate process roughly speaking solves the SDE $\dd X_s = \nabla \log u(t-s,X_s) \dd s + \dd B_s$, so to analyze $X$ we already need information about $u$.  Hence, we cannot rely on the results of \cite{CaCh18} and will therefore go via another route; we explain this now. 

The main difficulty to justify the above heuristic steps is to prove the assertion  $\log U(t)\sim\log U_{t}(t)$. 
To derive a proof we use a slight modification of the following procedure, which is often used for sufficiently well-behaved random potential $\xi$. We decompose the Feynman--Kac representation for $U(t)-U_{t}(t)$ into the sum over $k\in\N$ of the contributions to expectation of $\e^{\int_0^t \xi(B_s)\,{\rm d}s}$ coming from Brownian paths $(B_s)_{s\in[0,t]}$ that leave the box $Q_{t^k}$, but not the box $Q_{t^{k+1}}$. The $k$-th contribution then is estimated in terms of the maximum of the potential $\xi$ within $Q_{t^{k+1}}$ (using extreme-value analysis, say) times the probability that the Brownian motion leaves $Q_{ t^k}$ (using the reflection principle, say). 
In order to implement this strategy in our setting, we use a ``partial Girsanov transformation'', an idea  that was introduced in the setting of the KPZ equation by Gubinelli and Perkowski~\cite{GuPe17}. (We drop the mollification and renormalization and the dependence on the box from the notation.) Denote by  $Z=(1-\frac 12\Delta)^{-1}\xi$ the resolvent of the random potential and by $Y$ the solution to 
\begin{align*}
(\eta-\tfrac 12\Delta)Y=\tfrac 12 |\nabla Z|^2+\nabla Y\cdot \nabla Z,
\end{align*}
which exists and is unique for sufficiently large $\eta>0$. 
Now put $b=\nabla( Z+  Y)$ and consider the solution $X$ to the SDE ${\rm d} X_t= b (X_t)\,{\rm d}t+{\rm d} B_t$ for some Brownian motion $B$. Then it turns out that
\begin{align*}
U(t)=\E\Big[\e^{\int_0^t (Z+\eta-Y+\frac 12|\nabla Y|^2)(X_s)\dd s}\,\e^{(Y+Z)(X_0)-(Y+Z)(X_t)}\Big].
\end{align*}
While this \emph{Feynman--Kac-type representation} looks more complicated than~\eqref{FKform}, it has the advantage that it also works for the white-noise potential. 
To carry out the general proof strategy outlined above we then need to understand how $Y,Z$ grow with the box size, and we need bounds on the probability of $X$ leaving  the box $Q_{t^k}$. Since we do not need very precise bounds on $Y,Z$, the first part is quite easy and essentially follows from Gaussian hypercontractivity. 
From this we obtain that both $Y$ and $Z$ grow at most logarithmically with respect to the box size. 
The estimate on the probability of $X$ leaving the box $Q_{t^k}$ is shown in \cite[Corollary 1.2]{PevZ} and is a consequence of heat kernel estimates for the transition kernel of the diffusion $X$.


\subsection{Main results}\label{sec-results}

Let us precisely formulate our main results. For a smooth potential $V$, we denote by $\lambda_n(Q_t,V)$ the $n$-th largest (counted with multiplicity) Dirichlet eigenvalue of the operator $\frac12 \Delta+V$ in a box $Q_t = [-\frac{t}{2},\frac{t}{2}]^2$ and we write $\lambda(Q_t,V)=\lambda_1(Q_t,V)$ for the principal one. 
We introduce the following variational formula: 
\begin{equation}
\label{eqn:chi_var_expr}
\chi=8
\sup_{ \substack{ \psi \in C_{\rm c}^\infty(\R^2) \colon \\ \|\psi \|_{L^2} =1} } 
 \|\psi\|_{L^4}^2 - \int | \nabla \psi|^2. 
\end{equation}
$\chi$ is finite and $\frac12 \chi$ equals the smallest $C>0$ such that $\|f\|_{L^4}^4 \le C \|\nabla f\|_{L^2}^2 \|f\|_{L^2}^2$ for all $f\in H^1(\R^2)$ (this is Ladyzhenskaya's inequality, which is a special case of the Gagliardo-Nirenberg inequality), see \cite[Theorem 2.6]{ChvZ21}. 

The following theorem is our main result. 
For $A\subset \R$ we write 
\begin{align*}
\mbox{``}a_t \sim b_t , \quad t \in A, t\rightarrow \infty\mbox{''}
\quad \mbox{ instead of } \quad
\mbox{``}\lim_{\substack{ t\in A \\ t\rightarrow \infty}} \frac{a_t}{b_t} = 1\mbox{''}.
\end{align*}

\begin{theorem} 
\label{theorem:total_mass_asymptotics}
\begin{enumerate}
\item 
\label{item:total_mass_asymptotics}
\textnormal{[Asymptotics of the total mass]}
Almost surely, 
\begin{align*}
\log U(t) \sim \chi t \log t \qquad t\in \Q, t\to\infty. 
\end{align*}
\item 
\label{item:infimum_and_supremum_asymptotics}
\textnormal{[Asymptotics of the min and max]} For all $a\in (0,1)$, almost surely
\begin{align*} 
\log \left( \min_{x \in Q_{t^a}} u (t, x) \right) \sim \log \left(\sup_{x \in \mathbb{R}^2} u (t, x) \right) \sim \chi t \log t \qquad t\in \Q, t\rightarrow \infty. 
\end{align*}
\end{enumerate}
\end{theorem}

Theorem~\ref{theorem:total_mass_asymptotics}~\ref{item:infimum_and_supremum_asymptotics} means that on the level of the logarithmic asymptotics we see no intermittency effect, the $L^{\infty}$-norm of the solution to the PAM almost surely has the same logarithmic asymptotics as the $L^1$-norm.

Theorem~\ref{theorem:total_mass_asymptotics}~\ref{item:total_mass_asymptotics} is an immediate consequence of the following three main results, Proposition \ref{proposition:comparison_total_mass_full_space_and_box}, \ref{theorem:convergence_smooth_eigenvalues_and_asymptotics} and  \ref{theorem:comparison_box_with_t_times_first_eigenvalue}.

The first main step is a \lq compactification\rq, a reduction of the solution to some ($t$-dependent) box with Dirichlet boundary condition.   We write $U_L$ for the total mass of the solution of the parabolic Anderson model on the box $Q_L$ with Dirichlet boundary conditions.

\begin{proposition}[Reduction to a box]
\label{proposition:comparison_total_mass_full_space_and_box}
Almost surely, 
\begin{align*}
\log U(t) \sim \log U_{t} (t) \qquad \mbox{ as } t\rightarrow \infty. 
\end{align*}
\end{proposition}

The following is a strengthening of the main result of \cite{ChvZ21}, in that we develop and use a renormalization technique here that does not depend on the box considered. Recall our notation of eigenvalues from the beginning of this section. 

\begin{theorem}[Eigenvalue asymptotics: renormalization and large boxes]
\label{theorem:convergence_smooth_eigenvalues_and_asymptotics}
Let $\xi_\epsilon = \psi_\epsilon * \xi$ be a mollification of the white noise for a mollifying function $\psi$ (see Section~\ref{section:representation_on_boxes}) for $\epsilon>0$. 
There exists a $C$ in $\R$ that only depends on $\psi$ such that with $c_\epsilon := \frac{1}{2\pi} \log \frac{1}{\epsilon}+C$, we have for all $L\in [1,\infty)$ and $n\in \N$,  in probability 
\begin{align*}
\lim_{\epsilon\downarrow 0} \lambda_n(Q_L,\xi_\epsilon) - c_\epsilon = 
 \blambda_{n}(Q_L), 
\end{align*}
where  $\blambda_{n}(Q_L)$  is the $n$-th largest eigenvalue (counting multiplicities) of the Anderson Hamiltonian defined by the enhanced white-noise potential (see \eqref{obs:notation_H_L_etc}).  
Moreover, for all $n\in\N$ and for any countable unbounded set $\I \subset (e,\infty)$,  almost surely, 
\begin{align*}
\lim_{L\in \I, L\rightarrow \infty}
\frac{ \blambda_{n}(Q_L)}{\log L} = \chi. 
\end{align*}
\end{theorem} 

The last step is to use the eigenvalue expansion in \eqref{eigenvalueexpansion} in a suitable way to derive that the large-$t$ exponential rate of the total mass in a large box is asymptotically equivalent to the principal eigenvalue:

\begin{theorem}
\label{theorem:comparison_box_with_t_times_first_eigenvalue}
Almost surely 
\begin{align*}
\tfrac 1t\log U_{t}(t) \sim \blambda_1(Q_{t})
\qquad  t\in \Q , t\rightarrow \infty. 
\end{align*}
\end{theorem}

Combining Proposition~\ref{proposition:comparison_total_mass_full_space_and_box} with the second assertion of Theorem~\ref{theorem:convergence_smooth_eigenvalues_and_asymptotics} and Theorem~\ref{theorem:comparison_box_with_t_times_first_eigenvalue} implies Theorem~\ref{theorem:total_mass_asymptotics}.

\subsection{Remarks on more regular Gaussian potentials}\label{sec-discussion}

It is interesting to note that an interchange of the two limits $t\to\infty$ and $\eps\to 0$ (the mollification parameter) does not even phenemonologically or heuristically  explain our main result; the asymptotics of Theorem~\ref{theorem:total_mass_asymptotics} lie much deeper. Indeed, if we use for example the mollification $\xi_\epsilon=p_\epsilon * \xi$ (with $p_\epsilon$ the Gaussian kernel with variance $\epsilon$), then we obtain a smooth centered Gaussian field with covariance function $p_{2\epsilon}$, using the convolution property of the Gaussian kernel. The almost-sure large-$t$ asymptotics of the total mass $U_\epsilon(t)$ of the solution $u_\epsilon(t,\cdot)$ of \eqref{PAM} with potential $\xi_\epsilon-c_\epsilon$ are given in G\"artner, K\"onig and Molchanov~\cite{GaKoMo00} by
\begin{equation}\label{smoothGaussasy}
\frac 1t \log U_\epsilon(t)=-c_\epsilon+\frac1{\sqrt{\pi\epsilon}}(\log t)^{\frac 12}-\frac1{\sqrt{2\pi}\epsilon}(\log t)^{\frac 14}(1+o(1)),\qquad t\to\infty.
\end{equation}
While the first term, $c_\epsilon$, comes simply from subtracting it from the potential $\xi_\epsilon$, the other two terms are the asymptotics of $\lambda(Q_{L_t},\xi_\epsilon)$ (we used that $\log L_t\sim\log t$). Indeed, they come from a second-order extreme-value analysis of the Gaussian potential in $Q_{L_t}$: its absolute height is of order $(\log t)^{\frac 12}$, due to Gaussian tails, and the last term comes from its geometry in the local peak in an island of diameter $\asymp (\log t)^{-\frac 14}$ in which the potential approaches an explicit parabola. In contrast, the asymptotics of the principal eigenvalue of the limiting (enhanced) potential, $\blambda(Q_{L_t})$, are actually of order $\log t$, due to {\em exponential} tails of $\blambda(Q_{L_t})$ instead of the Gaussian tails of $\lambda(Q_{L_t},\xi_\epsilon)$, see \cite[Theorem 2.7]{ChvZ21}.

A closely related work to ours is Chen~\cite{Ch14}, where the Gaussian potential $\xi$ in $\R^d$ was assumed to (be centred and stationary and) have a covariance matrix $B$ that behaves like $B(x)=\sigma^2|x|^{-\alpha}$ as $x\to 0$ with $\sigma^2\in(0,\infty)$ and $\alpha\in(0,\min\{2,d\})$. Note that $\xi$ is not a function, although it is still more regular than the white noise in $d=2$. Then the almost-sure asymptotics of the total mass are identified as 
\begin{equation}\label{Chenasy}
\lim_{t\rightarrow \infty} 
 \frac{\log U(t)}{t (\log t)^{\frac{2}{4-\alpha}}} 
  = (2d\sigma^2 )^{\frac {2}{4-\alpha}} \chi_\alpha ,
\end{equation}
where (see \cite[Lemma A.4]{Ch14})
\begin{align*}
\chi_\alpha=\sup_{g\in H^1(\R^d)\colon\|g\|_{L^2}=1}\Big[\Big(\int_{\R^d}\int_{\R^d}\frac{g^2(x)g^2(y)}{|x-y|^\alpha} \dd x \dd y\Big)^{\frac 12}-\frac12 \|\nabla g\|_{L^2}^2\Big].
\end{align*}
Like in Theorem~\ref{theorem:total_mass_asymptotics}, the asymptotic behavior is given by an interesting variational formula, and the powers of the logarithm coincide for the boundary case $\alpha=2$ in $d=2$, although in that case $\chi_2 = \infty$ (contrary to our situation where $\chi$ is finite).

In the recent work by Lamarre \cite{GL}, a $t$-dependent Gaussian regularization of the eigenvalues in the box $(-t,t)^d$ and of the total mass of the PAM at time $t$ is asymptotically analysed in dimensions $d\in\{1,2,3\}$ for $t\to\infty$. A critical scale for the size of the  vanishing regularization parameter is identified: If it vanishes less fast than this scale, then the first-order asymptotics are identical to the ones for a fixed regular Gaussian field, and if it vanishes faster than this scale, then this asymptotics is the same as the ones for the unregularized Gaussian white noise in $d\in\{1,2\}$ (proved in \cite{DuLa17} in $d=1$, in \cite{ChvZ21} and the present paper, respectively, in $d=2$), with an analogous result in $d=3$. The techniques of \cite{GL} do not allow for the investigation of the unregularized white-noise case, which we consider in the present paper. Our main result appears there as \cite[Conjecture 1.19]{GL}.

\subsection{Outline of the paper}

The rest of the paper is built up as follows. 
In Section~\ref{section:representation_on_boxes} we introduce fundamental notation and auxiliary material involving
representations of the solution of the (Cauchy problem for the) heat equation with general potential on a box, in particular both the spectral representation and a Feynman--Kac-type representation.
In Section~\ref{section:asymptotics_large_box} we prove Theorem 1.4, that is, we derive the asymptotic behavior of the total mass of the solution with white-noise potential
on a large box. 
In Section~\ref{section:separating_space_in_boxes} we prove Proposition~\ref{proposition:comparison_total_mass_full_space_and_box}, 
that is, we show
that the logarithm of the total mass at time $t$ is asymptotically the same as the one of the restriction to a box with diameter $t$.
In Section~\ref{section:asymptotics_of_sup_and_inf_boxes} we prove Theorem~\ref{theorem:total_mass_asymptotics}~\ref{item:infimum_and_supremum_asymptotics}. 
Section~\ref{section:gaussian_calcs} is 
dedicated to proving that
the enhanced noise that we obtain here via mollification on the entire $\R^d$, after projection to
a box, is the same as in the setting of \cite{ChvZ21} where a box-dependent approximation is carried out.

\subsection{Notation}

We write $\N= \{1,2,\dots\}, \N_0= \{0\}\cup \N$ and $\N_{-1} = \{-1\} \cup \N_0$. 
For families $(a_i)_{i \in \I}, (b_i)_{i\in \I}$ in $\R$ for an index set $\I$, we write $a_i \lesssim b_i$ to denote the existence of a $C>0$ such that $a_i \le C b_i$ for all $i\in \I$.

\section{Two representations for the heat equation with potential on a box}
\label{section:representation_on_boxes}

In this section we present two different representations of solutions to the heat equation with a (deterministic) irregular potential on a box with Dirichlet boundary conditions. 
First, in Section~\ref{subsection:spectral_representation}, we use the Anderson Hamiltonian with Dirichlet boundary conditions for irregular potentials constructed in \cite{ChvZ21} to obtain a spectral representation. 
Then we present a new Feynman--Kac-type representation, as explained in Section~\ref{sec-purpose}. 
Due to the irregularity of the potential the classical Feynman--Kac representation does not make sense. 
We are however able to adapt to irregular potentials using a trick which is inspired by the `partial Girsanov transform' of \cite{GuPe17}. 
We do this for a deterministic potential in Section~\ref{subsection:feynman_kac_type} and apply it to the white noise potential in Section~\ref{subsection:feynman_kac_for_pam_with_white_noies}.

First we introduce the heat equation with smooth potentials and for enhanced potentials, see \ref{obs:solution_smooth_pam_full_space_and_box} and \ref{obs:multiplicative_heat_eq_with_enhanced_potentials}. 

\bigskip

\textbf{In this section, for $\epsilon>0$, $\theta_\epsilon$ is a smooth function on $\R^2$ and $c_\epsilon \in \R$.}
More assumptions will be made in \ref{obs:convergence_assumptions_theta} and in \ref{obs:assumption_alpha_and_L}.

\begin{obs}[Heat equation with smooth potential]
\label{obs:solution_smooth_pam_full_space_and_box}
We write $Q_L = [-\frac{L}{2},\frac{L}{2}]^2$ for $L\in (0,\infty)$ and $Q_\infty = \R^2$. 
For $L\in (0,\infty]$ and $\epsilon>0$ 
we let $u_{L,\epsilon}^\phi$ be the solution to the heat equation with potential $\theta_\epsilon - c_\epsilon$ and Dirichlet boundary conditions on the box $Q_L $ with initial condition $\phi$:
\begin{align}
\label{eqn:pam_box_smooth}
\begin{cases}
\partial_t u_{L,\epsilon}^\phi = \tfrac12 \Delta u_{L,\epsilon}^\phi + (\theta_\epsilon - c_\epsilon  ) u_{L,\epsilon}^\phi \qquad \mbox{ on } (0,\infty) \times Q_L, \\
u_{L,\epsilon}^\phi(0,\cdot) = \phi, \qquad \mbox{ and } \qquad 
u_{L,\epsilon}^\phi|_{\partial Q_L}(t,\cdot) = 0.
\end{cases}
\end{align}
($u_{\infty,\epsilon}^\phi |_{ \partial Q_\infty}=0$ is interpreted as an empty condition.) 
\end{obs}

 In the case that $\theta_\epsilon$ represents mollified white noise and $c_\epsilon = \frac{1}{\pi} \log \frac{1}{\epsilon}$,  
we know that $u_{\infty,\epsilon}^\phi$ converges locally uniformly in probability as $\epsilon\to 0$ to $u_\infty^\phi$, the solution to the continuous heat equation with renormalised potential and initial condition $\phi$; see~\cite[Theorem 4.1]{HaLa15}. 
To derive an analogous result for $u_{L,\epsilon}^\phi$, we use the approach to paracontrolled distributions with Dirichlet boundary conditions from \cite{ChvZ21}. 

\begin{obs}[Notation]
\label{obs:notation_neumann}
For a function $f \colon [0,L]^2 \rightarrow \R$ we write $\overline f :\R^2 \rightarrow \R$ for the even extension of $f$ with period $2L$. 
This means $\overline f( \fq_1 x_1, \fq_2 x_2) = f(x_1,x_2)$ for $\fq_1,\fq_2\in \{-1,1\}$ and $x_1,x_2\in [0,L]$ and $\overline f(x) = \overline f(x+2Lk)$ for $x\in [-L,L]^2$ and $k\in \Z^2$. 

If $g$ is a function on $Q_L$, which is a translation of the set $[0,L]^2$ by  $y = (\frac{L}{2},\frac{L}{2})$, we write $\overline g$ for $\overline{g(\cdot-y)}(\cdot+y)$ and also call this the even extension of $g$. 

The following notation we borrow from \cite[Section~4]{ChvZ21}, see that reference for more details. 
We write $\cS_\fn(Q_L)$ for the space of all $f\colon Q_L \to \R$ such that $\overline f \in C^\infty(\R^2)$. 
Let $B_{p,q}^{\fn,\alpha}(Q_L)$ be the Neumann Besov space defined as in  \cite[Section~4]{ChvZ21}. We abbreviate $\cC^\alpha_\fn = B_{\infty,\infty}^{\fn,\alpha}$. We write $(\fn_{k,L})_{k\in\N_0^2}$ for the Neumann basis of $L^2(Q_L)$ and, slightly abusing notation, we also write ``$\fn_{k,L}$'' for the extension of $\fn_{k,L}$ to $\R^2$ that equals $0$ outside $Q_L$. 
We write $\overline \fn_{k,L}$ for the even extension of $\fn_{k,L}$. 
For an even function $\sigma: \R^2 \rightarrow \R$ and $f\in \cS_\fn'(Q_L)$ we write $\sigma(\rD) f$ for the Fourier multiplier
\begin{align*}
\sigma(\rD) f = \sum_{k\in \N_0^2} \sigma(\tfrac{k}{L}) \langle f , \fn_{k,L} \rangle \fn_{k,L}. 
\end{align*} 
\end{obs}

\begin{obs}[Notation]
\label{obs:notation_para_and_reso}
We write $u \para v$  for the paraproduct between $u$ and $v$ (with the low frequencies of $u$ and the high frequencies of $v$), and $u\reso v$ for the resonance product; we adopt the notation from \cite{MaPe19} and refer to \cite{BaChDa11} for background material.  
\end{obs}

\begin{definition}
\label{def:enhanced_neumann_potentials}
Let $L>0$ and $\sigma(x) = (1+\frac12\pi^2|x|^2)^{-1}$.  
For $\beta \in \R$, we define the 
space of \emph{enhanced Neumann potentials}, written ``$\fX_\fn^\beta(Q_L)$'' or just ``$\fX_\fn^\beta$'',  as the closure in $\cC_\fn^\beta \times \cC_\fn^{2\beta+2}$ of the set $\{ (\zeta, \zeta \reso \sigma(\rD) \zeta - c) : \zeta \in \cS_\fn, c\in \R\}$, equipped with the relative topology with respect to $\cC_\fn^\beta \times \cC_\fn^{2\beta +2}$. 
\end{definition}

Let 
\begin{align}
\label{eqn:theta_L_eps_fn_projection}
\theta_{L,\epsilon} = \sum_{k\in\N_0^2} \langle \theta_\epsilon, \fn_{k,L} \rangle \fn_{k,L}, 
\qquad 
\mbox{ so that } \quad 
\overline{\theta_{L,\epsilon}} = \sum_{k\in\N_0^2} \langle \theta_\epsilon, \fn_{k,L} \rangle \overline \fn_{k,L}.
\end{align} 
Since $\theta_\epsilon$ is almost everywhere equal to $\theta_{L,\epsilon}$ and thus to $\overline{\theta_{L,\epsilon}}$ on $Q_L$, the solution $u_{L,\epsilon}^\phi$ to~\eqref{eqn:pam_box_smooth} also solves the same equation with  ``$\theta_\epsilon $'' replaced by ``$\theta_{L,\epsilon}$'' or ``$\overline{\theta_{L,\epsilon}}$''.

\begin{obs}
\label{obs:convergence_assumptions_theta}
\textbf{For the rest of this section we assume that for all $L>0$ there exists a $\btheta_L  \in \bigcap_{\alpha<-1} \fX_\fn^\alpha$ and $c_\epsilon\in \R$ (not depending on $L$) for $\epsilon>0$ such that the following convergence holds in $\fX_\fn^\alpha$ for all $\alpha <-1$} 
\begin{align}
\label{eqn:limit_of_enhanced_theta_L_eps_fn}
(\theta_{L,\epsilon}, \theta_{L,\epsilon} \reso \sigma(\rD) \theta_{L,\epsilon} - c_\epsilon ) \arroweps \btheta_L. 
\end{align}
\end{obs}

\begin{obs}[Heat equation with enhanced potentials]
\label{obs:multiplicative_heat_eq_with_enhanced_potentials}
We write $\cC_p^{\fd, \beta}(Q_L)= B_{p,\infty}^{\fd,\beta}(Q_L)$, where the latter is the Dirichlet Besov space defined as in \cite[Section 4]{ChvZ21}, and $\cC^{\fd,\alpha}(Q_L) = \cC^{\fd,\alpha}_\infty(Q_L)$. By taking an odd extension of $u_{L,\epsilon}^\phi$ and an even extension of the noise (as in \cite{ChvZ21}) we obtain a periodic solution on the torus of length $2L$. 
Therefore the convergence shown in \cite[Theorem~5.4]{GuImPe15} implies the following. 
For $\beta \in (\frac23,1)$ and for all $T>0$ and $\phi \in \cC^{\fd, \beta}(Q_L)$ the solution $u_{L,\epsilon}^\phi$ converges (in $C([0,T], \cC^{\fd, \beta}(Q_L))$ and therefore) uniformly on $[0,T] \times Q_L$ in probability to $u_L^\phi$, where $u_L^\phi$ solves 
\begin{align}
\label{eqn:pam_xi_dir_bc}
\begin{cases}
\partial_t u_L^\phi = \tfrac12 \Delta u_L^\phi + \btheta_L \diamond u_L^\phi \qquad \mbox{ on } (0,\infty) \times Q_L, \\
 u_L^\phi(0,\cdot) = \phi, \qquad \mbox{ and } \qquad u_{L,\epsilon}^\phi|_{\partial Q_L} = 0,
\end{cases}
\end{align}
interpreted in the sense of paracontrolled distributions (i.e., as in \cite[Section~5]{GuImPe15}). 
Let us point out that we have some differences, namely, we have a factor $\frac12$ in front of the Laplacian, which is also the reason why our renormalisation constant differs by a factor $2$ (we explain this in the proof of Theorem~\ref{theorem:results_ChvZ21}). 
\end{obs}

Combining arguments as in \cite[Proposition 2.4]{PevZ} with arguments as in \cite[Section~6]{GuPe17} in the context of singular initial conditions, 
 we can extend this result to more general initial conditions, and also we get the continuous dependence of the solution on the initial condition.

\begin{obs}[Notation]
\label{obs:notation_solutions}
In the next lemma we write ``$u_{L,0}^\phi$'' for ``$u_L^\phi$'', where $u_L^\phi$ is as in \eqref{eqn:pam_xi_dir_bc} and we write $u_{L,\epsilon}^\phi$ for the solution of \eqref{eqn:pam_box_smooth}. 
\end{obs}

\begin{lemma}
\label{lemma:solution_map_box_continuous_in_inintial}
Let $\alpha \in (-\frac{4}{3},-1)$ and $\gamma \in (-1, 2 + \alpha)$, $p \in [1,\infty]$ and $L>0$. 
The solution map $[0,\infty) \times \cC_p^{\fd,\gamma}(Q_L) \to C([0,\infty), \cC_p^{\fd,\gamma}(Q_L))$ given by $(\epsilon,\phi) \mapsto u_{L,\epsilon}^\phi$  is continuous, where $u_{L,0}^\phi := u_L^\phi$.  Moreover, for all $\beta < 2 + \alpha$, $t>0$, $L< \infty$ and $\epsilon \in [0,\infty)$  we have $u_{L,\epsilon}(t,\cdot) \in \cC^{\fd,\beta}$ and the map $[0,\infty) \times \cC_p^{\fd,\gamma}(Q_L) \rightarrow  \cC^{\fd,\beta}(Q_L)$ given by  $(\epsilon,\phi) \mapsto u_{L,\epsilon}^\phi(t,\cdot)$ is continuous. 
\end{lemma}
\begin{proof}
The continuity is shown similarly to the continuity statement in \cite[Proposition 2.4]{PevZ}. Also the better integrability, i.e., that $u_{L,\epsilon}(t,\cdot) \in \cC^\beta$ and not just $u_{L,\epsilon}(t,\cdot) \in \cC^\beta_p$, can be shown by a bootstrap argument similar to \cite[Proposition 2.4]{PevZ}.
\end{proof}

\subsection{The spectral representation}
\label{subsection:spectral_representation}

Here we give a representation of $u_{L}$ in terms of the eigenfunctions of the Anderson Hamiltonian from~\cite{ChvZ21}. Later we will use this representation to study the logarithmic asymptotics of the total mass  of the PAM on the box $Q_L$.

\begin{theorem}
\label{theorem:dirichlet_summary}
\textnormal{\cite[Theorem 5.4]{ChvZ21}}
Let $\alpha \in (-\frac{4}{3},-1)$, $L>0$ and $\btheta \in \fX_\fn^\alpha(Q_L)$. 
There exists a domain $\fD_\btheta^{\fd} \subset L^2(Q_L)$ such that $\sH_\btheta$ defined on $ \fD_\btheta^{\fd} $ by 
\begin{align*}
\sH_\btheta u = \tfrac12 \Delta u + \btheta \diamond u 
\end{align*}
 is a closed and self-adjoint operator with values in $L^2(Q_L)$. 
 $\sH_\btheta$ has a pure point spectrum $\sigma(\sH_\btheta)$ consisting of eigenvalues
$\lambda_1(Q_L,\btheta) > \lambda_2(Q_L,\btheta) \ge  \lambda_3(Q_L,\btheta)  \ge \cdots $ 
(counting multiplicities) 
 with $\lim_{n\to \infty} \lambda_n(Q_L, \btheta) = -\infty$, 
and such that
\begin{align*}
\fD_\btheta^\fd = \bigoplus_{\lambda \in \sigma(\sH_\btheta)} \ker ( \lambda - \sH_\btheta). 
\end{align*} 
Moreover, the map $\btheta \mapsto \lambda_n(Q_L,\btheta)$ is (locally Lipschitz) continuous and $\fD_\btheta^\fd$ is dense in $L^2(Q_L)$.
\end{theorem}

\begin{obs}
\label{obs:notation_H_L_etc}
For abbreviation, we write ``$\sH_L$'' and ``$\blambda_{n,L}$'' for ``$\sH_{\btheta_L}$'' and ``$\lambda_n(Q_L,\btheta_L)$''. 
We let $(v_{n,L})_{n \in \N}$ be an orthonormal basis of $L^2(Q_L)$, such that $v_{n,L}$ is an eigenvector with eigenvalue $\blambda_{n,L}$. 
\end{obs}


By the properties of the semigroup of the Anderson Hamiltonian, we obtain the following representation for the solution to \eqref{eqn:pam_xi_dir_bc} with initial condition $\phi \in L^2(Q_L)$. 

\begin{lemma}
\label{lemma:spectral_rep_sol_box_smooth_initial}
For $L,t >0$, $y\in Q_L$  and $\phi \in L^2(Q_L)$ 
\begin{align}
\label{eqn:solution_for_L_2_initial_on_box}
u_L^\phi (t,y) = [e^{t\sH_L} \phi](y) = \sum_{n\in\N} e^{t\blambda_{n,L}} \langle v_{n,L} ,\phi \rangle_{L^2} v_{n,L}(y). 
\end{align}
\end{lemma}

\begin{proof}
For $\phi \in \fD_{\btheta_L}$ the identity holds since the semigroup $S(t)= e^{t\sH_L}$ satisfies $\partial_t S(t) \phi = \sH_L S(t) \phi$, see \cite[Theorem 2.4]{Pa83}. 
As $\fD_{\btheta_L}$ is dense in $L^2(Q_L)$ and convergence in $L^2(Q_L)$ implies convergence in $\cC_2^{\fd,-\beta}(Q_L)$ for all $\beta>0$ (combine for example \cite[Theorem 2.71]{BaChDa11} with \cite[Theorem 4.7]{ChvZ21}), Lemma~\ref{lemma:solution_map_box_continuous_in_inintial} shows that \eqref{eqn:solution_for_L_2_initial_on_box} holds for all $\phi \in L^2(Q_L)$.
\end{proof}

In particular, Lemma~\ref{lemma:spectral_rep_sol_box_smooth_initial} yields $v_{n,L} = e^{-t\blambda_{n,L}} u_L^{v_{n,L}} (t,\cdot)\in \cC^{\fd,\beta}(Q_L)$ for all $\beta < 2+\alpha$.

\begin{theorem}
\label{theorem:spectral_rep_sol_box_delta_initial}
Let $L>0$. 
For $x,y\in Q_L$ and $\phi =\delta_x$ the solution to \eqref{eqn:pam_xi_dir_bc} is given by 
\begin{align}
\label{eqn:eigenvalue_rep_on_box}
u_L^{\delta_x} (t,y)= [e^{t\sH_L} \delta_x](y) =  \sum_{n\in\N} e^{t\blambda_{n,L} } v_{n,L}(x) v_{n,L}(y).
\end{align}
\end{theorem}
\begin{proof}
Let $\psi \in C_c^\infty$ be such that $\int \psi =1$ and $\psi_\epsilon^x(y) := \psi(\frac{y-x}{\epsilon})$. 
Then $\psi_\epsilon^0 * f \rightarrow f$ in $L^p$ for $f\in L^p$ (see \cite[Proposition 1.2.32]{HyvNVeWe16}). 
Therefore $\psi_\epsilon^x \rightarrow \delta_x$ in $B_{1,\infty}^{\gamma}$ for all $\gamma < 0$. 
\begin{calc}
Because $\Delta_i \psi_\epsilon^0 = \cF^{-1}(\rho_i) * \psi_\epsilon^0 \rightarrow \cF^{-1}(\rho_i) = \Delta_i \delta_0$ in $L^1$ for all $i\in \N_{-1}$. 
\end{calc}
We have for all $\epsilon > 0$:
	\begin{equation}\label{eqn:spectral_rep_sol_box_delta_initial-pr}
	\begin{aligned}
		\left| u_L^{\delta_x} (t,y) -  \sum_{n\in\N} e^{t\blambda_{n,L} } v_{n,L}(x) v_{n,L}(y) \right| & \le \left| u_L^{\delta_x} (t,y) -  \sum_{n\in\N} e^{t\blambda_{n,L} } \langle v_{n,L}, \psi_\epsilon^x\rangle \langle v_{n,L}, \psi^y_\epsilon\rangle \right| \\ 
		& \quad + \sum_{n\in\N} e^{t\blambda_{n,L} } |  \langle v_{n,L}, \psi_\epsilon^x - \delta_x \rangle \langle v_{n,L}, \psi^y_\epsilon\rangle | \\
		& \quad + \sum_{n\in\N} e^{t\blambda_{n,L} } |  \langle v_{n,L}, \delta_x \rangle \langle v_{n,L}, \psi^y_\epsilon - \delta_y \rangle |.
	\end{aligned}	
	\end{equation}
	The first term on the right-hand side equals $| u_L^{\delta_x} (t,y)-  \langle u_L^{\psi^x_\epsilon} (t,\cdot), \psi^y_\epsilon\rangle|$, and since $u_L^{\psi^x_\epsilon} (t,\cdot)$ converges uniformly to $u_L^{\delta_x} (t,\cdot)$ by Lemma~\ref{lemma:solution_map_box_continuous_in_inintial}, it vanishes as $\epsilon \downarrow 0$. The second and third term on the right-hand side of~\eqref{eqn:spectral_rep_sol_box_delta_initial-pr} are similar, so we only argue for the last term. We saw that $v_{L,n} \in \cC^{\fn,\beta}(Q_L)$, so in particular $v_{L,n}$ is continuous and we have by Fatou's lemma:
	\begin{align*}
		\sum_{n\in\N} e^{t\blambda_{n,L} } |  \langle v_{n,L}, \delta_x \rangle \langle v_{n,L}, \psi^y_\epsilon - \delta_y \rangle | & \le \liminf_{\epsilon'\to 0} \sum_{n\in\N} e^{t\blambda_{n,L} } |  \langle v_{n,L}, \psi^x_{\epsilon'} \rangle \langle v_{n,L}, \psi^y_\epsilon - \psi^y_{\epsilon'} \rangle | \\
		& \le \liminf_{\epsilon'\to 0} |\langle u_L^{\psi_{\epsilon'}^x}(t,\cdot), \psi_{\epsilon'}^x \rangle|^{\frac12}
|\langle u_L^{\psi_\epsilon^y - \psi_{\epsilon'}^y}(t,\cdot), \psi_\epsilon^y - \psi_{\epsilon'}^y\rangle|^{\frac12} \\
		& = | u_L^{\delta_x}(t,x)|^{\frac12}
|\langle u_L^{\psi_\epsilon^y - \delta_y}(t,\cdot), \psi_\epsilon^y - \delta_y\rangle|^{\frac12}.
	\end{align*}
	Now another application of Lemma~\ref{lemma:solution_map_box_continuous_in_inintial} shows that the right-hand side vanishes for $\epsilon \downarrow 0$.
\end{proof}

\subsection{A Feynman--Kac-type representation}
\label{subsection:feynman_kac_type}

Here we give a Feynman--Kac-type representation for $u_L^\phi$ with continuous initial condition $\phi$ and enhanced potential as introduced in \ref{obs:multiplicative_heat_eq_with_enhanced_potentials}. 
The technique behind this representation is inspired by the partial Girsanov transform of~\cite{GuPe17}, for a heuristic explanation see \ref{obs:heuristics_feynman}. 

Let us start by recalling the classical Feynman--Kac representation. 
We write $\P^x$ and $\E^x$ for the probability and expectation on $C([0,\infty),\R^2)$ such that the coordinate process $X= (X_t)_{t\in [0,\infty)}$  
is a Brownian motion with $X_0= x$ for $x\in \R^2$. 
Later, by an application of the Girsanov transform, we will consider a different measure, under which $X$ is a diffusion with non-trivial drift. 
For a probability measure $\Q$ on $C([0,\infty),\R^2)$ we write $\E_\Q$ for the expectation with respect to $\Q$. 

In this section we often abbreviate $\cC_\fn^{-\gamma}(Q_L)$ by $\cC_\fn^{-\gamma}$ and $\cC^{-\gamma}(\R^2)$ by $\cC^{-\gamma}$.

\begin{lemma}[Feynman--Kac representation]
\label{lemma:feynman_kac_box_smooth}
Let $L \in (0,\infty]$ and $Q_\infty = \R^2$. 
For $\phi \in C_{\rm b}(Q_L)$, 
 $\epsilon>0$ and $(t,y) \in [0,\infty) \times Q_L$ we have
\begin{align}
u_{L,\epsilon}^\phi (t,y) 
\label{eqn:feynman_kac_rep_smooth_box}
& = 
\E^y \bigg[ \exp \Big( \int_0^t (\theta_\epsilon -   c_\epsilon )(X_s) \dd s \Big) \phi(X_t) \1_{[X_{[0,t]}\subset Q_L]} \bigg],
\end{align}
where $X_{[0,t]} = \{ X_s : s\in [0,t]\}$. 
\end{lemma}

\begin{proof}
See for example \cite{Fr75,KaSh91,Sz98}.  
\end{proof}


For the application we have in mind, namely $\theta_\epsilon $ being a mollification of a typical realisation of the white noise, the limit is a distribution and $c_\epsilon$ diverges. Hence one cannot naively take the limit $\epsilon \downarrow 0$ in the right-hand side of \eqref{eqn:feynman_kac_rep_smooth_box} as the limit does not make sense. 
In the next observation we reformulate the Feynman--Kac formula in terms of $\theta_{L,\epsilon}$, on which we imposed assumptions on its limit as $\epsilon \downarrow 0$ (see \ref{obs:convergence_assumptions_theta}).

\begin{obs}
\label{obs:potential_replaced_by_periodic}
Let $L\in (0,\infty)$. As $\theta_\epsilon = \overline{\theta_{L,\epsilon}}$ (see \eqref{eqn:theta_L_eps_fn_projection}) almost everywhere in $Q_L$, we also have 
\begin{align}
\label{eqn:feynman_kac_rep_smooth_box_with_periodic_potential}
u_{L,\epsilon}^\phi (t,y) =  \E^y \bigg[ \exp \Big( \int_0^t (\overline{\theta_{L,\epsilon}}(X_s) -  c_\epsilon ) \dd s \Big) \phi(X_t) \1_{[X_{[0,t]}\subset Q_L]} \bigg] .
\end{align}
 Let $\theta_{L,0} \in \bigcap_{\alpha<-1} \cC_\fn^\alpha$ be the first component of $\btheta_L$, i.e., $\theta_{L,0} = \lim_{\epsilon \downarrow 0} \theta_{L,\epsilon}$. Then
$\overline{\theta_{L,0}} = \lim_{\epsilon \downarrow 0} \overline{\theta_{L,\epsilon}}$, see Lemma~\ref{lemma:equivalence_full_space_and_neumann}. 
\end{obs}

We will use \eqref{eqn:feynman_kac_rep_smooth_box_with_periodic_potential} to obtain our Feynman--Kac representation for $u_L^\phi = u_{L,0}^\phi$ in Theorem~\ref{theorem:feynman_kac_like_representation_for_limit}.

\begin{obs}[Heuristic idea]
\label{obs:heuristics_feynman}
Let us first heuristically derive the Feynman--Kac-type representation for an irregular potential, which we will rigorously prove in Theorem~\ref{theorem:feynman_kac_like_representation_for_limit}. For this, we neglect the renormalisation by $c_\epsilon$, consider $L=\infty$ for simplicity and do formal calculations. 
As mentioned, the integral $\int_0^t \theta(X_s) \dd s$ makes no sense as $\theta$ is of regularity $-1-$ (meaning: of regularity $-1-\delta$ for all $\delta>0$), therefore we want to replace it by something that does make sense. 
We do this by a trick introduced in \cite{GuPe17} that uses the It\^o formula. 
When we define $Z = (1- \frac12 \Delta)^{-1} \theta$, which has regularity $1-$, then we have by the It\^o formula 
\begin{align*}
Z(X_t) - Z(X_0) 
\cand \begin{calc} 
= \int_0^t \nabla Z \cdot \dd X_s + \int_0^t \frac12 \Delta Z \dd s 
\end{calc} \cnewline
\cand \begin{calc}
 = \int_0^t \nabla Z(X_s) \cdot \dd X_s - \int_0^t (1-\tfrac12 \Delta) Z(X_s) \dd s + \int_0^t Z(X_s) \dd s  
\end{calc} \cnewline
&  = \int_0^t \nabla Z(X_s) \cdot \dd X_s - \int_0^t \theta(X_s) \dd s + \int_0^t Z(X_s) \dd s . 
\end{align*}
So that 
\begin{align}
\notag &  \E^y \bigg[ \exp \Big( \int_0^t \theta(X_s) \dd s \Big) \phi(X_t) \bigg] \\
& = 
 \E^y \bigg[ \exp \Big(Z(X_0) - Z(X_t) + \int_0^t Z(X_s) \dd s + \int_0^t \nabla Z(X_s) \cdot \dd X_s \Big) \phi(X_t) \bigg]. 
 \label{eqn:to_be_girsanoved}
\end{align}
Now Girsanov's theorem tells us that $\tilde B_t = X_t - \int_0^t \nabla Z(X_s) \dd s$ is a Brownian motion under the probability measure $\Q_{Z}^y = \exp( \int_0^t \nabla Z(X_s) \cdot \dd X_s
-  \int_0^t \frac12|\nabla Z(X_s)|^2 \dd s) \dd \P^y $. 
Hence we can rewrite \eqref{eqn:to_be_girsanoved} by 
\begin{align*}
\E_{\Q_Z^y}\bigg[
 \exp \Big(Z(X_0) - Z(X_t) + \int_0^t Z(X_s) \dd s + \int_0^t \tfrac12 |\nabla Z(X_s)|^2 \dd s \Big) \phi(X_t)
\bigg], 
\end{align*}
where the coordinate process $X$ under $\Q_Z^y$ solves the SDE $\dd X_t = \nabla Z(X_t) \dd t + \dd \tilde B_t$. 
This is an improvement, since $\nabla Z$ has regularity $0-$ whereas $\theta$ has regularity $-1-$. However, the integral $\int_0^t |\nabla Z(X_s)|^2 \dd s$ does not yet make sense, because $\nabla Z$ is not a function but a distribution of regularity $0-$ and so, heuristically, $|\nabla Z|^2$ is also of regularity $0-$. 
Therefore we introduce another function $Y$ with regularity $2-$ and apply the Girsanov transformation to $Z+Y$  instead of $Z$, such that the transformed expectation has functions instead of distributions in the exponent and therefore makes sense. 
Let $Y$ be the solution of
\begin{align*}
(\eta - \tfrac12 \Delta) Y  = \tfrac12 |\nabla Z|^2 + \nabla Y \cdot \nabla Z, 
\end{align*}
where $\eta\in \R$ is chosen large enough so that $Y$ actually exists (Proposition~\ref{prop:solution_v_lambda}). 
Then we have 
\begin{align*}
\theta = (1-\tfrac12 \Delta) Z 
= Z + \eta Y - \tfrac12 \Delta ( Z+Y) - \tfrac12 |\nabla (Z+Y)|^2 + \tfrac12 |\nabla Y|^2. 
\end{align*}
Then the It\^o formula applied to $Z+Y$ gives 
\begin{align*}
\int_0^t \theta(X_s) \dd s 
&  = \int_0^t (Z+\eta Y + \tfrac12 |\nabla Y|^2)(X_s) \dd s + 
(Z+Y)(X_0) - (Z+Y)(X_t) \\
& \qquad    +  \int_0^t \nabla(Z+Y)(X_s) \cdot \dd X_s - \frac12 \int_0^t |\nabla (Z+Y)(X_s)|^2 \dd s. 
\end{align*}
So that, similarly to what we have done above, we can rewrite \eqref{eqn:to_be_girsanoved} by 
\begin{align*}
\E_{\Q_{Z+Y}^y} \bigg[
 \exp \Big( \int_0^t (Z+\eta Y + \tfrac12 |\nabla Y|^2)(X_s) \dd s + 
(Z+Y)(X_0) - (Z+Y)(X_t)  \Big) \phi(X_t)
\bigg]. 
\end{align*}
As $Y$ is of regularity $2-$ and $Z$ of regularity $1-$ all the terms in the exponential now do make sense, that is, $Z,Y$ and $\nabla Y$ are nice functions, and $X$ under $\Q_{Z+Y}^y$ solves the SDE $\dd X_t = \nabla (Z+Y)(X_t)\dd t + \dd B_t$, where  $B_t = X_t - \int_0^t \nabla (Z+Y)(X_s) \dd s$ is a Brownian motion under $\Q_{Z+Y}^y$. 

The above heuristics cannot be carried out as such for $L=\infty$. 
Instead we do this for each $L$ separately, and will need a control on the norms of $Z_L$ and $Y_L$, which diverge as $L\rightarrow \infty$ (hence $L=\infty$ cannot be done immediately). 
And in order to extend the representation \eqref{eqn:feynman_kac_rep_smooth_box} for $\epsilon=0$, we first rewrite it for $\epsilon>0$. 
\end{obs}

\begin{obs}
\label{obs:assumption_alpha_and_L}
 \textbf{For the rest of this section we fix $\alpha \in (-\frac43,-1)$ and $L>0$.} 
\end{obs}

For $\epsilon\in [0,\infty)$   we define
\begin{align*}
 Z_{L,\epsilon} : = \sigma(\rD) \overline{\theta_{L,\epsilon}} = (1-\tfrac12 \Delta)^{-1} \overline{\theta_{L,\epsilon}} .
\end{align*}
Observe that $Z_{L,\epsilon}$ is at least $C^2$ as $\overline{\theta_{L,\epsilon}}$ is continuous. 
To rewrite the Feynman--Kac representation for the solution to the  heat equation on $Q_L$ with potential $\theta_{\epsilon} - 2 c_\epsilon$ we will need to solve the equation $(\eta - \tfrac12 \Delta) Y = \tfrac12 |\nabla Z_{L,\epsilon}|^2 -  c_\epsilon  + \nabla Y \cdot \nabla Z_{L,\epsilon}$, where $\eta > 0$ is sufficiently large. For that purpose we need  bounds on $\nabla Z_{L,\epsilon}$ and $\frac12 |\nabla Z_{L,\epsilon}|^2 -  c_\epsilon $  that are uniform in $\epsilon$, and those can be obtained from bounds for the enhanced Neumann potential,  see Lemma~\ref{lemma:bound_on_grad_Z_eps_k_squared_min_renorm_and_Z_eps_k}.

\begin{lemma}
\label{lemma:bound_on_grad_Z_eps_k_squared_min_renorm_and_Z_eps_k}
Let $\sigma$ be as in Definition~\ref{def:enhanced_neumann_potentials}. 
There exists a $C>0$ (independent of $L$) such that for all $\epsilon>0$ 
\begin{align}
\notag
& \| \tfrac12 |\nabla Z_{L,\epsilon}|^2 -   c_\epsilon  \|_{\cC^{2\alpha+2}} + \| \nabla Z_{L,\epsilon} \|_{\cC^{\alpha+1}} \\
\label{eqn:bound_on_grad_Z_eps_k_squared_min_renorm_and_Z_eps_k}
& \le C (\| \theta_{L,\epsilon}\|_{\cC_\fn^{\alpha}} +\| \theta_{L,\epsilon}\|_{\cC_\fn^{\alpha}}^2 + \| \theta_{L,\epsilon} \reso \sigma(\rD) \theta_{L,\epsilon} - c_\epsilon \|_{\cC_\fn^{2\alpha+2}}). 
\end{align}
Moreover, 
 with $\Theta_L$ the second component of $\btheta_L$   
and by using the notation 
\begin{align*}
\tfrac12 |\nabla Z_L|^{\diamond 2} 
:= \nabla Z_L \para \nabla Z_L
- (1- \tfrac14 \Delta) (Z_L \reso Z_L) 
 + \Theta_L ,
\end{align*} 
we have 
\begin{align}
\label{eqn:convergence_nabla_Z_and_square_min_c}
	( \nabla Z_{L,\epsilon},\tfrac12 |\nabla Z_{L,\epsilon}|^2 -   c_\epsilon ) \to (\nabla Z_L,\tfrac12 |\nabla Z_L|^{\diamond 2}) .
\end{align}
\end{lemma}

\begin{proof}
To shorten notation we write $Z=Z_{L,\epsilon}$,  $\theta = \theta_{L,\epsilon}$, $\overline \theta= \overline{\theta_{L,\epsilon}}$ and $c= c_\epsilon$. 
By Lemma~\ref{lemma:equivalence_full_space_and_neumann} we have $\| \theta\|_{\cC_\fn^{-\gamma}} = \| \overline{\theta} \|_{\cC^{-\gamma}} $ and $\| \theta \reso \sigma(\rD) \theta - c \|_{\cC_\fn^{-\gamma}} = \| \overline{\theta} \reso \sigma(\rD) \overline{\theta} - c \|_{\cC^{-\gamma}}$.  

We have $\frac12| \nabla Z|^2 =  \nabla Z \para \nabla Z + \frac12 \nabla Z \reso \nabla Z$, where the para- and resonant product are combined with the inner product on $\R^2$, e.g. $ \nabla Z \para \nabla Z = \partial_1 Z \para \partial_1 Z + \partial_2 Z \para \partial_2 Z$. 
For the resonant product we  use that $(1-\frac12 \Delta) Z = \overline{\theta}$ and  apply Leibniz's rule to obtain:
\begin{calc}
First of all 
\begin{align*}
\Delta (Z \reso Z) = 2 \nabla Z \reso \nabla Z + 2 (\Delta Z) \reso Z, 
\end{align*}
therefore 
\end{calc}
\begin{align*}
\tfrac12 \nabla Z \reso \nabla Z 
\notag & = \tfrac14 \Delta (Z \reso Z) - (\tfrac12 \Delta Z) \reso Z \\
\cand \begin{calc}\notag
 = \tfrac14 \Delta (Z \reso Z) +  ((1- \tfrac12 \Delta) Z) \reso Z - Z \reso Z 
 \end{calc} \cnewline
 \cand \begin{calc}\notag
 = (\tfrac14 \Delta-1) (Z \reso Z) +  \overline \theta \reso Z
\end{calc} \cnewline
& = - (1- \tfrac14 \Delta) (Z \reso Z) +  \overline{\theta} \reso \sigma(\rD) \overline{\theta}.
\end{align*}

Then the claim follows since by the H\"ormander-Mikhlin bound, Bernstein's inequality and the Bony estimates for the paraproduct and resonance product 
\cite[Lemma~2.1, Lemma~2.2, Theorem~2.82 and Lemma~2.85]{BaChDa11}
we have 
$	\| \nabla Z \para \nabla Z \|_{\cC^{2\alpha+2}} \lesssim  \| \nabla Z  \|_{\cC^{\alpha+1}}^2 \lesssim \| Z \|_{\cC^{\alpha+2}}^2\lesssim \| \theta \|_{\cC^{\alpha}}^2$, \quad
 $\|  (1- \frac14 \Delta) (Z\reso Z)\|_{\cC^{2\alpha+2}} \lesssim 
\|  Z\reso Z \|_{\cC^{2\alpha+4}} 
\lesssim \|  Z \|_{\cC^{\alpha+2}}^2 \lesssim \| \overline{\theta}  \|_{\cC^{\alpha}}^2$, \quad $\| Z  \|_{\cC^{\alpha+2}} \lesssim \| \overline{\theta}  \|_{\cC^{\alpha}}$ and $ \| \nabla Z \|_{\cC^{\alpha+1}} \lesssim \| Z  \|_{\cC^{\alpha+2}}$. 
 \eqref{eqn:convergence_nabla_Z_and_square_min_c} then follows from the convergence of \eqref{eqn:limit_of_enhanced_theta_L_eps_fn}. 
\end{proof}

We show in Lemma~\ref{prop:solution_v_lambda} in the appendix that there exists a $C>0$ such that for all $M>0$ for all $f \in \cC^{2\alpha+2}$ and $g \in \cC^{\alpha+1}$ with $\|f\|_{\cC^{2\alpha+2}}, \|g\|_{\cC^{\alpha+1}} \le M$, 
 for all $\beta \in(-\alpha, 2\alpha+4)$ and for all $\eta \ge C (1 + M)^{\frac{2}{2\alpha + 4 - \beta}}$
there exists a unique solution $v\in \cC^\beta(\R^d)$ to
\begin{align*}
	\left( \eta - \tfrac{1}{2} \Delta \right) v = f + \nabla v \cdot g,
\end{align*}
such that $\|v\|_{\cC^\beta} \le M$, and that $v$ depends continuously on $f$ and $g$. 

Fix $\beta \in(-\alpha, 2\alpha+4)$. Let 
\begin{align}
\label{eqn:def_M_L_eps}
	M_{L,\epsilon} &:= \max\{\|\tfrac12 |\nabla Z_{L,\epsilon}|^2 -  2 c_\epsilon  \|_{\cC^{2\alpha+2}}, \|\nabla Z_{L,\epsilon}\|_{\cC^{\alpha+1}}\}, \\
\label{eqn:def_M_L}
	M_L &:= \max\{\|\tfrac12 |\nabla Z_L|^{\diamond 2} \|_{\cC^{2\alpha+2}}, \|\nabla Z_L \|_{\cC^{\alpha+1}} \}.
\end{align}

By \eqref{eqn:convergence_nabla_Z_and_square_min_c} we have $M_{L,\epsilon} \rightarrow M_L$. 
For $\epsilon> 0$ we define $\eta_{L,\epsilon} := C(1+M_{L,\epsilon})^{\frac{2}{2\alpha + 4 - \beta}}$ and $\eta_L := C(1+M_{L})^{\frac{2}{2\alpha + 4 - \beta}}$. 
Therefore there exists a solution $Y_{L,\epsilon}$ to 
\begin{align}
\label{eqn:equation_for_Y_eps_L_eta}
(\eta_{L,\epsilon} - \tfrac12 \Delta) Y_{L,\epsilon}  = \tfrac12 |\nabla Z_{L,\epsilon}|^2 -  c_\epsilon  + \nabla Y_{L,\epsilon} \cdot \nabla Z_{L,\epsilon},
\end{align}
 with $\|Y_{L,\epsilon}\|_{\cC^\beta}\le M_{L,\epsilon}$
Moreover, in $\cC^\beta$ we have $Y_{L,\epsilon} \rightarrow Y_L$ as $\epsilon \downarrow 0$, where $Y_L$ solves
\begin{align}
\label{eqn:equation_for_Y_L_eta}
(\eta_L - \tfrac12 \Delta) Y_L  = \tfrac12 |\nabla Z_L|^{\diamond 2} + \nabla Y_L \cdot \nabla Z_L,
\end{align}
and $\|Y_L\|_{\cC^\beta} \le M_L$.
We will also write ``$ Z_{L,0},M_{L,0},\eta_{L,0}, Y_{L,0}$'' for ``$Z_{L},M_{L},\eta_{L}, Y_{L}$''.

Let us now discuss how $u_{L,\epsilon}^\phi $ 
 can be described in terms of $Z_{L,\epsilon}$ and $Y_{L,\epsilon}$. For $0\le r \le t$ and $\epsilon \in [0,\infty)$ we write 
\begin{align}
\label{eqn:cD_L_formula}
\cD_{L,\epsilon}(r,t) := e^{\int_r^t (Z+\eta Y + \frac12 |\nabla Y|^2)(X_s) \dd s + 
(Z+Y)(X_r) - (Z+Y)(X_t)},
\end{align}
where we have abbreviated ``$Z_{L,\epsilon},\eta_{L,\epsilon},Y_{L,\epsilon}$'' by ``$Z,\eta,Y$''.

\begin{lemma}
\label{lemma:feynman_kac_like_representation_on_box_with_epsilon}
Let $\epsilon > 0$ and $y \in \R^2$ and let $\Q_{L,\epsilon}^y$ be the probability measure on $C([0,\infty),\R^2)$ such that the coordinate process satisfies $\Q_{L,\epsilon}^y$-almost surely 
\begin{align*}
	X_t = y + \int_0^t \nabla (Z_{L,\epsilon}+Y_{L,\epsilon})( X_s) \dd s + B_t, \qquad t \ge 0,
\end{align*}
 for a Brownian motion $B_t$. 
 Then we have for all $\epsilon > 0$ and for any measurable and bounded (or positive) functional $F \colon C([0,t],\R^2) \to \R$
\begin{align*}
& \E^y \left[ e^{ \int_0^t (\overline{\theta_{L,\epsilon}}(X_s) -  c_\epsilon  ) \dd s } F(X|_{[0,t]})\1_{[ X_{[0,t]}  \subset Q_L]} \right] 
 = 
\E_{\Q_{L,\epsilon}^y} \left[ \cD_{L,\epsilon}(0,t) F(X|_{[0,t]})\1_{[ X_{[0,t]}  \subset Q_L]}\right].
\end{align*}
In particular, $u_{L,\epsilon}^\phi(t,y)= \E_{\Q_{L,\epsilon}^y} \left[ \cD_{L,\epsilon}(0,t)  \phi(X_t) \1_{[ X_{[0,t]} \subset Q_L]}\right]$.
\end{lemma}
\begin{proof}
To shorten notation we write $Z= Z_{L,\epsilon}$, $Y= Y_{L,\epsilon}$, $\theta = \overline{\theta_{L,\epsilon}}$ and $c= c_\epsilon$ in this proof. By definition of $Z$ and $Y$, we have $\theta = Z - \frac12 \Delta Z$ and 
\begin{align*}
 -  c  =  \eta Y - \tfrac12 \Delta Y  - \tfrac12 | \nabla Z|^2  - \nabla Y \cdot \nabla Z 
 =  \eta Y - \tfrac12 \Delta Y - \tfrac12 | \nabla (Z+Y)|^2 + \tfrac12 |\nabla Y|^2,
\end{align*}
so that 
\begin{align*}
\theta -  c  = Z+\eta Y - \tfrac12 \Delta (Z+Y) - \tfrac12 | \nabla (Z+Y)|^2  + \tfrac12 |\nabla Y|^2.
\end{align*}
Hence, on $[X_{[0,t]}\subset Q_L]$, by using It\^o's formula to rewrite $\int_0^t \frac12 \Delta (Z+Y)(X_s) \dd s$ we have
\begin{calc}
\begin{align*}
(Z+Y)(X_t) = (Z+Y)(X_0) + \int_0^t \tfrac12 \Delta (Z+Y)(X_s) \dd s + \int_0^t \nabla (Z+Y)(X_s) \cdot \dd X_s. 
\end{align*}
\end{calc}
\begin{align*}
\int_0^t \theta(X_s) -  c  \dd s 
&  = \int_0^t (Z+\eta Y + \tfrac12 |\nabla Y|^2)(X_s) \dd s + 
(Z+Y)(X_0) - (Z+Y)(X_t) \\
& \qquad    +  \int_0^t \nabla(Z+Y)(X_s) \cdot \dd X_s - \frac12 \int_0^t |\nabla (Z+Y)(X_s)|^2 \dd s. 
\end{align*}
Since
\begin{align}
\label{eqn:girsanov_density}
G_t := \exp ( \int_0^t \nabla(Z+Y)(X_s) \cdot \dd X_s - \frac12 \int_0^t |\nabla (Z+Y)(X_s)|^2 \dd s )
\end{align}
is a Radon-Nikodym density, the claim follows from Girsanov's theorem  (see e.g. Le Gall \cite[\S 5.6 Page 138]{LG16}).
The formula for $u_{L,\epsilon}^\phi$ then follows from \eqref{eqn:feynman_kac_rep_smooth_box_with_periodic_potential}. 
\end{proof}

\begin{theorem}
\label{theorem:feynman_kac_like_representation_for_limit}
Let $y \in Q_L^\circ$. 
There exists a unique probability measure $\Q_L^y$ on $C([0,\infty),\R^2)$ under which the coordinate process $X$ satisfies $X_0 = y$ almost surely, and $X$ solves the martingale problem for the following SDE with distributional drift (see for example \cite[Definition 2.1]{PevZ})
\begin{align}
\label{eqn:X_t_with_grad_Z_L_Y_L_as_drift}
\dd X_t = \nabla (Z_L + Y_L)( X_t) \dd t + \dd B_t
\end{align}
where $B$ is a Brownian motion.  
Moreover,
with  $\cD_L = \cD_{L,0}$ (see \eqref{eqn:cD_L_formula}), 
 with $\phi \in C(Q_L)$, and with $A= [ X_{[0,t]} \not \subset S,  X_{[0,t]} \subset Q_L]$ for $S=\emptyset$ or $S = Q_r$ for $r\in (0,L)$, we have 
\begin{align}
 & \E^y \left[ e^{ \int_0^t (\theta_\epsilon(X_s) -  c_\epsilon  ) \dd s } \phi(X_t) \1_{A} \right]  \xrightarrow{\epsilon \downarrow 0}
 \E_{\Q_L^y} \left[ \cD_L(0,t)  \phi(X_t) \1_{A}\right].
\label{eqn:feynman_kac_like_with_A}
\end{align} 
In particular, $ u_L^\phi(t,y)= \E_{\Q_L^y} \left[ \cD_L(0,t)  \phi(X_t) \1_{[ X_{[0,t]} \subset Q_L]}\right]$.
\end{theorem}

\begin{proof}
The existence and uniqueness of the solution $\Q_L^y$ to the martingale problem follows from \cite[Theorem 1.2]{CaCh18}. 
  It is not mentioned in that theorem, but from the proof of \cite[Theorem 4.3]{CaCh18} one can extract that the martingale depends continuously on the drift (in the space $\cC^{-\alpha}(\R^d,\R^d)$.   
Therefore, as $\nabla (Z_{L,\epsilon} + Y_{L,\epsilon})$ converges to $\nabla (Z_L+Y_L)$, 
we have that $\Q_{L,\epsilon}^y$ converges weakly to $\Q_L^y$ (where $C([0,\infty), \R^2)$ is equipped with the topology of uniform convergence on compacts). 
By Lemma~\ref{lemma:feynman_kac_like_representation_on_box_with_epsilon}
	\begin{align*}
		\E^y \left[ e^{ \int_0^t (\theta_\epsilon(X_s) -  c_\epsilon  ) \dd s } \phi(X_t) \1_{A} \right] = \E_{\Q_{L,\epsilon}^y} \left[ \cD_{L,\epsilon}(0,t) \phi(X_t) \1_{A}\right].
	\end{align*}
Since $Z_{L,\epsilon} \rightarrow Z_L$ and $Y_{L,\epsilon} \rightarrow Y_L$ uniformly, we get
	\begin{align*}
		\lim_{\epsilon \downarrow 0} \E_{\Q_{L,\epsilon}^y} \left[ \cD_{L,\epsilon}(0,t) \phi(X_t) \1_{A}\right] = \lim_{\epsilon \downarrow 0} \E_{\Q_{L,\epsilon}^y} \left[ \cD_L(0,t) \phi(X_t) \1_{A}\right].
	\end{align*}
	By the Portmanteau theorem \cite[Theorem 2.1]{Bi68}, the convergence in~\eqref{eqn:feynman_kac_like_with_A} follows once we show that $A$ is a $\Q_L^y$ continuity set, i.e., $\Q_L^y[\partial A] = 0$. We treat the case $S=\emptyset$, the case $S=Q_r$ follows from very similar but slightly more tedious arguments. 
The set $A$ consists of paths that stay inside the closed box $Q_L$ until time $t$, so $\partial A = [ X_{[0,t]}\not \subset Q_L^\circ , X_{[0,t]} \subset Q_L ]$. To prove that this set has probability zero we show that ``almost surely $X$ leaves the box right after it hits the boundary''. 

So let $\tau  = \inf \{ s\in [0,t] : X_s \in \partial Q_L\}$, with which we have $\partial A = [ \tau \le t , X_{[0,t]} \subset Q_L]$. First note that $\Q_L^y[\tau = t] \le \Q_L^y[X_t \in \partial Q_L] = 0$, since 
 the law of $X_t$ under $\Q_L^y$ is absolutely continuous to the Lebesgue measure by \cite[Proposition 2.9]{PevZ}. So it remains to show that $\Q_L^y[\tau < t , X_{[0,t]} \subset Q_L] = 0$. By the strong Markov property, we have
\begin{align*}
\Q_L^y[\tau < t , X_{[0,t]} \subset Q_L] 
\cand \begin{calc}
= \E_{\Q_L^y} [ \1_{\tau < t} \Q_L^y [ X_{[0,t]} \subset Q_L | \cF_{\tau} ]]
\end{calc} \cnewline
\cand \begin{calc}
= \E_{\Q_L^y} [ \1_{\tau < t} \Q_{X_\tau} [ X_{[0,t-\tau]} \subset Q_L ]]
\end{calc} \cnewline
& \le \E_{\Q_L^y}[ \1_{\tau < t} \Q_{X_\tau}[ \exists \epsilon > 0:  X{[0,\epsilon]} \subset Q_L]]. 
\end{align*}
Similar arguments as in~\cite[Theorem 13, estimate (48)]{DeDi16} show that $X-B$ is $\Q_{X_\tau}$-almost surely $\frac{2-\alpha-\epsilon}{2}$-H\"older continuous for all $\epsilon > 0$. 
We can choose $\epsilon>0$ small enough such that $\frac{2-\alpha-\epsilon}{2} > \frac12$. 
Therefore, by the law of the iterated logarithm for $B$ (by the strong Markov property $B$ is a Brownian motion under $\Q_{X_\tau}$) 
it follows that that $\Q_{X_\tau}$-almost surely there is a sequence of times $(t_n)_{n\in\N}$ such that $t_n \downarrow 0$, 
\begin{align*}
\limn \frac{|X_{t_n} - X_0 - B_{t_n}|}{|B_{t_n}|} 
\le  \limn \frac{t_n^{\frac{2-\alpha-\epsilon}{2}}}{\sqrt{2t_n \log \log \frac{1}{t}} } =0, 
\end{align*}
and also $X_0 + B_{t_n} \notin Q_L$ for all $n$. 
\begin{calc}
Indeed, the law of the iterated logarithm \cite[Theorem 5.5.14]{CoEl15} implies that almost surely for the first (and thus each) component of the multidimensional Brownian motion
\begin{align*}
\limsup_{t\downarrow 0} \frac{B_t^1}{ \sqrt{2t \log \log \frac{1}{t}}} =1 , \mbox{ and }
\liminf_{t\downarrow 0} \frac{B_t^1}{ \sqrt{2t \log \log \frac{1}{t}}} =-1 . 
\end{align*}
As $X_0$ lies on the boundary of $Q_L$, we can choose a coordinate $i$ and a positive or negative direction such that $X_0$ plus or minus $e_i$ is outside $Q_L$. 
\end{calc}
Therefore,
\begin{align*}
	X_{t_n} = X_0 + |B_{t_n}| \left(\frac{B_{t_n}}{|B_{t_n}|} + \frac{X_{t_n} - X_0 - B_{t_n}}{|B_{t_n}|}\right) \notin Q_L
\end{align*}
for all sufficiently large $n$. 
\begin{calc}
Let $i$ be such that $|\langle X_0 + B_{t_n}, e_i\rangle| > L$. 
We may assume that $\langle X_0 , e_i \rangle = \frac{L}{2}$. 
Then $\langle B_{t_n} , e_i \rangle >0$ for large $n$. 
If $n$ is taken large enough such that 
\begin{align*}
\frac{|X_{t_n} - X_0 - B_{t_n}|}{|B_{t_n}|}<1, \quad \mbox{then} \quad 
\langle \frac{B_{t_n}}{|B_{t_n}|} + \frac{X_{t_n} - X_0 - B_{t_n}}{|B_{t_n}|} , e_i \rangle >0 . 
\end{align*}
\end{calc}
This concludes the proof.
\end{proof}

\subsection{The Feynman--Kac-type representation for the PAM}
\label{subsection:feynman_kac_for_pam_with_white_noies}

In this section we show how the above representation techniques can be applied to the PAM with white-noise potential. 
We give the Feynman--Kac-type representation for the PAM in Theorem~\ref{theorem:feynman_kac_type_representation_for_PAM} and state estimates which are important for the next sections in Lemma~\ref{lemma:escape_probability_of_box_L_with_drift_Z_and_Y_L}. 

\begin{definition}
A \emph{white noise} on $\R^d$ is a random variable $\xi$ on a probability space $(\Omega,P)$ with values in $\cS'(\R^d)$ such that  for all $f\in \cS$ the random variable $\langle \xi,f \rangle$ is a centered Gaussian random with $E [ \langle \xi, f \rangle  \langle \xi, g \rangle ] = \langle f, g \rangle_{L^2(\R^d)}$ for $f, g \in \cS$. 
\end{definition}

Like in Section~\ref{subsection:feynman_kac_type} we will have the Brownian motion as another source of randomness, besides the noise term $\xi$ as the potential of the PAM. 
These two are independent sources of randomness, therefore we write ``$P$'' and ``$E$'' for the probability and expectation with respect to $\xi$, contrary to the notation ``$\P$'' and ``$\E$'' that we use for the probability and expectation with respect to the Brownian motion $B$ (as in Section~\ref{subsection:feynman_kac_type}).

\begin{obs}
\label{obs:setting_representation_section}
Let $\psi \in C_{\rm c}^\infty(\R^2)$ be given by $\psi (x) = \varphi(x_1) \varphi(x_2)$ for an even function $\varphi \in C_{\rm c}^\infty(\R, [0,\infty))$ with $\supp \varphi \subset [-\frac12,\frac12]$ and $\int_\R \varphi (x) \dd x= 1$.  
We write $\psi_\epsilon(x) = \frac{1}{\epsilon^2} \psi( \frac{x}{\epsilon})$ and $\varphi_\epsilon(x) = \frac{1}{\epsilon} \varphi(\frac{x}{\epsilon})$.  Let $\xi$ be a white noise on $\R^2$ and let $\xi_\epsilon = \psi_\epsilon * \xi$ be its mollification. 
\end{obs}

\begin{obs}
\label{obs:pam_smooth_to_white_noise}
Let
\begin{align}
\label{eqn:xi_L_eps_fn_projection}
\xi_{L,\epsilon}^\fn = \sum_{n\in\N_0^2} \langle \xi_\epsilon, \fn_{k,L} \rangle \fn_{k,L}. 
\end{align} 
Now we want to apply the deterministic results of Section~\ref{subsection:feynman_kac_type} 
to the white-noise potential, in the sense that we apply the above for $\theta_{L,\epsilon}$ being a typical realisation of $\xi_{L,\epsilon}^\fn$.  
In order to do so, we will need to have a convergence for $\xi_{L,\epsilon}^\fn$ as we had assumed for $\theta_{L,\epsilon}$ in \ref{obs:convergence_assumptions_theta}.  
We show in Theorem~\ref{theorem:convergence_mollified_noise}  that there exists a $C\in\R$ such that for 
\begin{align*}
c_\epsilon = \tfrac{1}{\pi} \log \tfrac{1}{\epsilon} + C , 
\end{align*}
there  exists a $ \bxi_L = (\xi_L,  \Xi_L) \in \bigcap_{\alpha < -1} \fX_\fn^{\alpha}$ such that the following convergence holds in $\fX_\fn^{\alpha}$ for all $\alpha < -1$: 
\begin{align}
\label{eqn:limit_of_enhanced_xi_L_eps_fn}
(\xi_{L,\epsilon}^\fn, \xi_{L,\epsilon}^\fn \reso \sigma(\rD) \xi_{L,\epsilon}^\fn - c_\epsilon  ) \xrightarrow[\epsilon \downarrow 0]{P} \bxi_L. 
\end{align}
\end{obs}

Analogously to the notation in Section~\ref{subsection:feynman_kac_type}, 
we let $\alpha \in (-\frac43, -1)$ and $\beta \in (-\alpha, 2\alpha +4)$;  
for $L>0$ and $\epsilon>0$ we write $Z_{L,\epsilon}$ and $Z_L$ to be the (random) functions given by 
\begin{align*}
	Z_{L,\epsilon} = \sigma(\rD) \overline{\xi_{L,\epsilon}^\fn} = (1- \tfrac12 \Delta)^{-1} \overline{\xi_{L,\epsilon}^\fn}, \qquad Z_L = \sigma(\rD) \overline{\xi_{L,0}} = (1-\tfrac12\Delta)^{-1} \overline{\xi_{L,0}},
\end{align*}
where $\xi_{L,0}$ is the first component of $\bxi_L$; 
$M_{L,\epsilon}$ is as in \eqref{eqn:def_M_L_eps};
$M_L$ as in \eqref{eqn:def_M_L};
$\eta_{L,\epsilon} := C(1+M_{L,\epsilon})^{\frac{2}{2\alpha + 4 - \beta}}$;
$\eta_L := C(1+M_{L})^{\frac{2}{2\alpha + 4 - \beta}}$;
$Y_{L,\epsilon}$ is the solution to \eqref{eqn:equation_for_Y_eps_L_eta}; 
and, $Y_L$ the solution to \eqref{eqn:equation_for_Y_L_eta}. 
Again, we will also write ``$ Z_{L,0},M_{L,0},\eta_{L,0}, Y_{L,0}$'' for ``$Z_{L},M_{L},\eta_{L}, Y_{L}$''. 
For $L>0$ and $\epsilon\ge 0$ we let $\cD_{L,\epsilon}$ be as in \eqref{eqn:cD_L_formula} and write ``$\cD_L$'' for ``$\cD_{L,0}$''.

\begin{theorem}
\label{theorem:feynman_kac_type_representation_for_PAM}
We have 
\begin{align}
& \notag 
Z_{L,\epsilon} \xrightarrow[\epsilon \downarrow 0]{P} Z_L, 
\quad 
M_{L,\epsilon} \xrightarrow[\epsilon \downarrow 0]{P} M_L, 
\quad 
\eta_{L,\epsilon} \xrightarrow[\epsilon \downarrow 0]{P} \eta_L, 
\quad 
Y_{L,\epsilon} \xrightarrow[\epsilon \downarrow 0]{P} Y_L. 
\end{align}
Furthermore, for $A= [ X_{[0,t]} \not \subset S,  X_{[0,t]} \subset Q_L]$ with $S=\emptyset$ or $S = Q_r$ for $r\in (0,L)$, we have 
\begin{align}
 & \E^y \left[ e^{ \int_0^t (\xi_\epsilon(X_s) -  c_\epsilon  ) \dd s } \phi(X_t) \1_{A} \right]  
 \xrightarrow[\epsilon \downarrow 0]{P}
 \E_{\Q_L^y} \left[ \cD_L(0,t)  \phi(X_t) \1_{A}\right],
\label{eqn:feynman_kac_like_with_A_for_white_noise}
\end{align} 
in particular $u_L^\phi = \E_{\Q_L^y} \left[ \cD_L(0,t)  \phi(X_t) \1_{[ X_{[0,t]} \subset Q_L]}\right]$.
\end{theorem}
\begin{proof}
These convergences follow from the deterministic convergences of Section~\ref{subsection:feynman_kac_type} and by using the fact that a sequence converges in probability to a random variable if and only if each subsequence has a further subsequence that converges almost surely to that random variable. 
\end{proof}

In the following section, we will use some estimates that we present below in Lemma~\ref{lemma:escape_probability_of_box_L_with_drift_Z_and_Y_L}. 
We could estimate the right-hand side of \eqref{eqn:feynman_kac_like_with_A_for_white_noise} by estimating the $Z_{L,\epsilon}, \eta_{L,\epsilon}, Y_{L,\epsilon}$ in terms of $M_{L,\epsilon}$. 
We will need a control on the growth of $M_{L,\epsilon}$, which is provided in the following lemma. 

\begin{lemma}
\label{lemma:a_epsilon_bound_of_M_L_eps}
For $\epsilon \ge 0$ the following random variable is almost surely finite 
\begin{align}
\label{eqn:a_eps}
a_{\epsilon}
& := 1\vee \sup_{L\in \N, L > e}  \frac{M_{L,\epsilon} }{\log L},
\end{align}
and there exists an $h_0>0$ such that $\sup_{\epsilon \ge 0 } E[ e^{h a_{\epsilon}}]<\infty$ for all $h \in [0,h_0]$.
\end{lemma}
\begin{proof}
By Lemma~\ref{lemma:bound_on_grad_Z_eps_k_squared_min_renorm_and_Z_eps_k} we can bound for $\epsilon > 0$
\begin{align*}
	M_{L,\epsilon} \le C (\| \xi_{L,\epsilon}\|_{\cC_\fn^{\alpha}} +\| \xi_{L,\epsilon}\|_{\cC_\fn^{\alpha}}^2 + \| \xi_{L,\epsilon} \reso \sigma(\rD) \xi_{L,\epsilon} - c_\epsilon \|_{\cC_\fn^{2\alpha+2}}),
\end{align*}
and similarly for $M_{L,0}$. So the claim follows from Lemma~\ref{lemma:growth_of_noise_in_terms_of_box_size} below.
\end{proof}

\begin{obs}\label{obs:bounds_on_ZL_YL_thetaL}
So far $\alpha \in (-\tfrac43,-1)$ and $\beta \in (-\alpha, 2\alpha+4)$ were arbitrary.
\begin{calc}
Observe that this implies $-\beta \in (- 2\alpha -4, \alpha)$ and $2\alpha +4 \in (4-\frac{8}{3}, 2)$ and thus $2\alpha + 4 - \beta \in (0, 2+ \alpha)\subset (0,1)$.
\end{calc} 
But we can and do from now assume that $\alpha$ ``is close enough to $-1$'' and $\beta$ ``is close enough to $-\alpha$'' so that
\begin{calc}
(observe that for $\alpha' =  - \frac98$, we have $2\alpha' +4 = \frac{14}{8}= \frac74$ and thus for $\beta' =\frac{10}{8} = \frac54$ we have $2\alpha' + 4 - \beta = \frac{4}{8} = \frac12$, so take for example $\alpha \in (-\frac98,-1)$ and $\beta \in (-\alpha, \frac54)$)
\end{calc}
\begin{align}
\label{eqn:assumption_alpha_beta}
	\frac{2}{2\alpha+4-\beta} \le 4.
\end{align}
Note that with this assumption, Lemma~\ref{lemma:a_epsilon_bound_of_M_L_eps} implies for all $\epsilon \ge 0$ and $L> e$ (so that $a_\epsilon \log L > 1$), 
\begin{align*}
	\|Z_{L,\epsilon}\|_{\cC^{\alpha+2}} \le a_\epsilon \log L , \qquad \|Y_{L,\epsilon}\|_{\cC^\beta} \le  a_\epsilon \log L,\qquad \eta_{L,\epsilon} \le 2C (a_\epsilon \log L)^4.
\end{align*}
As $\| |\nabla Y_{L,\epsilon}|^2 \|_{\cC^{\beta-1}} \le \| Y_{L,\epsilon} \|_{\cC^\beta}^2$ we can estimate all terms under the integral in the definition of $\cD_{L,\epsilon}$  \eqref{eqn:cD_L_formula} as $\|\cdot\|_{L^\infty} \lesssim \|\cdot\|_{\cC^\delta}$ for $\delta>0$. The estimate we obtain from this for $\cD_{L,\epsilon}$ and another estimate will be presented in the following lemma. 
This lemma will play a key role in Section~\ref{section:separating_space_in_boxes}. 
\end{obs}

\begin{lemma}
\label{lemma:escape_probability_of_box_L_with_drift_Z_and_Y_L}
Let $\Q_{L,\epsilon}^y$ be the probability measures from Lemma~\ref{lemma:feynman_kac_like_representation_on_box_with_epsilon} and let us write $\Q_{L,0}^y$ for the probability measure $\Q_L^y$ from Theorem~\ref{theorem:feynman_kac_like_representation_for_limit}. 
Let $a_\epsilon$ be as in \eqref{eqn:a_eps}. 
There exists a $C>1$ such that for all $L\in\N$ with $L>e$, $\epsilon \ge 0$, $r \in (0,L)$ and $t\ge 1$:
\begin{align}
\label{eqn:bound_cD_L}
e^{-C  a_\epsilon^5 t (\log L)^{5}} 
\1_{[X_{[0,t]}\subset Q_L]}
 \le \cD_{L,\epsilon}(0,t) \1_{[X_{[0,t]}\subset Q_L]}
 & \le e^{C a_\epsilon^5 t (\log L)^{5}}, \\
\label{eqn:escape_probability_of_box_L_with_drift_Z_and_Y_L}
 \Q_{L,\epsilon}^0 [ X_{[0,t]}\not \subset Q_r] 
&\le  C e^{C a_\epsilon^5 t (\log L)^5 - \frac{r^2}{Ct}}, \\
\label{eqn:bound_feynman_part}
 \E_{\Q_{L,\epsilon}^0} \left[ \cD_{L,\epsilon}(0,t) 
 \1_{[ X_{[0,t]}\not \subset Q_r, X_{[0,t]} \subset Q_L ] }\right] 
& \le C e^{C a_\epsilon^5 t (\log L)^{5} - \frac{r^2}{Ct}}.
\end{align}
\end{lemma}
\begin{proof}
\eqref{eqn:bound_cD_L} follows immediately from the definition of $\cD_{L,\epsilon}(0,t)$ in~\eqref{eqn:cD_L_formula} together with Observation~\ref{obs:bounds_on_ZL_YL_thetaL}. 

By \cite[Corollary 1.2]{PevZ}, where we take $b = \nabla (Z_{L,\epsilon} + Y_{L,\epsilon})$ of regularity $\alpha+1$, $\delta>0$ and use that
\begin{align*}
\|\Delta_{-1} b\|_{L^\infty}^2 + \|\Delta_{\ge 0 } b \|_{B^{1+\alpha-\delta}_{\infty,1}}^{\frac{2}{2+\alpha-\delta}} \lesssim 1 + \|b \|_{B^{1+\alpha-\delta}_{\infty,1}}^{\frac{2}{2+\alpha-\delta}} \lesssim 1+  \|b \|_{\cC^{1+\alpha}}^{\frac{2}{2+\alpha-\delta}},
\end{align*}
we have 
\begin{align*}
	\Q_{L,\epsilon}^0 [ X_{[0,t]}\not \subset Q_r] \le C \exp \bigg( C t(1+\| \nabla (Z_{L,\epsilon} + Y_{L,\epsilon}) \|_{\cC^{\alpha+1}}^{\frac{2}{2+\alpha-\delta}}) - \frac{r^2}{C t} \bigg). 
\end{align*}
Since $\alpha>-\tfrac43$ we can choose $\delta>0$ small enough so that $\frac{2}{2+\alpha-\delta} \le 4$. 
Then~\eqref{eqn:escape_probability_of_box_L_with_drift_Z_and_Y_L} follows from another application of Observation~\ref{obs:bounds_on_ZL_YL_thetaL}. Finally, \eqref{eqn:bound_feynman_part} follows from \eqref{eqn:bound_cD_L} and \eqref{eqn:escape_probability_of_box_L_with_drift_Z_and_Y_L}.
\end{proof}

\section{Asymptotic behaviour of the mass on a box: Proof of Theorem~\ref{theorem:comparison_box_with_t_times_first_eigenvalue}}
\label{section:asymptotics_large_box}

In this section we prove Theorem~\ref{theorem:comparison_box_with_t_times_first_eigenvalue}. 
We separate the proof in two parts by showing that almost surely $\limsup_{t\to \infty} \frac{U_{t}(t)}{t \blambda_{1, t}} \le 1$ (Lemma~\ref{lemma:limsup_log_mass_box_L_t}) and $\liminf_{ t\in \Q, t\to \infty} \frac{U_{t}(t)}{t \blambda_{1, t}} \ge 1$ (Lemma~\ref{lemma:liminf_bound}).

We use the notation $u_L^x = u_L^{\delta_x}$ and $U_L^x(t)= \int_{Q_L} u_L^x(t,y) \dd y$. 
We let $\Gamma_t^L(x,\cdot)$ be the transition probability kernel of $X_t$ under $\Q_{L}^x$. 
We will rely on the following estimates which follow from Lemma~\ref{lemma:a_epsilon_bound_of_M_L_eps}, Lemma~\ref{lemma:escape_probability_of_box_L_with_drift_Z_and_Y_L} and 
\cite[Theorem 1.1 and Corollary 1.2]{PevZ}: There exists a random variable $C$ (depending on $\xi$, not on $X$) with values $>1$  such that for $L\in \N$ with $L>e$, for $s,t\in [0,\infty)$ with $s<t$ and $t-s>1$ 
\begin{align}
\label{eqn:cD_bound}
e^{- C (t-s) (\log L)^5} \1_{[X_{[s,t]}\subset Q_L]}
\le \cD_L(s , t) \1_{[X_{[s,t]}\subset Q_L]} \le e^{C (t-s) (\log L)^5},
\end{align}
as well as for $L,r\in \N$ with $L>r>e$ and $t\ge 1$, 
\begin{calc}
 for $c>1$
\begin{align*}
 \tfrac{1}{2\pi ct} e^{-\frac{L^2}{ t}} 
\le \tfrac{1}{2\pi ct} e^{-\frac{L^2}{2ct}} 
\le p(ct,x-y)  
\le \tfrac{1}{2\pi ct}, 
\end{align*} 
and as 
\begin{align*}
e^{- t( \log L)^4}
\le \tfrac{1}{t}\le e^{ t( \log L)^4}, 
\end{align*}
\end{calc}
\begin{align}
\label{eqn:Gamma_bound_below_box_L}
\Gamma_t^L(x,y) & \ge \tfrac{1}{C} e^{-C t (\log L)^5 - \frac{r^2}{Ct}}, \qquad x,y \in Q_r, \\
\label{eqn:Gamma_bound_above_box_L}
 \Gamma_t^L(x,y) 
& \le C  e^{C t (\log L)^5}, \qquad x,y \in \R^2, \\
\label{eqn:escape_box_L}
\Q_L^0 [ X_{[0,t]}\not \subset Q_L] 
& \le  C e^{Ct (\log L)^5 - \frac{L^2}{Ct}}. 
\end{align}

\begin{obs}
\label{obs:total_mass_as_sol_with_1_as_initial_condition}
We use the notation as in Lemma~\ref{lemma:solution_map_box_continuous_in_inintial} and write $u_\infty^\phi$ for the solution of the PAM with white-noise potential  on $\R^2$ with initial condition $\phi$ as in \cite{HaLa15}. 

By Theorem~\ref{theorem:spectral_rep_sol_box_delta_initial},  Lemma~\ref{lemma:spectral_rep_sol_box_smooth_initial} and Lemma~\ref{lemma:feynman_kac_box_smooth} we have for all $L\in (0,\infty]$, $\epsilon >0$ and $x,y\in \R^2$ 
\begin{align}
\label{eqn:u_x_y_change}
 u_{L,\epsilon}^{x}(t,y) 
&  = u_{L,\epsilon}^{y}(t,x),  \\
\notag 
\int_{\R^2} u^{x}_{L,\epsilon}(t,y) \dd y
&= \langle u^{x}_{L,\epsilon}(t,\cdot) , \1 \rangle_{L^2(\R^2)} 
 = u_{L,\epsilon}^\1(t,x)  \\
\label{eqn:total_mass_in_terms_of_one_as_initial_condition}
& = \E^x \left[ \exp \left( \int_0^t (\xi_\epsilon(X_s) -  c_\epsilon ) \dd s \right) \1_{[X_{[0,t]} \subset Q_L]}  \right].
\end{align}
By the continuity in $\epsilon$, see Lemma~\ref{lemma:solution_map_box_continuous_in_inintial}, for $L >0$
\begin{align}
\label{eqn:total_mass_box_as_initial_condition_1}
U_{L}^x(t) := \int_{Q_L} u_L^{x}(t,y) \dd y = u_L^\1(t,x) 
= \E_{\Q_L^x}[ \cD_L(0,t) \1_{[X_{[0,t]} \subset Q_L]} ] . 
\end{align}
We also write ``$U_L(t)$'' for ``$U_L^0(t)$''. 
\end{obs}

Before we turn to the upper bound, we prove the following lemma that will be used for both the upper and the lower bound. 

\begin{lemma}
\label{lemma:u_r_y_bounded_from_above_by_factor_times_total_mass}
For all $r\in \N$ with $r>e$, $t>1$, $\delta \in (1,t)$ and $x,y\in Q_r$ we have with $C$ as above
\begin{align}
\label{eqn:bound_u_r_x_t_x}
u_r^y (t,x) = u_r^x (t,y) 
\le 
C e^{2C \delta (\log r)^5} 
U_{r}^x(t-\delta). 
\end{align}
\end{lemma}

\begin{proof}
 The first equality in \eqref{eqn:bound_u_r_x_t_x} follows from \eqref{eqn:u_x_y_change}.  
We have $u_r^y (t,x) = \lim_{\epsilon \downarrow 0} u_r^{\psi_\epsilon^y}(t,x)$ by Lemma~\ref{lemma:solution_map_box_continuous_in_inintial}, where $\psi_\epsilon^y(z) = \psi_\epsilon(z-y)$. By definition we have
$\cD_r(0,t) =\cD_r(0,t-\delta) \cD_r(t-\delta,t)$. 
By using subsequently \eqref{eqn:cD_bound} and the tower-property; the definition of $\Gamma_t^r$,
\eqref{eqn:Gamma_bound_above_box_L} and $\int \psi_\epsilon^y=1$,
and finally \eqref{eqn:total_mass_box_as_initial_condition_1} we obtain 
\begin{align*}
 u_r^{y}(t,x) 
\cand \begin{calc}
=\lim_{\epsilon \downarrow 0} u_r^{\psi_\epsilon^y}(t,x)
= \lim_{\epsilon\downarrow 0} \E_{\Q_{r,x}} \left[ \cD_r(0,t)  \psi_\epsilon^y(X_t) \1_{[X_{[0,t]}\subset Q_r]} \right]
\end{calc} \cnewline
& \le \lim_{\epsilon \downarrow 0} e^{C \delta (\log r)^5} \E_{\Q_{r,x}} \left[ \cD_r(0,t-\delta)  \1_{[X_{[0,t-\delta]}\subset Q_r]} \E_{\Q_{r,x}}[\psi_\epsilon^y(X_{t}) | X_{t-\delta} ] \right]\\
& = \lim_{\epsilon \downarrow 0} e^{C \delta (\log r)^5} \E_{\Q_{r,x}} \left[ \cD_r(0,t-\delta)  \1_{[X_{[0,t-\delta]}\subset Q_r]} \int_{\R^2} \Gamma_\delta^r( X_{t-\delta},z) \psi_\epsilon^y(z) \dd z \right]\\
& \le C e^{C \delta (\log r)^5 + C \delta (\log r)^5} \E_{\Q_{r,x}} \left[ \cD_r(0,t-\delta)  \1_{[X_{[0,t-\delta]}\subset Q_r]} \right] \\
& \le C e^{2 C \delta (\log r)^5} U_{r}^x(t-\delta). 
\end{align*}
\end{proof}

\subsection{Upper bound}

We will use the spectral representation given in Theorem~\ref{theorem:spectral_rep_sol_box_delta_initial}. 
Observe that for $L,t>0$ and $\phi \in C(Q_L)$, we have by Lemma~\ref{lemma:spectral_rep_sol_box_smooth_initial} 
\begin{align}
\notag \int_{Q_L} u^\phi_L(t,x) \dd x 
& = 
\sum_{n\in\N} e^{t\blambda_{n,L}} \langle v_{n,L} ,\phi \rangle_{L^2} \langle v_{n,L} ,\1_{Q_L} \rangle_{L^2} \\
\notag 
& \le 
\Big( \sum_{n\in\N} e^{2t\blambda_{n,L}} \langle v_{n,L} ,\phi \rangle_{L^2}^2 \Big)^{\frac12} \| \1_{Q_L} \|_{L^2} \\\
& \le 
\label{eqn:upper_estimate_smooth_initial_total_mass}
e^{t\blambda_{1,L}} 
\| \phi \|_{L^2}
\| \1_{Q_L} \|_{L^2}. 
\end{align}
As our initial condition $\phi = \delta_0$ is not in $L^2$, we use that by the Chapman-Kolmogorov equation $u_L^0(t,x) = u_L^\phi(t-3,x)$ for $t>3$, with $\phi = u_L^0(3,x)$ and show that $u_L^0(3,x)$ is in $L^2$ and that its $L^2$-norm can be bounded as follows.  

\begin{lemma}
\label{lemma:q_th_moment_of_solution_on_box}
Let $q \in [1,\infty]$, $L\in\N$, $L>e$, $y \in Q_L$ and $t>2$. 
Then $u_L^y(t,\cdot) \in L^q(Q_L)$ and $\|u_L^y(t,\cdot)\|_{L^q} \le C e^{3Ct (\log L)^5} $.
\end{lemma}
\begin{proof}
We treat the case $q<\infty$, the case $q=\infty$ is similar but slightly easier.
By Lemma~\ref{lemma:u_r_y_bounded_from_above_by_factor_times_total_mass} we have for $t>2$ and $\delta = \frac{t}{2}$, 
\begin{align*}
\|u_L^y(t,\cdot)\|_{L^q} \le 
C e^{C t (\log L)^5} 
\Big( \int_{Q_L} |U_L^x(\tfrac{t}{2})|^q \dd x \Big)^{\frac1q}. 
\end{align*}
As $U_L^x(\frac{t}{2}) = u_L^\1(\frac{t}{2},x) \le e^{C t (\log L)^5}$ by \eqref{eqn:total_mass_box_as_initial_condition_1} and \eqref{eqn:cD_bound}, and as $L^{\frac{2}{q}}\le L^2 \le e^{C t (\log L)^5}$ for $t>2$, the desired estimate follows. 
\end{proof}

\begin{lemma}
\label{lemma:limsup_log_mass_box_L_t}
Let $L_t = t (\log t)^5$ or $L_t =t$. Then almost surely 
\begin{align}
\label{eqn:limsup_log_mass_box_L_t}
\limsup_{t\rightarrow \infty} \frac{\log U_{L_t}(t)}{ t \blambda_{1,L_t}  } \le 1. 
\end{align}
\end{lemma}
\begin{proof}
Let $\fL_t= \lceil L_t \rceil$. Then $U_{L_t} \le U_{\fL_t}$ and $\blambda_{1,L_t} \sim \blambda_{1,\fL_t}$. 
For $t$ large enough we have by \eqref{eqn:upper_estimate_smooth_initial_total_mass} and Lemma~\ref{lemma:q_th_moment_of_solution_on_box} (as $\|\1_{Q_L}\|_{L^2} = L \le e^{C (\log L)^5}$):
\begin{align*}
U_{\fL_t}(t) 
\le e^{(t-3) \blambda_{1,\fL_t}} C e^{10 C (\log \fL_t)^5},
\end{align*}
so the claim easily follows. 
\end{proof}

\subsection{Lower bound}

By the spectral representation given in Theorem~\ref{theorem:spectral_rep_sol_box_delta_initial} we obtain 
\begin{align}
\label{eqn:e_t_lambda_bound_from_below_by_int}
e^{t \blambda_{1,L}} 
\le \sum_{n\in\N} e^{t \blambda_{n,L}} \int_{Q_L} |v_{n,L}(x)|^2 \dd x
= \int_{Q_L} u_L^x(t,x) \dd x. 
\end{align}
So we would already have established the lower bound if ``$u_L^x$'' on the right-hand side was replaced by ``$u_L^0$''. The idea is therefore to use \eqref{eqn:e_t_lambda_bound_from_below_by_int} with ``$t-\delta$'' instead of ``$t$'', for $\delta$ chosen appropriately (depending on $t$), and to combine this with lower bounds for the transition density $\Gamma^{L}_t$ in order to obtain a lower bound for $U_{L}(t)$. 

Inequality \eqref{eqn:cD_bound} together with the Markov property gives for  $\delta \in (1,t)$
\begin{align}
\label{eqn:U_L_ge_e_factor_times_expectation_of_delta_to_t}
U_L(t) 
\ge e^{- C \delta (\log L)^5} 
\E_{\Q_L^0} \Big[ \1_{[X_{[0,\delta]} \subset Q_L]} U_{L}^{X_\delta}(t-\delta) \Big].
\end{align}
Now we take a smaller box inside the box of size $L$, i.e., let $r\in (0,L]$. 
Then 
\begin{align}
\notag & \E_{\Q_L^0} \Big[ \1_{[X_{[0,\delta]} \subset Q_L]} U_L^{X_\delta}(t-\delta) \Big]
 \ge \E_{\Q_L^0} \Big[ \1_{[X_{[0,\delta]}  \subset Q_L]} \1_{[X_\delta \in Q_r]} U_r^{X_\delta}(t-\delta) \Big] \\
\label{eqn:expectation_X_stays_in_box_till_delta_then_rest_mass_from_below}
& = \E_{\Q_L^0} \Big[ \1_{[X_\delta \in Q_r]} U_r^{X_\delta}(t-\delta) \Big] -
\E_{\Q_L^0} \Big[ \1_{[X_{[0,\delta]}  \not \subset Q_L]} \1_{[X_\delta \in Q_r]} U_r^{X_\delta}(t-\delta) \Big].
\end{align}
The second term on the right-hand side of \eqref{eqn:expectation_X_stays_in_box_till_delta_then_rest_mass_from_below} can be estimated as follows, using \eqref{eqn:escape_box_L}, \eqref{eqn:total_mass_box_as_initial_condition_1} and \eqref{eqn:cD_bound} 
\begin{align*}
& \E_{\Q_L^0} \Big[ \1_{[X_{[0,\delta]}  \not \subset Q_L]} \1_{[X_\delta \in Q_r]} U_r^{X_\delta}(t-\delta) \Big]
 \le 
\P_{\Q_L^0} [X_{[0,\delta]}  \not \subset Q_L] \sup_{x\in Q_r} U_r^x(t-\delta)  \\
& \le C e^{C \delta (\log L)^5 - \frac{L^2}{C\delta}}
e^{C (t-\delta) (\log r)^5 } 
\le  C e^{C t (\log L)^5 - \frac{L^2}{C\delta}}.
\end{align*}
By the heat kernel bound \eqref{eqn:Gamma_bound_below_box_L}, we can estimate the first term on the right-hand side of \eqref{eqn:expectation_X_stays_in_box_till_delta_then_rest_mass_from_below} by
\begin{align*}
 \E_{\Q_L^0} \Big[ \1_{[X_\delta \in Q_r]} U_r^{X_\delta}(t-\delta) \Big]
& = \int_{Q_r} \Gamma^L_\delta(0,x) U_r^x(t-\delta) \dd x \\
& \ge \frac{1}{C}  e^{-C \delta (\log L)^5 - \frac{r^2}{C\delta}} 
\int_{Q_r} U_r^x(t-\delta) \dd x.
\end{align*}

Now we use the lower estimate for $U_{t}^x(t-\delta)$ from Lemma~\ref{lemma:u_r_y_bounded_from_above_by_factor_times_total_mass} which allows us to use \eqref{eqn:e_t_lambda_bound_from_below_by_int}. 
We sum up all the estimates in the following lemma. 
\begin{calc}
\begin{align*}
& U_L(t) \\
& \ge e^{-C\delta (\log L)^p} \left( 
 \E_{\Q_L^0} \Big[ \1_{[X_\delta \in Q_r]} U_r^{X_\delta}(t-\delta) \Big] -
\E_{\Q_L^0} \Big[ \1_{[X_{[0,\delta]}  \not \subset Q_L]} \1_{[X_\delta \in Q_r]} U_r^{X_\delta}(t-\delta) \Big]
 \right) \\
& \ge \frac{1}{C}  e^{-C \delta (\log L)^4 - \frac{r^2}{C\delta}} 
\int_{Q_r} U_r^x(t-\delta) \dd x
-  C e^{C t (\log L)^5 - \frac{L^2}{C\delta}}  \\
& \ge 
\frac{1}{C^2}  e^{-3C \delta (\log L)^5 - \frac{r^2}{C\delta}} 
\int_{Q_r} u_{r}^x(t-\delta,x) \dd x
-  C e^{C t (\log L)^p - \frac{L^2}{C\delta}} . 
\end{align*}
\end{calc}

\begin{lemma}
\label{lemma:lower_bound_U_L_t_with_r_and_delta}
For $L,r\in\N$ with $L\ge r>e$ and $t>\delta>1$ (with $C$ as above) we have
\begin{align*}
U_L(t) \ge
\frac{1}{C^2}  e^{-3C \delta (\log L)^5 - \frac{r^2}{C\delta}} e^{(t-\delta)\blambda_{1,r}}
-  C  e^{C t (\log L)^5 - \frac{L^2}{C\delta}}.
\end{align*}
\end{lemma}

Now we tune $L$, $r$ and $\delta$, i.e., we choose them depending on $t$ in such a way that the lower bound for the total mass follows.

\begin{lemma}
\label{lemma:liminf_bound}
Let $b\in (\frac12,1]$ and $L_t \ge t^b$ for all $t\ge 0$ and such that $\log L_t \sim b \log t$. 
Let 
$\I \subset (e,\infty)$ be a countable unbounded set. Then almost surely
\begin{align*}
\liminf_{t \in \I, t\rightarrow \infty} \frac{\log U_{L_t}(t)}{t \blambda_{1,L_t} } \ge  1.
\end{align*}
\end{lemma}
\begin{proof}
Let $\fL_t = \lfloor L_t \rfloor$. Then $U_{L_t} \ge U_{\fL_t}$. 
Let $b' \in (0,b)$, $a\in (2b'-1,2b-1)$, $\delta_t = t^{a}$ and $r_t = \lfloor t^{b'} \rfloor$. 
Write 
\begin{align*}
v_t = \frac{1}{C^2}  e^{-3C \delta_t (\log \fL_t)^5 - \frac{r_t^2}{C\delta_t}} e^{(t-\delta_t)\blambda_{1,r_t}}, 
\qquad 
w_t =  C e^{C t (\log \fL_t)^5 - \frac{L_t^2}{C\delta_t}}. 
\end{align*}
By Lemma~\ref{lemma:lower_bound_U_L_t_with_r_and_delta} we have $\log U_{L_t}(t) \ge \log U_{\fL_t}(t) \ge \log( v_t - w_t) = \log v_t + \log (1- \frac{w_t}{v_t})$. 
Observe that, as $\delta_t <t$, 
\begin{align*}
\frac{w_t}{v_t} \le C^3 \exp \left( 
4 C t (\log \fL_t)^5 - \frac{\fL_t^2}{C\delta_t} 
+ \frac{r_t^2}{C\delta_t} - (t- \delta_t) \blambda_{1,r_t}
\right). 
\end{align*}
By Theorem~\ref{theorem:convergence_smooth_eigenvalues_and_asymptotics} we have $\blambda_{1,r_t} > 0$ for large $t$, and as $\fL_t^2 \delta_t^{-1} \ge t^{2b-a}$ and $r_t^2/\delta_t = t^{2b'-a}$ and $2b-a>\max\{1,2b'-a\}$, we deduce that $\frac{w_t}{v_t} \rightarrow 0$ and thus $\log (1- \frac{w_t}{v_t}) \rightarrow 0$. 
Therefore, we obtain from Theorem~\ref{theorem:convergence_smooth_eigenvalues_and_asymptotics}:
\begin{align*}
	\liminf_{t\in \I, t\rightarrow \infty} \frac{\log U_{L_t}(t)}{t \blambda_{1,L_t}} 
	\ge \liminf_{t\in \I, t\rightarrow \infty} \frac{\log v_t}{t \blambda_{1,r_t}}  \frac{\blambda_{1,r_t}}{\blambda_{1,L_t}} 
	= \liminf_{t\in \I, t\rightarrow \infty} \frac{\blambda_{1,r_t}}{\blambda_{1,L_t}} = \liminf_{t\in \I, t\rightarrow \infty} \frac{\chi \log r_t}{\chi \log L_t} \ge \frac{b'}{b} ,
\end{align*}
where we used that $\log r_t \sim b' \log t$ and $\log L_t \sim b \log t$. 
\begin{calc}
\begin{align*}
\liminf_{t\in \I, t\rightarrow \infty} \frac{\log v_t}{t \blambda_{1,r_t}}  
& = \liminf_{t\in \I, t\rightarrow \infty} \frac{-3C \delta_t (\log L_t)^5 - \frac{r_t^2}{C\delta_t}}{t \blambda_{1,r_t}} + \frac{t-\delta_t}{t} \\
& = 1- \liminf_{t\in \I, t\rightarrow \infty} \frac{r_t^2}{C\delta_t t \chi \log r_t} \\
& = 1- \liminf_{t\in \I, t\rightarrow \infty} \frac{t^{2b' -1-a}}{C  \chi \log r_t} =1. 
\end{align*}
\end{calc}
Since $b'\in (0,b)$ was arbitrary, the claim follows.
\end{proof}

Now Lemma~\ref{lemma:limsup_log_mass_box_L_t} and Lemma~\ref{lemma:liminf_bound} imply Theorem~\ref{theorem:comparison_box_with_t_times_first_eigenvalue}.

\section{Splitting the PAM into boxes: Proof of Proposition~\ref{proposition:comparison_total_mass_full_space_and_box}}
\label{section:separating_space_in_boxes}

In this section, we prove Proposition~\ref{proposition:U_t_and_U_L_t_equivalence}, that is, the fact that the total mass of the solution is well approached by the total mass in a sufficiently large box. As we indicated in Section~\ref{sec-purpose}, we will be following a standard strategy that decomposes $\R^d$ into many large boxes and estimates the contribution from each box to the Feynman--Kac formula. 
We will be using the Feynman--Kac-type representation that we derived in Theorem~\ref{theorem:feynman_kac_type_representation_for_PAM} and the estimates of Lemma~\ref{lemma:escape_probability_of_box_L_with_drift_Z_and_Y_L}.

Let $U(t) = \int_{\R^2} u_\infty(t,x) \dd x$ for $t\in [0,\infty)$, where $u_\infty$ is the solution to the PAM on $\R^2$ with $u_\infty(0,\cdot) = \delta_0$. 
We will choose $(L_t)_{t\in [0,\infty)}$ in $\N$ later, in such a way that $L_t \uparrow \infty$. We assume $t$ is large enough such that $L_t >1$. 
We define the events 
\begin{align*}
A_0 = \left[ X_{[0,t]} \subset Q_{L_t} \right], 
\qquad 
A_k = \left[ X_{[0,t]} \not \subset Q_{L_t^k},  X_{[0,t]} \subset Q_{L_t^{k+1}} \right] \quad \mbox{ for } k\in \N. 
\end{align*}
As $X$ has continuous paths, we have $\P^0( \bigcup_{k\in\N_0} A_k) =1$. 
 We consider the setting as in Section~\ref{subsection:feynman_kac_for_pam_with_white_noies}. 
We define for $k\in \N_0$ and $\epsilon \ge 0$
\begin{align*}
\fU_{k,\epsilon}(t) := \E_{\Q_{L_t^{k+1},\epsilon}^0 } \left[ \cD_{L_t,\epsilon}(0,t)   \1_{A_k}\right].
\end{align*}
We will abbreviate $\fU_k = \fU_{k,0}$. Observe that $U_{L_t}(t) = \fU_{0}(t)$. 
In Lemma~\ref{lemma:total_mass_in_sum} below we prove that $U(t) = \sum_{k\in \N_0} \fU_k(t)$, and in Lemma~\ref{lemma:bound_U_k_t} we derive bounds on $\fU_k(t)$. Combining these, we will then prove Proposition~\ref{proposition:U_t_and_U_L_t_equivalence}, which in turn implies Proposition~\ref{proposition:comparison_total_mass_full_space_and_box}.

\begin{obs}
By the continuity in $\epsilon$, see \cite[Theorem 4.1]{HaLa15}, $U(t) = u_\infty^\1(t,0)$ (almost surely), where $u_\infty^\1$ is the solution to the PAM on $\R^2$ with $u_\infty^\1(0,\cdot) = \1$. 
Observe that by
Lemma~\ref{lemma:feynman_kac_box_smooth}, 
 Lemma~\ref{lemma:feynman_kac_like_representation_on_box_with_epsilon} and \eqref{eqn:total_mass_box_as_initial_condition_1} we have for $L>0$, and $\epsilon>0$
\begin{align}
\label{eqn:total_mass_epsilon_equals_sum_k_fk_U}
u_{\infty,\epsilon}^\1(t,0) = \sum_{k\in \N_0} \fU_{k,\epsilon}(t).
\end{align}
\end{obs}

\begin{lemma}
\label{lemma:bound_U_k_t}
Let 
$L_t = \lfloor t (\log t)^5 \rfloor$ and let $a_\epsilon$ be as in Lemma~\ref{lemma:a_epsilon_bound_of_M_L_eps}. There exist $C>1$ and $T > 0$ such that for all $t \ge T$, $k \in \N$ and $\epsilon \ge 0$
\begin{align*}
\fU_{k,\epsilon}(t)
& \le 
\begin{cases}
C\exp \Big( t (\log t)^5  \Big[C a_\epsilon^5  (k+1)^5   -  \frac{ t^{2k-2} (\log t)^{5(2k-1)} }{ C  4^k   } \Big] \Big),
& k \ge 1, \\
C\exp \Big(  (k+1)^5 t (\log t)^5 \Big[C a_\epsilon^5    -  \frac{k t}{ C } \Big] \Big),
& k \ge 2. 
\end{cases}
\end{align*}
\end{lemma}

\begin{proof}
By Lemma~\ref{lemma:escape_probability_of_box_L_with_drift_Z_and_Y_L} there exists a $C>1$ such that for all $t\ge 1$ with $L_t >e$, for all $k\in \N$ and for all $\epsilon \ge 0$ 
\begin{align}
\label{eqn:fU_k_eps_bound}
\fU_{k,\epsilon}(t) \le C \exp \Big( C a_\epsilon^5 t (k+1)^5 (\log L_t)^5 -  \frac{L_t^{2k}}{C t} \Big) .
\end{align}
Let $T > 0$ be such that $\log t > 5\log \log t$ and $\lfloor t (\log t)^5\rfloor \ge \tfrac12 t (\log t)^5$ for $t \ge T$. Then $\log L_t \le \log t + 5 \log \log t \le 2 \log t$ for all $t \ge T$, and therefore the first bound follows from~\eqref{eqn:fU_k_eps_bound} (with a new $C>0$). 
\begin{calc}
\begin{align*}
C a_\epsilon^5 t (k+1)^5 (\log L_t)^5 -  \frac{L_t^{2k}}{C t}
\le 
C a_\epsilon^5 t (k+1)^5 2^5 (\log t)^5 -  \frac{t^{2k-1} (\log t)^{10k}}{2^{2k} C t}
\end{align*}
\end{calc}
The second bound follows from the first one by choosing $T>0$ large enough so that 
\begin{align*}
T^{2k -3} (\log T)^{5(2k-1)} \ge  4^k  k (k+1)^5
\end{align*}
 for all $k \ge 2$ (and then of course the same inequality holds for all $t \ge T$).
\end{proof}

These bounds allow us to derive the series expansion of the total mass $U(t)$:

\begin{lemma}
\label{lemma:total_mass_in_sum}
Let $L_t = \lfloor t (\log t)^5 \rfloor$.
There exists a $T>0$ such that for all $t\ge T$,  we have almost surely $U(t) = \sum_{k\in \N_0} \fU_k(t)$.
\end{lemma}
\begin{proof}
Since each term $\fU_{k}(t)$ is positive the series $\sum_{k\in \N_0} \fU_k(t)$ converges almost surely to some $V(t)$ with values in $[0,\infty]$. So if we can show that the series converges in probability to $U(t)$, then almost surely $V(t) = U(t)$ and the claimed identity holds.

By \cite[Theorem 1.4]{HaLa15} we know that $u^\1_{\infty,\epsilon}(t,0) \xrightarrow{\epsilon\downarrow 0} u^\1_{\infty,0}(t,0) = U(t)$ in probability. By Theorem~\ref{theorem:feynman_kac_like_representation_for_limit} we have $\fU_{k,\epsilon}(t) \xrightarrow{\epsilon\downarrow 0} \fU_k(t)$ in probability, for all $k \in \N_0$.
Moreover, $u^\1_{\infty,\epsilon}(t,0) = \sum_{k \in \N_0} \fU_{k,\epsilon}(t)$ for $\epsilon>0$, see \eqref{eqn:total_mass_epsilon_equals_sum_k_fk_U}. 
Therefore, we have for $K\in \N$, $\delta \in (0,1)$ and $t > 0$ such that $L_t > e$:
\begin{align*}
&	P \Big( \Big| U(t) - \sum_{k =0 }^K \fU_k(t) \Big| > \delta \Big) 
 \le \limsup_{\epsilon \downarrow 0} P\left( \Big| u_\infty^\1(t,0) - u^\1_{\infty,\epsilon}(t,0) \right| > \tfrac{\delta}{3} \Big) \\
& \quad  + \limsup_{\epsilon \downarrow 0} P \Big( \Big| u^\1_{\infty,\epsilon}(t,0) - \sum_{k=0}^K \fU_{k,\epsilon}(t) \Big| > \tfrac{\delta}{3} \Big)  + \limsup_{\epsilon \downarrow 0} \sum_{k=0}^K P\left( \left| \fU_{k,\epsilon}(t) - \fU_{k}(t) \right| > \tfrac{\delta}{3K} \right)  \\
& \le  \limsup_{\epsilon \downarrow 0} 
P \Big( \Big| \sum_{k=K+1}^\infty \fU_{k,\epsilon}(t) \Big| > \tfrac{\delta}{3}, \ a_\epsilon^5 \le \frac{KT}{2C^2} \Big)
+ \limsup_{\epsilon \downarrow 0} P \Big( a_\epsilon^5 > \frac{KT}{2C^2} \Big).  
\end{align*}
Let $T > e$ and $C$ be as in Lemma~\ref{lemma:bound_U_k_t} and let $t\ge T$ and $K\ge 1$.  On the event $[a_\epsilon^5 \le \frac{KT}{2C^2}]$ we have by the second estimate of Lemma~\ref{lemma:bound_U_k_t}  
\begin{align*}
\sum_{k =K+1 }^\infty 
\fU_{k,\epsilon}(t)
& \le 
\sum_{k =K+1 }^\infty C\exp \Big(  (k+1)^5 t (\log t)^5 \Big[C a_\epsilon^5    -  \frac{K t}{ C } \Big] \Big)\\
&  \le C \sum_{k = K+1}^\infty \exp \Big( - (k+1)^5 T (\log T)^5  \frac{K T}{2 C } \Big),
\end{align*}
\begin{calc}
Let $f : (0,1) \rightarrow (0,\infty)$ be given by $f(r) = \frac{r}{1-r}$. 
Then $f$ is increasing on $(0,\frac12)$ and the maximum of $f$ is attained at $\frac12$ and equal to $1$. 
There exists an $r_\delta$ such that $f(r_\delta) = \frac{\delta}{3}$. 
As $\delta \in (0,1)$ we have $f(r) > \frac{\delta}{3}$ only if $r > r_\delta$. 
\end{calc}
and the right-hand side vanishes as $K\to \infty$. 
Therefore 
\begin{align*}
\limsup_{\epsilon \downarrow 0} P \Big( \Big| \sum_{k=K+1}^\infty \fU_{k,\epsilon}(t) \Big| > \tfrac{\delta}{3}, \ a_\epsilon^5 \le \frac{KT}{2C^2} \Big) \xrightarrow{\epsilon \downarrow 0} 0. 
\end{align*}
By Markov's inequality we have 
\begin{align*}
\sup_{\epsilon\ge 0 } P \left(a_\epsilon^5 > \frac{KT}{2C^2} \right)
\le  \exp\left( -h \left( \frac{KT}{2C^2} \right)^{\frac15} \right) \sup_{\epsilon\ge 0} E[e^{h a_\epsilon} ],
\end{align*}
so that by Lemma~\ref{lemma:a_epsilon_bound_of_M_L_eps} it follows that 
$\sup_{\epsilon\ge 0} P(a_\epsilon^5 > \frac{KT}{2C^2})$ converges to zero as $K \to \infty$. 
This completes the proof.
\end{proof}

\begin{obs}
\label{obs:limit_log_sum_to_max}
In the proof of the next proposition we use the following observation. 
Let $K\in \N$, $q_{t,k} \in (0,\infty)$ and $R_t \in (0,\infty)$ for $t\in (0,\infty)$ and $k \in \{0,\dots, K\}$. Suppose $R_t \xrightarrow{t\rightarrow \infty} \infty$. 
As 
\begin{align*}
\max_{k\in \{0,\dots, K\}} \frac{\log  q_{t,k}}{R_t} \le \frac{\log \sum_{k=0}^K q_{t,k}}{R_t} \le \frac{\log K}{R_t} + \max_{k\in \{0,\dots, K\}} \frac{\log  q_{t,k}}{R_t},
\end{align*}
we have
\begin{align*}
\lim_{t\rightarrow \infty} \left| \frac{\log \sum_{k=0}^K q_{t,k}}{R_t} - \max_{k\in \{0,\dots, K\}} \frac{\log  q_{t,k}}{R_t} \right| = 0.
\end{align*}
\end{obs}

\begin{proposition}
\label{proposition:U_t_and_U_L_t_equivalence}
For  $L_t = t (\log t)^5$ we have almost surely 
\begin{align*}
\lim_{t\in \Q, t\rightarrow \infty} 
\left|
\frac{\log U(t)}{ t \log L_t} - \frac{\log U_{L_t} (t)}{ t \log L_t} 
\right| =0 .
\end{align*}
\end{proposition}

\begin{proof}
By \eqref{eqn:total_mass_box_as_initial_condition_1} for all $\epsilon>0$ the function $L \mapsto u_{L,\epsilon}^\1(t,0)$ is increasing. 
This implies, by Theorem~\ref{theorem:feynman_kac_like_representation_for_limit}, that also $L \mapsto U_{L}(t)$ is increasing. 
Therefore, we may consider $L_t = \lfloor t (\log t)^5 \rfloor$ instead, so that $U_{L_t}(t) = \fU_0(t)$. We will use $U(t) = \sum_{k \in \N_0} \fU_k(t)$ (see Lemma~\ref{lemma:total_mass_in_sum}).

Let $T$ and $C$ be as in Lemma~\ref{lemma:bound_U_k_t}. Then we have for $k = 1$ almost surely
\begin{align*}
\lim_{t \to \infty} \frac{\log \fU_1(t) }{ t \log L_t  } 
\le \lim_{t\to \infty} \frac{t (\log t )^5}{t \log L_t}  \Big[C (2a)^5  2^5   -  \frac{(\log t)^5 }{ C } \Big] 
\xrightarrow{t\rightarrow \infty} = - \infty. 
\end{align*}
Moreover, if $t \ge T$ is large enough so that $C (2a)^5  - \frac{2 t}{C} \le -1$ and $\exp(-t(\log t)^5) \le \frac12$, then the second bound of Lemma~\ref{lemma:bound_U_k_t} gives
\begin{align*}
\sum_{k=2}^\infty \fU_k(t) 
&\lesssim
\sum_{k=2}^\infty   \exp \Big( (k+1)^5 t (\log t)^5  \Big[C (2a)^5  - \frac{k t}{C} \Big] \Big) \\
& \le
\sum_{k=2}^\infty  \exp ( - t (\log t)^5 )^{k+1} 
\begin{calc}
= \sum_{k=3}^\infty r^k = \frac{r^3}{1-r} < \frac{2^{-3}}{1-\frac12} 
\end{calc}
\le 1,
\end{align*}
and therefore $\frac{\log (\sum_{k=2}^\infty \fU_k(t) ) }{ t \log L_t  }  \to -\infty$ as $t \to \infty$. 
Now it suffices to apply \ref{obs:limit_log_sum_to_max} with $R_t = t \log L_t$, $q_{t,k} = \fU_k(t)$ for $k =0,1$, and $q_{t,2} = \sum_{k=2}^\infty \fU_k(t)$. 
  As 
\begin{align*}
\liminf_{t\in \Q, t \to \infty} \frac{\log U_{L_t}(t)}{t \log t} > -\infty
\end{align*}
by Lemma~\ref{lemma:liminf_bound}, this completes the proof. 
 
\end{proof}

\begin{proof}[Proof of Proposition~\ref{proposition:comparison_total_mass_full_space_and_box}]
By Proposition~\ref{proposition:U_t_and_U_L_t_equivalence} we have $\log U(t) \sim \log U_{L_t}(t)$. 
By Lemma~\ref{lemma:limsup_log_mass_box_L_t} and Lemma~\ref{lemma:liminf_bound} we have $\log U_{L_t}(t) \sim t \lambda_{1,L_t}$ and $\log U_{t}(t) \sim t \lambda_{1,t}$. But $\lambda_{1,L_t} \sim \lambda_{1,t}$ as $ t \sim L_t$ by Theorem~\ref{theorem:convergence_smooth_eigenvalues_and_asymptotics}. 
\end{proof}

\section{Asymptotics of the supremum and infimum of the PAM}
\label{section:asymptotics_of_sup_and_inf_boxes}

In this section we prove Theorem~\ref{theorem:total_mass_asymptotics}~\ref{item:infimum_and_supremum_asymptotics} by proving that for all $a \in (0, 1)$: 
\begin{align*} 
\lim_{t\in \Q, t \rightarrow \infty} \frac{\log U (t)}{t \log t} 
& \le \lim_{t\in \Q, t \rightarrow \infty} \frac{\log \left( \inf_{x \in Q_{t^a}} u (t, x) \right)}{t \log t} \\
& \le \lim_{t\in \Q, t \rightarrow \infty} \frac{\log (\sup_{x \in \mathbb{R}^2} u (t, x))}{t \log t}
\le \lim_{t\in \Q, t \rightarrow \infty} \frac{\log U (t)}{t \log t} .
\end{align*}
To prove this, we first have to show that $u(t,x)$ is comparable to $u_{L_t}(t,x)$, and for that purpose we need the following lemma, in which $A_k$ is as in Section~\ref{section:separating_space_in_boxes}, i.e.
\begin{align*} 
A_k = [ X_{[0, t]} \not \subset Q_{L_t^k}, X_{[0, t]} \subset Q_{L_t^{k + 1}} ] . 
\end{align*}
  
\begin{obs}
Let us write $\Gamma_t^{L,\epsilon}(x,\cdot)$ for the transition probability kernel of $X_t$ under $\Q_{L,\epsilon}^x$. Let $c>1$. Then the heat kernel bound in \cite[Theorem 1.1]{PevZ} and Lemma~\ref{lemma:escape_probability_of_box_L_with_drift_Z_and_Y_L} imply the existence of a $C>1$ and a $\kappa \in (0,1)$ such that for all $L\in \N$ with $L> e$, $\epsilon \ge 0$ and $t\ge 1$
\begin{calc}
Observe that as $a_\epsilon \log L \ge 1$ we can choose $C'$ large enough to have $C  e^{C t (a_\epsilon \log L)^5} \le  e^{C' t (a_\epsilon \log L)^5}$ and also $\frac{1}{Ct}e^{-C t (a_\epsilon \log L)^5} \ge  e^{-C' t (a_\epsilon \log L)^5}$. 
\end{calc}
\begin{align}
\label{eqn:cD_L_eps_bound}
& e^{-  C t (a_\epsilon \log L)^5} \1_{[X_{[0,t]}\subset Q_L]} \le \cD_{L,\epsilon}(0 , t) \1_{[X_{[0,t]}\subset Q_L]} \le e^{C t (a_\epsilon \log L)^5}, \\
\label{eqn:Gamma_L_on_box_from_above_with_xy_dependence}
&
 e^{-C t (a_\epsilon \log L)^5} e^{ -\frac{| x-y|^2}{2 \kappa t}}
\le 
 \Gamma_t^{L,\epsilon}(x,y) 
\le 
  e^{C t (a_\epsilon \log L)^5} \frac{1}{t} e^{ -\frac{| x-y|^2}{2 ct}}
\qquad x,y \in \R^2. 
\end{align}
Before we turn to the next lemma, let us make a few observations that will be used a couple of times. 
Let $\phi \in C_c (\R^2)$ be a positive function. As $\E_{\Q_{L,\epsilon}^z}[ \phi(X_t)] = \int_{\R^2} \Gamma_t^{L,\epsilon}(z,y) \phi(y) \dd y$ we obtain the following bounds from \eqref{eqn:Gamma_L_on_box_from_above_with_xy_dependence} for all $L>e, \epsilon \ge 0, t\ge 1$ and $z \in \R^2$ 
\begin{align}
\notag 
\|\phi\|_{L^1} e^{-C  t(a_\epsilon \log L)^5} \inf_{y \in \supp \phi } e^{-\frac{|z-y|^2}{2 \kappa t}}
& \le 
\E_{\Q_{L,\epsilon}^z} [ \phi (X_t)] \\
\label{eqn:bound_exp_phi_X_t}
&\le \|\phi\|_{L^1} e^{C t(a_\epsilon \log L)^5} \sup_{y \in \supp \phi } e^{-\frac{|z-y|^2}{2 ct}}. 
\end{align}
Let $r\in (0,L]$ and $\tau_r = \inf \{ t \ge 0 : X_t \in \partial Q_r\}$. 
If $0<\delta <t$, then we have by the above  
\begin{align*}
\notag 
\E_{\Q_{L,\epsilon}^0} [ \1_{[ t- \delta < \tau_r < t]} \phi(X_t)]
& = \E_{\Q_{L,\epsilon}^0} \Big[ \1_{[t-\delta < \tau_r < t]} \E_{\Q_{L,\epsilon}^{X_{\tau_r}}} [\phi(X_{t-{\tau_r}})] \Big] \\
\cand \begin{calc} \notag 
\le \E_{\Q_{L,\epsilon}^0} \Big[ \1_{[t-\delta < \tau_r < t]} \sup_{s\in [0,t-q]} \sup_{z\in \partial Q_r} \E_{\Q_{L,\epsilon}^{z}} [\phi(X_{s})] \Big]
\end{calc}\cnewline
\cand \begin{calc}
 \le \|\phi\|_{L^1} e^{C  t( a_\epsilon \log L)^5} \sup_{z\in\partial Q_r}  \sup_{y \in \supp \phi } \sup_{s\in [0,\delta]} \frac{1}{s} e^{-\frac{|z-y|^2}{2c s}}
\end{calc} \cnewline
&  \le \|\phi\|_{L^1} e^{C t(a_\epsilon\log L)^5}  \sup_{s\in [0,\delta]} \frac{1}{s} e^{-\frac{d(\partial Q_r, \supp \phi)^2}{2c s}}, 
\end{align*}
where $d(A,B)$ is the distance between two sets $A,B \subset \R^2$. 
By computing its derivative we see that the function $f (s) = s^{- 1}  e^{- \frac{m}{s}}$ is increasing on $[0,  m ]$. 
\begin{calc}
\begin{align*}
\frac{\dd}{\dd s} s^{- 1}  e^{- \frac{m}{s}} 
= -s^{- 2}  e^{- \frac{m}{s}} + s^{- 1}  e^{- \frac{m}{s}} \frac{m}{s^2} = (\frac{m}{s} -1 ) s^{- 2}  e^{- \frac{m}{s}}. 
\end{align*}
\end{calc}  
Therefore, if $1\le  \delta \le \frac{d(\partial Q_r, \supp \phi)^2}{2c}$, then we have 
\begin{align}
\E_{\Q_{L,\epsilon}^0} [ \1_{[ t- \delta < \tau_r < t]} \phi(X_t)]
\label{eqn:bound_exp_phi_X_t_with_stopping_times}
& \le \|\phi\|_{L^1} e^{C  t(a_\epsilon \log L)^5}  e^{-\frac{d(\partial Q_r, \supp \phi)^2}{2c \delta }}. 
\end{align}
\end{obs}

\begin{lemma}
\label{lemma:estimate_fU_k_eps_phi_x}
  Let $L_t = \lfloor t (\log t)^5 \rfloor$, let $\epsilon \ge 0$, and
  let $a_{\epsilon}$ be as in Lemma~\ref{lemma:a_epsilon_bound_of_M_L_eps}. Let $x\in \R^2$ and $\phi \in C_c
  (\mathbb{R}^2)$ be positive and such that $\operatorname{supp} \phi \subset B (x, 1) =
  \{ y \in \R^2 : | x - y | \le 1 \}$. Let
  \begin{align*} 
  \fU_{k, \epsilon}^{\phi,x} (t) =\mathbb{E}_{\mathbb{Q}_{L_t^{k + 1}, \epsilon}^0} [\cD_{L_t^{k + 1}, \epsilon} (0, t) \1_{A_k} \phi (X_t)] . 
  \end{align*}
  There exist $C > 1$ and $T > 0$ such that for all $k \in \mathbb{N}$,  $\epsilon \ge 0$,  $t
  \ge T$: 
  \begin{align*} 
  \fU_{k, \epsilon}^{\phi,x} (t) \le 
  \begin{cases}
C \| \phi \|_{L^1} \exp \left( t (\log t)^5 \left[ C a_{\epsilon} ^5 (k + 1)^5 - \frac{t^{2 k - 2} (\log t)^{5 (2 k -   1)}}{C} \right] \right), & k \ge 1,\\
       C \| \phi \|_{L^1} \exp \left( (k + 1)^5 t (\log t)^5
       \left[ C  a_{\epsilon}^5 - \frac{k t}{C} \right] \right), & k
       \ge 2.
   \end{cases}
  \end{align*}
  In other words, $\mathfrak{U}^{\phi}_{k, \epsilon} (t)$ satisfies the
  same bound as $\mathfrak{U}_{k, \epsilon} (t)$, except for the factor $\|
  \phi \|_{L^1}$.
\end{lemma}

\begin{proof}
  By the proof of Lemma~\ref{lemma:bound_U_k_t} it is sufficient to prove that there exists a $C>0$ such that for all $k \in \N$, $\epsilon \ge 0$
\begin{align}
\label{eqn:fU_k_eps_phi_bound}
\fU_{k, \epsilon}^{\phi,x}(t) \le  C \|\phi \|_{L^1} e^{C a_{\epsilon}^5 t (k + 1)^5 (\log L_t)^5}  e^{- \frac{L_t^{2 k}}{C t}}. 
\end{align}
Let $T\ge 1$ be large enough such that $L_T\ge 4 c T$ and let $t\ge T$.
We prove \eqref{eqn:fU_k_eps_phi_bound} by distinguishing the two cases $| x | > \frac{L_t^k}{2}$ and $| x | \le \frac{L_t^k}{2}$. 
In the first case, $|x|>\frac{L_t^k}{2}$, observe that as $1 \le \frac{L_t^k}{4}$
\begin{align*}
\inf_{y \in B(x,1)} |y| = |x|-1 \ge \frac{L_t^k}{2} -1 \ge \frac{L_t^k}{4}.
\end{align*}
Therefore, by bounding $\1_{A_k}$ by $\1_{[X_{[0,t]} \subset Q_{L_t^{k+1}}]}$ and using \eqref{eqn:cD_L_eps_bound} and \eqref{eqn:bound_exp_phi_X_t} we estimate:   
  \begin{align*}
 \fU_{k, \epsilon}^{\phi,x}(t) 
\cand \begin{calc}
 \le e^{C a_\epsilon^5 t (k + 1)^5 (\log L_t)^5} \|\phi \|_{L^1} e^{C  a_{\epsilon}^5 t (k + 1)^5 (\log L_t)^5} e^{- \frac{(| x | - 1)^2}{2c t}}
\end{calc} \cnewline 
& \le e^{2 C a_{\epsilon}^5 t (k + 1)^5 (\log L_t)^5} \|\phi \|_{L^1} e^{- \frac{L_t^{2 k}}{32 c t}} , 
  \end{align*}
   If on the other hand $| x | \le \frac{L_t^k}{2}$, then the distance between $\partial Q_{L_t^k } $ and $B(x,1)$ is bounded as follows
\begin{align*}
d( \partial Q_{L_t^k}, B(x,1)) \ge L_t^k - |x|-1 \ge \frac{L_t^k}{4} \ge t.  
\end{align*}   
Therefore with $\tau = \inf \{ t \ge 0 : X_t \not\in Q_{L_t^k} \}$, by bounding $\1_{A_k}$ by $\1_{[\tau <t]}$ and using \eqref{eqn:cD_L_eps_bound} and \eqref{eqn:bound_exp_phi_X_t} 
  \begin{align*}
 \fU_{k, \epsilon}^{\phi,x}(t) 
\cand \begin{calc}
 \le e^{C a_{\epsilon}^5 t (k + 1)^5 (\log L_t)^5}
\mathbb{E}_{\mathbb{Q}_{L_t^{k + 1}, \epsilon}^0} \left[ \1_{[\tau < t]} \phi (X_{t}) \right] 
\end{calc} \cnewline
& \le  \| \phi \|_{L^1 (\mathbb{R})} e^{2 C a_{\epsilon}^5 t (k + 1)^5 (\log L_t)^5} e^{- \frac{L_t^{2 k}}{32 c t}}.
  \end{align*}
   
\end{proof}

\begin{corollary}
\label{cor:series-expansion-u-phi}
There exists a $T > 0$ such that for all $x\in \R^2$, for all positive $\phi \in C_c (\mathbb{R}^2)$ with $\operatorname{supp} \phi \subset B (x,1)$ and for all $t \ge T$ we have almost surely 
\begin{align*}
u^{\phi} (t, x) 
= \int_{\mathbb{R}^2} u^{\delta_x} (t, y) \phi (y) \dd y 
= \sum_{k \in \mathbb{N}_0} \fU_{k,0}^{\phi,x} (t) . 
\end{align*}
\end{corollary}

\begin{proof}
  This follows from the same arguments as Lemma~\ref{lemma:total_mass_in_sum}. 
\end{proof}

Now we apply the above corollary to the mollifier function $\psi_\epsilon$ (see~\ref{obs:setting_representation_section}). As is done before, we write $\psi_\epsilon^x$ for the shifted function $\psi_\epsilon^x(y) = \psi_\epsilon(y-x)$.

\begin{corollary}
  \label{cor:supremum-on-box}For $L_t = t (\log t)^5$ we have almost surely
  \begin{align*} \lim_{t\in \Q, t \rightarrow \infty} \left| \frac{\log (\sup_{x \in \mathbb{R}^2}  u (t, x))}{t \log t} - \frac{\log \left( \sup_{x \in Q_{L_t}} u_{L_t} (t, x) \right)}{t \log t} \right| = 0. \end{align*}
\end{corollary}

\begin{proof}
  Since $u (t, x) = u^{\delta_0} (t, x) = \lim_{\epsilon \downarrow 0}
  u^{\psi_{\epsilon}} (t, x)$ and the $L^1$-norm of
  $\psi_{\epsilon}$ is uniformly bounded in $\epsilon$, this
  follows by 
Lemma~\ref{lemma:estimate_fU_k_eps_phi_x},
  Corollary~\ref{cor:series-expansion-u-phi} together with the
  arguments from the proof of Proposition~\ref{proposition:U_t_and_U_L_t_equivalence}.
\begin{calc}
That 
\begin{align*}
\liminf_{t\in \Q, t \rightarrow \infty} 
\frac{\log \left( \sup_{x \in Q_{L_t}} u_{L_t} (t, x) \right)}{t \log t}
 > - \infty,
\end{align*}
follows as 
\begin{align*}
\sup_{x \in Q_{L_t}} u_{L_t} (t, x) 
\ge \frac{1}{L_t^2} \int_{Q_{L_t}} u_{L_t}(t,x) \dd x = \frac{U_{L_t}(t)}{L_t^2}. 
\end{align*}
\end{calc}  
  
\end{proof}

\begin{theorem}
For $a \in (0, 1)$ we have almost surely
\begin{align*} 
\lim_{t\in \Q, t \rightarrow \infty} \frac{\log U (t)}{t \log t} 
& \le \liminf_{t\in \Q, t \rightarrow \infty} \frac{\log \left( \inf_{x \in Q_{t^a}} u(t, x) \right)}{t \log t} \\
& \le \limsup_{t\in \Q, t \rightarrow \infty} \frac{\log (\sup_{x \in \mathbb{R}^2} u (t, x))}{t \log t} 
\le \lim_{t\in \Q, t \rightarrow \infty} \frac{\log U (t)}{t \log t} . 
\end{align*}
\end{theorem}

\begin{proof}
Let $C>1$ be such that~\eqref{eqn:cD_L_eps_bound}, \eqref{eqn:bound_exp_phi_X_t} and \eqref{eqn:bound_exp_phi_X_t_with_stopping_times} hold.
Let us write $M= C a_0^5$ and as before take $L_t = t (\log t)^5$. 
 To derive the upper bound we apply
  Corollary~\ref{cor:supremum-on-box} and Lemma~\ref{lemma:u_r_y_bounded_from_above_by_factor_times_total_mass} with $\delta
  = 2$ to obtain
  \begin{align*}
    \limsup_{t\in \Q, t \rightarrow \infty} \frac{\log (\sup_{x \in \mathbb{R}^2} u (t,
    x))}{t \log t} & = \limsup_{t\in \Q, t \rightarrow \infty} \frac{\log \left(
    \sup_{x \in Q_{L_t}} u_{L_t} (t, x) \right)}{t \log t}\\
    & \le \limsup_{t\in \Q, t \rightarrow \infty} \frac{\log (C) + 4 C (\log
    L_t)^5 + \log U_{L_t} (t - 2)}{t \log t}\\
    & = \lim_{t\in \Q, t \rightarrow \infty} \frac{\log U_{L_t} (t)}{t \log t} =
    \lim_{t\in \Q, t \rightarrow \infty} \frac{\log U (t)}{t \log t},
  \end{align*}
  where the last step follows from Proposition~\ref{proposition:U_t_and_U_L_t_equivalence}.
  
  The argument for the lower bound is more technical. 
Let $a\in (0,1)$ and $b\in (a,1)$.  
We estimate for $x \in Q_{t^a}$ and for $\delta_t = r_t = t^b$, $\eta \in (0,1)$:  
\begin{align}\label{eqn:infimum_bound_pr1} 
\nonumber
& \mathbb{E}_{\mathbb{Q}_{L_t}^0} \left[ \cD_{L_t} (0, t)
\1_{[X_{[0, t]} \subset Q_{L_t}]} \psi^x_{\eta} (X_t) \right]\\ 
\cand \begin{calc}
\nonumber
 \ge e^{- M \delta_t (\log L_t)^5} \mathbb{E}_{\mathbb{Q}_{L_t}^0} \left[ \cD_{L_t} (0, t - \delta_t) \1_{[X_{[0, t]} \subset Q_{L_t}]} \psi^x_{\eta} (X_t) \right]
\end{calc} \cnewline 
& 
\nonumber
 \ge e^{- M \delta_t (\log L_t)^5} \mathbb{E}_{\mathbb{Q}_{L_t}^0} \left[ \cD_{L_t} (0, t - \delta_t) \1_{[X_{[0, t - \delta_t]} \subset Q_{r_t}]} \1_{[X_{[t - \delta_t, t]}  \subset Q_{L_t}]} \psi^x_{\eta} (X_t) \right]
\\
\nonumber
& \ge e^{- M \delta_t (\log L_t)^5} \mathbb{E}_{\mathbb{Q}_{L_t}^0} \left[\cD_{L_t} (0, t - \delta_t) \1_{[X_{[0, t - \delta_t]}\subset Q_{r_t}]} \psi^x_{\eta} (X_t) \right]\\ 
& \quad - e^{- M \delta_t (\log L_t)^5} \mathbb{E}_{\mathbb{Q}_{L_t}^0}
\left[ \cD_{L_t} (0, t - \delta_t) \1_{[X_{[0, t - \delta_t]} \subset Q_{r_t}]} \1_{[X_{[t - \delta_t, t]} \not \subset Q_{L_t}]}\psi^x_{\eta} (X_t) \right].
\end{align}
For $\epsilon > 0$ we have
\begin{align*}
\mathbb{E}_{\mathbb{Q}_{L_t, \epsilon}^0} \left[ \cD_{L_t,
\epsilon} (0, t - \delta_t) \1_{[X_{[0, t - \delta_t]} \subset
Q_{r_t}]} \right] & =\mathbb{E}_0 \left[
e^{\int_0^{t - \delta_t} (\xi_{\epsilon} (X_s) -  c_\epsilon  )
\dd s} \1_{[X_{[0, t - \delta_t]} \subset Q_{r_t}]}
 \right]\\
& = U_{r_t,\epsilon}(t-\delta_t),
\end{align*}
so after passing to the limit the same is true for $\epsilon = 0$.
Therefore, for $\eta \in (0,1)$ and large enough $t$   by using \eqref{eqn:bound_exp_phi_X_t} 
\begin{align*}
\mathbb{E}_{\mathbb{Q}_{L_t}^0} & \left[ \cD_{L_t} (0, t -
\delta_t) \1_{[X_{[0, t - \delta_t]} \subset Q_{r_t}]}
\psi^x_{\eta} (X_t) \right] \\
& 
\ge 
\mathbb{E}_{\mathbb{Q}_{L_t}^0} \left[ \cD_{L_t} (0, t -
\delta_t) \1_{[X_{[0, t - \delta_t]} \subset Q_{r_t}]}
\mathbb{E}_{\mathbb{Q}_{L_t}^{X_{t-\delta_t}} }[\psi^x_{\eta} (X_{\delta_t})] \right]
\\
& 
 \ge U_{r_t} (t - \delta_t)  e^{-M \delta_t (\log L_t)^5} \inf_{z \in Q_{r_t} y \in Q_{r_t+1} } e^{-\frac{|z-y|^2}{2 \kappa t}} \\
& \ge U_{r_t} (t - \delta_t) e^{- M \delta_t (\log  L_t )^5} e^{- \frac{ 2 r_t^2}{\kappa \delta_t}}.
\end{align*}
 
This, as we will see below, is of the right order. If $t$ is large enough so that $L_t > 2 (r_t + 1)$, then
the expectation appearing in the negative term in~\eqref{eqn:infimum_bound_pr1} can be bounded from above as follows:   With $\tau = \inf \{ t> 0 : X_t \in \partial Q_{L_t} \}$, using \eqref{eqn:bound_exp_phi_X_t_with_stopping_times}  and that 
\begin{align*}
d(\partial Q_{L_t}, \supp \psi_\eta^x)\ge L_t - r_t -1 > \frac{L_t}{2}, 
\end{align*}
we have for $t$ large enough such that $\delta_t \le \frac{L_t^2}{8c}$
\begin{align*}
 \mathbb{E}_{\mathbb{Q}_{L_t}^0} & \left[ \cD_{L_t} (0, t -\delta_t) \1_{[X_{[0, t - \delta_t]} \subset Q_{r_t}]} \1_{[ X_{[t - \delta_t, t]} \not \subset Q_{L_t}]} \psi^x_{\eta} (X_t)
\right]\\
& \le e^{M (t - \delta_t) (\log L_t)^5} \mathbb{E}_{\mathbb{Q}_{L_t}^0} \left[ \1_{[ t-\delta_t < \tau < t ]} \psi^x_{\eta} (X_t) \right]\\
& \le e^{2M t  (\log L_t)^5}   e^{-\frac{L_t^2}{8 c \delta_t}} .
\end{align*}
 
Therefore we have, using that $r_t = \delta_t = t^b$:
  \begin{align*}
 \mathbb{E}_{\mathbb{Q}_{L_t}^0} & \left[ \cD_{L_t} (0, t) \1_{[X_{[0, t]} \subset Q_{L_t}]} \psi^x_{\eta} (X_t) \right]  \\
& \ge 
U_{r_t} (t - t^b) e^{- M t^b (\log L_t  )^5} e^{-\frac{2t^{b}}{\kappa}} - C e^{- 2M t (\log L_t)^5 -  \frac{L_t^2}{8c\delta_t} }. 
  \end{align*}
  and letting $\eta \downarrow 0$ we deduce that also
  \begin{align*} 
  \inf_{x \in Q_{t^a}} u (t, x) \ge 
  U_{r_t} (t - t^b) e^{- M t^b (\log L_t  )^5} e^{-\frac{2t^{b}}{\kappa}} - C e^{- 2M t (\log L_t)^5 -  \frac{L_t^2}{8c\delta_t} }
     \end{align*}
  From here we obtain by Lemma~\ref{lemma:liminf_bound} and Theorem~\ref{theorem:convergence_smooth_eigenvalues_and_asymptotics}.
  \begin{align*}
    \liminf_{t\in \Q, t \rightarrow \infty} \frac{\log \left( \inf_{x \in Q_{t^a}} u (t, x) \right)}{t \log t} & \ge \liminf_{t\in \Q, t \rightarrow \infty}
    \frac{\log U_{r_t} (t - t^b)}{t \log t} =  b \lim_{t\in \Q, t \rightarrow \infty}     \frac{\log U (t)}{t \log t}. 
  \end{align*}
Since $b \in (a, 1)$ was arbitrary, this concludes the proof.
\end{proof}

\section{Gaussian calculations: Proof of Theorem~\ref{theorem:convergence_smooth_eigenvalues_and_asymptotics}}
\label{section:gaussian_calcs}

We saved some Gaussian calculations for this last section. 
In Section~\ref{section:same_eigenvalues} we prove that our enhanced mollified noise converges up to a constant to the same limit as the differently mollified noise considered in \cite{ChvZ21}. 
In Section~\ref{section:bounds_on_noise_terms} we prove the bounds on the noise terms and its enhancement (for $\btheta = (\theta,\Theta) \in \fX_\fn^\alpha$, $\Theta$ is called the enhancement), which are used to prove Lemma~\ref{lemma:a_epsilon_bound_of_M_L_eps} to control the growth of the $L^\infty$-norms of $Z_{L,\epsilon}$ and $Y_{L,\epsilon}$ with respect to $L$.

\subsection{The limiting eigenvalues of a mollification on the full space}
\label{section:same_eigenvalues}

In this section we prove Theorem~\ref{theorem:convergence_smooth_eigenvalues_and_asymptotics}. 
We state the theorems from \cite{ChvZ21} on which it relies in Theorem~\ref{theorem:results_ChvZ21}. 
The proof then follows by the Gaussian convergence that is stated in Theorem~\ref{theorem:convergence_mollified_noise} (see also \ref{obs:proof_theorem_eigenvalues}).

Since $\xi \in \cS'$, the map $f\mapsto \langle \xi, f\rangle$ is linear and as $\| \langle \xi, f \rangle \|_{L^2(\Omega,P)} = \|f\|_{L^2(\R^d)}$, it extends to a bounded linear operator $\cW \colon L^2(\R^d) \rightarrow L^2(\Omega,P)$ such that for all $f\in L^2(\R^d)$, $\cW f$ is a centered Gaussian random variable and $E [ \cW f \cW g ] = \langle f, g \rangle_{L^2}$ for all $f, g\in L^2(\R^d)$. 
We will make abuse of notation and write ``$\langle \xi, \fn_{k,L} \rangle$'' instead of ``$\langle \cW, \fn_{k,L} \rangle$'' in the following (where $\fn_{k,L}$ is as in \ref{obs:notation_neumann}).

\begin{theorem}
\label{theorem:results_ChvZ21}
\begin{enumerate}
\item 
\label{item:fourier_multiplier_convergence}
\textnormal{\cite[Theorem 6.4]{ChvZ21}}
Let $\tau$ be in $C_{\rm c}^\infty(\R^2,[0,1])$ and equal to $1$ around $0$. 
Let $\xi_{L,\epsilon}$ be given by
\begin{align}
\label{eqn:xi_L_eps}
\xi_{L,\epsilon} : = \sum_{k\in\N_0^2} \tau(\tfrac{\epsilon}{L} k) \langle \xi, \fn_{k,L} \rangle \fn_{k,L}.
\end{align}
There exist $(\tilde C_L)_{L \in [1,\infty)}$ and $c_\tau$ in $\R$ with $\tilde C_L \xrightarrow{L\rightarrow \infty} 0$ and a  $\tilde \bxi_L= (\xi_L, \tilde \Xi_L) \in \fX_\fn^\alpha$ such that almost surely in $\fX_\fn^\alpha$ for $\alpha <-1$ 
\begin{align}
\label{eqn:tilde_bxi_L}
(\xi_{L,\epsilon}, \xi_{L,\epsilon} \reso \sigma(\rD) \xi_{L,\epsilon} - \tfrac{1}{\pi} \log \tfrac{1}{\epsilon} - \tilde C_L - c_\tau)
\rightarrow 
 \tilde \bxi_L,
\end{align}
where $\sigma(x) = (1+ \frac12 \pi^2 |x|^2)^{-1}$. 
\item \textnormal{\cite[Theorem 2.9]{ChvZ21}}
\label{item:eigenvalue_asymptotics_as_in_ChvZ21}
Let $\lambda_n(Q_L,\tilde \bxi_L)$ be the $n$-th eigenvalue of the operator $\sH_{\tilde \bxi_L}$ as in Theorem~\ref{theorem:dirichlet_summary}. 
Then for all unbounded countable sets $\I \subset (e,\infty)$,  almost surely, for all $n\in \N$, 
\begin{align*}
\lim_{L\in \I, L\rightarrow \infty}
\frac{\lambda_n(Q_L,\tilde \bxi_L)}{ \log L} = \chi. 
\end{align*}
\end{enumerate}
\end{theorem}
\begin{proof}
Actually, in \cite{ChvZ21} the setting is slightly different. 
Towards \ref{item:fourier_multiplier_convergence}: 
Let us write $\zeta_a: \R \rightarrow \R$ for the function $\zeta_a(x) = (a+ \pi^2 |x|^2)^{-1}$ with $a>0$. 
Then $\sigma = 2\zeta_2$. 
Instead of $\sigma$, in \cite{ChvZ21} $\zeta_1$ is considered in the sense that one shows the convergence of $\xi_{L,\epsilon} \reso \zeta_1(\rD) \xi_{L,\epsilon} - \frac{1}{2\pi} \log \frac{1}{\epsilon}- \tilde C_L' -c_\tau'$. 
The choice of $a=1$ in \cite{ChvZ21} is rather arbitrary, as one finds a different limit/enhancement but the operator stays the same due to cancellations in the definition of the product $\diamond$ which involves the enhancement (all the statements in \cite{ChvZ21} follow in an analogous manner). 
The difference of $\frac{1}{2\pi} \log \frac{1}{\epsilon}$ and $\frac{1}{\pi}\log \frac{1}{\epsilon}$ is then directly explained by the fact that the enhancement defined with $\sigma$ is two times the enhancement defined with $\zeta_2$, as $\sigma = 2\zeta_2$. 
Towards \ref{item:eigenvalue_asymptotics_as_in_ChvZ21}: In \cite{ChvZ21} the Hamiltonian is not defined with the factor $\frac12$. 
But it is shown that the eigenvalues of $\Delta + (2 \xi_L,4\tilde \Xi_L) \diamond$ (which is the limit of $\Delta + 2\xi_\epsilon - 4 c_\epsilon$) converge to $2\chi$, which then proves the above convergence to $\chi$ immediately. 
\end{proof}

The major difference between the setting in Theorem~\ref{theorem:results_ChvZ21} and our setting in \ref{obs:pam_smooth_to_white_noise} is that we have $\xi_{L,\epsilon}^\fn$ instead of $\xi_{L,\epsilon}$:
The $\xi_{L,\epsilon}^\fn$ is a projection of the mollified white noise $\xi_\epsilon$ onto the Neumann space on the box $Q_L$. 
Whereas, $\xi_{L,\epsilon}$ is an approximation of the white noise on the box $Q_L$ based on Fourier multipliers. 
The latter is dependent on the size of the box, whereas for $L>r$ we have that $\xi_{L,\epsilon}^\fn$ equals $\xi_{r,\epsilon}^\fn$ almost everywhere on the box $Q_r$. 
The following theorem shows that both approximations lead to the same limit, up to a constant.

\begin{theorem}
\label{theorem:convergence_mollified_noise}
Let $\xi_{L,\epsilon}$ be as in \eqref{eqn:xi_L_eps} and $\xi_{L,\epsilon}^\fn$ be as in \eqref{eqn:xi_L_eps_fn_projection}. 
There exist $(C_L)_{L\in [1,\infty)}$ and $C$ in $\R$ such that $C_L \xrightarrow{L\rightarrow \infty} 0$ and the following holds. 
In $\fX_\fn^\alpha(Q_L)$ for $\alpha <-1$ as $\epsilon \downarrow 0$, the following convergence in probability holds 
\begin{align}
\label{eqn:limit_of_mollified_noise}
(\xi_{L,\epsilon}^\fn - \xi_{L,\epsilon}, \xi_{L,\epsilon}^\fn \reso \sigma(\rD) \xi_{L,\epsilon}^\fn - \xi_{L,\epsilon} \reso \sigma(\rD) \xi_{L,\epsilon}) \xrightarrow{P} (0,  C+C_L).
\end{align}
\end{theorem}

\begin{obs}
\label{obs:proof_theorem_eigenvalues}
As a consequence of \eqref{eqn:limit_of_mollified_noise}, \eqref{eqn:limit_of_enhanced_xi_L_eps_fn} holds with $\bxi_L = (\xi_L, \tilde \Xi_L  + C_L + \tilde C_L +c_\tau )$ and $c_\epsilon = \frac{1}{2\pi} \log \frac{1}{\epsilon}+C$. 
Moreover, by the following identity 
\begin{align*}
\lambda_n(Q_L, \bxi_L) = \lambda_n(Q_L, \tilde \bxi_L) +C_L + \tilde C_L + c_\tau,
\end{align*}
we see that Theorem~\ref{theorem:convergence_mollified_noise} and Theorem \ref{theorem:results_ChvZ21} indeed imply Theorem~\ref{theorem:convergence_smooth_eigenvalues_and_asymptotics}. 
\end{obs}

We prove Theorem~\ref{theorem:convergence_mollified_noise} in Theorem~\ref{theorem:two_different_enhanced_convergences}. 
This is done using the following theorem. 

\begin{theorem}\textnormal{\cite[Theorem 11.2]{ChvZ21}}
\label{theorem:convergence_enhanced_pair_theta}
For $r\ge 1$ we let 
$X_{k,r}^\epsilon$ and $Y_{k,r}^\epsilon$ be centered Gaussian variables for $k\in \N_0^2$, $\epsilon>0$ such that every finite subset of $\{Y_{k,r}^\epsilon:  k\in \N_0^2, \epsilon>0\} \cup \{ X_{k,r}^\epsilon : k\in \N_0^2, \epsilon>0\}$ is jointly Gaussian for all $r\ge 1$. 
We write 
\begin{align}
\xi_{r,\epsilon} & = \sum_{k\in \N_0^2}  Y_{k,r}^\epsilon \fn_{k,r}, \qquad 
\theta_{r,\epsilon} = \sum_{k\in \N_0^2} X_{k,r}^\epsilon \fn_{k,r}, \\
\notag  \Theta_{r,\epsilon} & = \theta_{r,\epsilon}\reso \sigma(\rD) \theta_{r,\epsilon} - \E[ \theta_{r,\epsilon}\reso \sigma(\rD) \theta_{r,\epsilon}  ] , \\
\notag  \Xi_{r,\epsilon} & = \xi_{r,\epsilon}\reso \sigma(\rD) \xi_{r,\epsilon} - \E[ \xi_{r,\epsilon}\reso \sigma(\rD) \xi_{r,\epsilon}  ] .
\end{align}
Let $I\subset [1,\infty)$. 
We write $\fR = \{ (k,l)  \in \N_0^2 \times \N_0^2 : k_1 \ne l_1, k_2 \ne l_2 \}$. 
Let $G_{r,\epsilon}(k,l) = \E[X_{k,r}^\epsilon X_{l,r}^\epsilon - Y_{k,r}^\epsilon Y_{l,r}^\epsilon] $ and $F_{r,\epsilon} : \N_0^2\times \N_0^2 \rightarrow \R$. 
Consider the following conditions. 
\begin{align}
\label{eqn:convergence_X_k_eps_to_Z_k}
& \forall k \in \N_0^2 \ \forall r \in I; \quad \E[ |X_{k,r}^\epsilon - Y_{k,r}^\epsilon|^2 ] \xrightarrow{\epsilon\downarrow 0} 0 \\
& \forall \delta>0 \  \exists C>0\  \forall r \in I \  \forall k,l \in \N_0^2 \ \forall \epsilon >0  :
\label{eqn:bound_F}
 |F_{r,\epsilon}(k,l)| 
 \le  C \prod_{i=1}^2  (1+ |k_i - l_i | )^{\delta-1}, \\
\notag &  \forall r \in I \  \forall \delta >0 \ \exists C>0  \ \exists \epsilon_0 >0 \ \forall \epsilon \in (0,\epsilon_0) \    \forall k,l \in \N_0^2 : \ 
 \\
\label{eqn:bound_G} 
& \quad 
 |G_{r,\epsilon}(k,l)|   
\le C \begin{cases}
\prod_{i=1}^2 \left( \frac{1}{1+ |k_i - \frac{r}{\epsilon}|)^{1-\delta}} +  \frac{1}{1+ |l_i - \frac{r}{\epsilon}|)^{1-\delta}} \right)
& (k,l) \in \fR,  \\
\sum_{i=1}^2 \left( \frac{1}{1+ |k_i - \frac{r}{\epsilon}|)^{1-\delta}} +  \frac{1}{1+ |l_i - \frac{r}{\epsilon}|)^{1-\delta}} \right)
& (k,l) \in \N_0^2 \times \N_0^2 \setminus \fR.
\end{cases}
\end{align}
\begin{enumerate}
\item 
\label{item:convergence_Theta_min_Xi_eps}
Suppose that \eqref{eqn:convergence_X_k_eps_to_Z_k} holds and that \eqref{eqn:bound_F} holds for $F_{r,\epsilon}(k,l)$ being either $\E[X_{k,r}^\epsilon X_{l,r}^\epsilon]$, 
$\E[X_{k,r}^\epsilon Y_{l,r}^\epsilon]$ or $\E[Y_{k,r}^\epsilon Y_{l,r}^\epsilon]$.
 Then for $r\in I$, $\alpha<-1$, in $\fX_n^\alpha$ we have 
\begin{align*}
(
\theta_{r,\epsilon} 
 - 
 \xi_{r,\epsilon} 
,
\Theta_{r,\epsilon}
 - 
 \Xi_{r,\epsilon} 
) \xrightarrow{P} 0. 
\end{align*}
\item 
\label{item:convergence_expectation_difference}
Suppose \eqref{eqn:bound_G} holds. Then $\E[ \theta_{r,\epsilon}\reso \sigma(\rD) \theta_{r,\epsilon}  -  \xi_{r,\epsilon} \reso \sigma(\rD)  \xi_{r,\epsilon} ] \rightarrow 0$ in $\cC_\fn^{-\gamma}$ for all $\gamma>0$ and $r\in I$. 
\end{enumerate}
Consequently, if the conditions in \ref{item:convergence_Theta_min_Xi_eps} and \ref{item:convergence_expectation_difference} hold, then with $c=0$ for $r\in I$, $\alpha<-1$, in $\fX_n^\alpha$
\begin{align}
\label{eqn:difference_first_and_second_chaos}
(
\theta_{r,\epsilon} 
 - 
 \xi_{r,\epsilon} 
,
\theta_{r,\epsilon}\reso \sigma(\rD) \theta_{r,\epsilon}
 - 
 \xi_{r,\epsilon} \reso \sigma(\rD)  \xi_{r,\epsilon}
) \xrightarrow{P} (0,c). 
\end{align}
\end{theorem}

Remember that the mollifying function $\psi$ is given by $\psi(x) = \varphi(x_1) \varphi(x_2)$, for a mollifying function $\varphi$ see \ref{obs:setting_representation_section}. 

\begin{theorem}
\label{theorem:two_different_enhanced_convergences}
Let $\tau$ be in $ C_c^\infty(\R^2,[0,1])$ and equals $1$ around $0$. 
We write $\check \varphi(z) = \varphi(-z)$. 
Let $\rho : \R \rightarrow \R$ be given by 
\begin{align}
\label{eqn:our_tilde_tau}
\rho(x) = 2 \int_0^\pi [\varphi * \check \varphi]( z) \cos( x z) \dd z. 
\end{align}
Let $\tilde \tau(x) = \rho(x_1) \rho(x_2)$. 
\begin{enumerate}
\item 
\label{item:tau_and_tilde_tau}
With 
$X_{k,r}^\epsilon = \tau( \frac{\epsilon}{r} k) \langle \xi, \fn_{k,r} \rangle$,  $Y_{k,r}^\epsilon = \tilde \tau( \frac{\epsilon}{r} k) \langle \xi, \fn_{k,r} \rangle$ and $\theta_{r,\epsilon}$, $\xi_{r,\epsilon}$ as in Theorem~\ref{theorem:convergence_enhanced_pair_theta}, there exist $C>0, C_r>0$ with $C_r \xrightarrow{r\rightarrow \infty} 0$ such that \eqref{eqn:difference_first_and_second_chaos} holds with $c= C+ C_r$. 
\item 
\label{item:pairing_with_xi_eps_and_tilde_tau}

With $X_{k,r}^\epsilon = \langle \xi_\epsilon, \fn_{k,r} \rangle$ and $Y_{k,r}^\epsilon = \tilde \tau( \frac{\epsilon}{r} k) \langle \xi, \fn_{k,r} \rangle$ \eqref{eqn:difference_first_and_second_chaos} holds for $\theta_{r,\epsilon}$, $\xi_{r,\epsilon}$ as in Theorem~\ref{theorem:convergence_enhanced_pair_theta}. 
\end{enumerate}
As a consequence of \ref{item:tau_and_tilde_tau} and \ref{item:pairing_with_xi_eps_and_tilde_tau} we obtain \eqref{eqn:limit_of_mollified_noise}. 
\end{theorem}
\begin{proof}
\ref{item:tau_and_tilde_tau}
That \eqref{eqn:convergence_X_k_eps_to_Z_k} holds is immediate. 
For $F_{r,\epsilon}(k,l)$ being either $\E[X_{k,r}^\epsilon X_{l,r}^\epsilon], \E[X_{k,r}^\epsilon Y_{l,r}^\epsilon]$ or $\E[Y_{k,r}^\epsilon Y_{l,r}^\epsilon]$ we have $|F_{r,\epsilon}(k,l)| =0$ in case $k\ne l$, hence the conditions for Theorem~\ref{theorem:convergence_enhanced_pair_theta}\ref{item:convergence_Theta_min_Xi_eps} hold and so it is sufficient to prove that $\E[ \theta_{r,\epsilon}\reso \sigma(\rD) \theta_{r,\epsilon}  -  \xi_{r,\epsilon} \reso \sigma(\rD)  \xi_{r,\epsilon} ] \rightarrow (0,C+C_r)$ in $\cC_\fn^{-\gamma}$ for all $\gamma>0$. This follows from Lemma~\ref{lemma:specific_rho_satisfies_conditions} and Lemma~\ref{lemma:renormalisation_constant}. 

\ref{item:pairing_with_xi_eps_and_tilde_tau} 
For this we use Theorem~\ref{theorem:convergence_enhanced_pair_theta}. 
That \eqref{eqn:convergence_X_k_eps_to_Z_k} holds is immediate. 
The other conditions in Theorem~\ref{theorem:convergence_enhanced_pair_theta}\ref{item:convergence_Theta_min_Xi_eps} and \ref{item:convergence_expectation_difference} 
follow from (the stronger conditions in) Lemma~\ref{lemma:bounds_on_F_and_G}. 
\end{proof}

\begin{lemma}
\label{lemma:renormalisation_constant}
Let $\tilde \tau : \R^2 \rightarrow [0,1]$ be given by $\tilde \tau(x) = \rho(x_1) \rho(x_2)$, where $\rho \in C^2(\R,[0,1])$ is an even function such that $\rho'$ is bounded, $\rho(0)=1$, $\rho'(0)=0$ and $|\rho(x)| \lesssim (1+|x|)^{-\gamma}$ for some $\gamma>0$. 
Write $\tilde \xi_{L,\epsilon} = \sum_{k\in \N_0^2} \tilde \tau( \frac{\epsilon}{L}k ) \langle \xi, \fn_{k,L} \rangle \fn_{k,L}$. 
There exists a $C>0$ and for all $L \ge 1$, there exists a $C_L$ with $C_L \xrightarrow{L\rightarrow \infty} 0$ such that 
\begin{align*}
\frac14 \E[\tilde \xi_{L,\epsilon}\reso \sigma(\rD) \tilde \xi_{L,\epsilon}](0) - \tfrac{1}{\pi} \log \tfrac{1}{\epsilon}
 \xrightarrow{\epsilon \downarrow 0} C+ C_L. 
\end{align*} 
Moreover, there exists an $M>0$ such that for all $L\ge 1$ and $\epsilon \in (0,1)$ 
\begin{align*}
\left| \tfrac14 \E[\tilde \xi_{L,\epsilon}\reso \sigma(\rD) \tilde \xi_{L,\epsilon}](0) - \tfrac{1}{\pi} \log \tfrac{1}{\epsilon} \right|
\le M.
\end{align*}
\end{lemma}
\begin{proof}
The proof is similar to \cite[Lemma 6.16]{ChvZ21}.
We define $\lfloor y \rfloor = (\lfloor y_1 \rfloor , \lfloor y_2 \rfloor )$
and $h_L(y) = ( L^2 + \frac12 \pi^2 |y|^2)^{-1}$
 for $y\in \R^2$. Then (see also \cite[Section 6, equations (60), (61)]{ChvZ21})
\begin{align*}
c_{L,\epsilon} := 
 \E[\tilde \xi_{L,\epsilon}\reso \sigma(\rD) \tilde \xi_{L,\epsilon}](0)
 =   \sum_{k\in  \Z^2} \frac{\tilde \tau(\frac{\epsilon}{L} k)^2}{ L^2 + \frac12 \pi^2 |k|^2}
 =   \int_{\R^2} \tilde \tau( \tfrac{\epsilon}{L}  \lfloor  y \rfloor )^2 h_L ( \lfloor  y \rfloor) \dd y.
\end{align*}
Let 
\begin{align*}
\fD_{L,\epsilon} = c_{L,\epsilon} -   \int_{\R^2} \tilde \tau( \tfrac{\epsilon}{L}   y  )^2 h_L (   y ).
\end{align*}
We first show that there exists $a_L\in \R$ with $a_L \rightarrow 0$ such that 
\begin{align}
\label{eqn:renormalisation_constant_similar_to_integral}
\fD_{L,\epsilon} \xrightarrow{\epsilon \downarrow 0} a_L,
\end{align}
and that $\fD_{L,\epsilon}$ is uniformly bounded in $L$ and $\epsilon$. 
To shorten notation, we write $\delta = \frac{\epsilon}{L}$. 
Then $\fD_{L,\epsilon} = \fD_{L,\epsilon}^1 + \fD_{L,\epsilon}^2$, where 
\begin{align*}
\fD_{L,\epsilon}^1 =\int_{\R^2} \tilde \tau(\delta \lfloor y \rfloor)^2 (h_L(\lfloor y \rfloor) - h_L(y)) \dd y, \quad 
\fD_{L,\epsilon}^2 =\int_{\R^2} 
[\tilde \tau( \delta  \lfloor  y \rfloor )^2 
- \tilde \tau( \delta    y  )^2] h_L (   y ) \dd y .
\end{align*}
\begin{calc}
$| \lfloor y\rfloor -y| \le \sqrt{2}$. Hence $|\lfloor y \rfloor| \le \sqrt{2}+|y| \le 2 |y|$ for $|y|\ge \sqrt{2}$. On the other hand $|y| \le \sqrt{2} +|\lfloor y\rfloor|$ and so for $|y| \ge 2$, as $2- \sqrt{2} \ge \frac12$ we have $|\lfloor y \rfloor | \ge \frac12$ and so $\frac14 |y| \le  |\lfloor y\rfloor| \le 2 |y| $ for $|y|\ge 2$. 
\end{calc}
As $h_L(\lfloor y \rfloor) - h_L(y) = h_L(\lfloor y \rfloor) h_L(y) \frac12 \pi^2 ( |y|^2 - |\lfloor y \rfloor|^2)$, $h_L(\lfloor y \rfloor) \lesssim h_L(y)$ and $( |y|^2 - |\lfloor y \rfloor|^2) \lesssim 1+ |y|$, we have 
$h_L(\lfloor y \rfloor) - h_L(y) \lesssim (1+|y|) h_L(y)^2$. As the latter function is integrable over $\R^2$, it follows by Lebesgue's dominated convergence theorem that as $\delta \downarrow 0$
\begin{align*}
\fD_{L,\epsilon}^1 
 \rightarrow 
  \int_{\R^2} h_L(\lfloor y \rfloor) - h_L(y) \dd y =:  a_L, 
\end{align*}
and  $a_L \rightarrow 0$. 
As $\tilde \tau $ has values in $[0,1]$, $\fD_{L,\epsilon}^1$ is bounded by $a_L$ for all $\epsilon>0$ and therefore $\fD_{L,\epsilon}^1$ is uniformly bounded in $L$ and $\epsilon$. 
By the (multivariate) mean-value theorem, as $\nabla \tilde \tau$ is bounded, 
\begin{align*}
|\tilde \tau(\delta \lfloor y \rfloor ) - \tilde \tau( \delta y) | \lesssim \delta |\lfloor y \rfloor -y| \lesssim \delta. 
\end{align*}
As $\tilde \tau(y) \lesssim (1+|y|)^{-\gamma}$ uniformly in $y$ and $\delta (1+|\delta \lfloor y \rfloor|)^{-\gamma} \lesssim (1 + |y|)^{-\gamma}$ for $|y| \ge 2$ and $\delta \le 1$, we get (uniformly in $y$) 
\begin{align*}
|\tilde \tau( \delta  \lfloor  y \rfloor )^2 
- \tilde \tau( \delta    y  )^2| \lesssim (1+|y|)^{-\gamma }  . 
\end{align*}
Therefore $\fD_{L,\epsilon}^2$ is bounded by $\int_{\R^2} (1+|y|)^{-\gamma} (1+\pi^2|y|^2)^{-1} \dd y$ for $L\ge 1$ and $\epsilon\in (0,1)$, and
 by Lebesgue's dominated convergence theorem
\begin{align*}
\int_{\R^2} 
[\tilde \tau( \delta  \lfloor  y \rfloor )^2 
- \tilde \tau( \delta    y  )^2] h_L (   y ) \dd y \xrightarrow{\delta \downarrow 0} 0. 
\end{align*}
So we have obtained \eqref{eqn:renormalisation_constant_similar_to_integral} and therefore study the asymptotics of the following integral (which does not depend on $L$, so that the boundedness follows from the convergence)
\begin{align*}
& \int_{\R^2} \tilde \tau( \delta   y )^2 h_L (   y ) \dd y 
 \begin{calc}
= \int_{\R^2}  \frac{\tilde \tau( \frac{\epsilon}{L}   y )^2}{L^2 +\frac12 \pi^2 |y|^2} \dd y 
\end{calc}
= \int_{\R^2}  \frac{\tilde \tau(  x )^2}{ \epsilon^2 +\frac12 \pi^2  |x|^2} \dd x \\
\cand \begin{calc} = \int_{\R^2 \setminus B(0,1) }  \frac{\tilde \tau( x )^2}{ \epsilon^2 +\frac12 \pi^2  |x|^2} \dd x
+ 
\int_{B(0,1)} \frac{\tilde \tau( x )^2}{ \epsilon^2 +\frac12 \pi^2  |x|^2} \dd x 
\end{calc} \cnewline
& = \int_{\R^2 \setminus B(0,1) }  \frac{\tilde \tau( x )^2}{ \epsilon^2 +\frac12 \pi^2  |x|^2} \dd x
+ 
\int_{B(0,1)} \frac{\tilde \tau( \epsilon x )^2}{ 1 +\frac12 \pi^2  |x|^2} \dd x
+ 
\int_{A(\epsilon,1)} \frac{\tilde \tau( x )^2}{ \epsilon^2 +\frac12 \pi^2  |x|^2} \dd x,
\end{align*} 
where $A(\epsilon,1) = B(0,1) \setminus B(0,\epsilon)$. 
Because of the assumption that $\tilde \tau(y) \lesssim (1+|y|)^{-\gamma}$, the sum of the first two integrals converges in $\R$ (and the limit only depends on $\tilde \tau$): 
\begin{align*}
& \int_{\R^2 \setminus B(0,1) }  \frac{\tilde \tau( x )^2}{ \epsilon^2 +\frac12 \pi^2  |x|^2} \dd x
+ 
\int_{B(0,1)} \frac{\tilde \tau( \epsilon x )^2}{ 1 +\frac12 \pi^2  |x|^2} \dd x \\
& \xrightarrow{\epsilon \downarrow 0} \int_{\R^2 \setminus B(0,1) }  \frac{\tilde \tau( x )^2}{ \frac12 \pi^2  |x|^2} \dd x + \int_{B(0,1)} \frac{1}{ 1 +\frac12 \pi^2  |x|^2} \dd x
\end{align*}
On the other hand we have $\int_{A(\epsilon,1)} \frac{1}{\frac12 \pi^2 |x|^2} \dd x = \frac{4}{\pi} \log \frac{1}{\epsilon}$ and thus 
\begin{align*}
\int_{A(\epsilon,1)} \frac{\tilde \tau( x )^2}{ \epsilon^2 + \frac12 \pi^2  |x|^2} \dd x
- \frac{4}{\pi} \log \frac{1}{\epsilon}
= \int_{A(\epsilon,1)} \frac{(\tilde \tau( x )^2-1) \frac12 \pi^2 |x|^2 - \epsilon^2}{ (\epsilon^2 + \frac12\pi^2  |x|^2) \frac12\pi^2  |x|^2} \dd x. 
\end{align*}
Now 
\begin{align*}
\int_{A(\epsilon,1)} \frac{\epsilon^2}{ (\epsilon^2 +\frac12 \pi^2  |x|^2)\frac12 \pi^2  |x|^2} \dd x
\xrightarrow{\epsilon \downarrow 0} \int_{\R^2 \setminus B(0,1)} \frac{1}{ (1 + \frac12 \pi^2  |y|^2) \frac12 \pi^2  |y|^2} \dd y, 
\end{align*}
and 
\begin{align*}
\int_{A(\epsilon,1)} \frac{(\tilde \tau( x )^2-1) \frac12 \pi^2 |x|^2 }{ (\epsilon^2 +\frac12 \pi^2  |x|^2)\frac12 \pi^2  |x|^2} \dd x \rightarrow 
\int_{B(0,1)} \frac{\tilde \tau( x )^2-1 }{ \frac12 \pi^2  |x|^2} \dd x .
\end{align*}
Let us prove that the latter integral is finite by showing that the integrand is bounded. First observe that as $\tilde \tau(x) = \rho(x_1) \rho(x_2)$ 
\begin{calc}
and by using that $ab - 1 = (a-1)(b-1) + a- 1 + b-1$, 
\end{calc}
\begin{align*}
 \frac{\tilde \tau( x )^2-1 }{   |x|^2} 
\cand \begin{calc}
\left|\frac{(\rho(x_1)^2-1)(\rho(x_2)^2-1)+(\rho(x_1)^2-1) +(\rho(x_2)^2-1) }{x_1^2+x_2^2} \right|
\end{calc} \cnewline
& 
\le \left| \frac{\rho(x_1)^2-1}{x_1} \frac{\rho(x_2)^2-1}{ x_2 } \right|
 + \left|\frac{\rho(x_1)^2-1}{x_1^2 }\right| 
 +\left|\frac{\rho(x_2)^2-1}{x_2^2 }\right|.
\end{align*}
Now, by using that  and that $\rho(0)=1$ and $\rho'(0)=0$ and that $\rho$ is twice differentiable at $0$, the latter is bounded for $x \in B(0,1)$. 
\begin{calc}
Indeed, by L'H\^opital 
\begin{align*}
\lim_{y\rightarrow 0} \frac{\rho(y)^2-1}{y^2 }
= \lim_{y\rightarrow 0} \frac{\rho(y) \rho'(y)}{y}
=  \rho(0) \rho''(0) + \rho'(0)^2 = \rho''(0). 
\end{align*}
\end{calc}
\end{proof}

\begin{remark}
Instead of taking $\tilde \tau$ to be of the product form as in Lemma~\ref{lemma:renormalisation_constant}, one could take $\tilde \tau(x) = \rho(|x|)$ with a $\rho$ satisfying the properties and obtain the same statement. 
The idea of the proof is similar and easier as one can apply the substitution to polar coordinates. 
\end{remark}
\begin{calc}
So let us do the proof here in case $\tilde \tau (x) = \rho(|x|)$. 
By the same arguments, we obtain \eqref{eqn:renormalisation_constant_similar_to_integral}. 
The integral then equals
\begin{align*}
\int_{\R^2} \tilde \tau( \epsilon   y )^2 h_L (   y ) \dd y = 2\pi \int_0^\infty \frac{\rho(\epsilon r)^2 r}{L^2 + \pi^2 r^2} \dd r. 
\end{align*} 
As $\frac{1}{\pi^2}\log \frac{1}{\epsilon} = \int_1^{\frac{1}{\epsilon}} \frac{1}{\pi^2 r} \dd r $, it suffices to show the convergence of the following integral
\begin{align*}
\int_0^\infty \frac{\tilde \tau(\epsilon r)^2 r}{L^2 + \pi^2 r^2} - \frac{1}{\pi^2 r}\1_{[1,\frac{1}{\epsilon}]}(r) \dd r. 
\end{align*}
As $\tilde \tau$ is continuous at $0$, the integral over $[0,1]$ we have 
\begin{align*}
\int_0^1 \frac{\tilde \tau(\epsilon r)^2 r}{L^2 + \pi^2 r^2} \dd r \rightarrow \int_0^1 \frac{ r}{L^2 + \pi^2 r^2} \dd r. 
\end{align*}
Moreover, (take $s= \epsilon r$)
\begin{align*}
\int_{\frac{1}{\epsilon}}^\infty \frac{\tilde \tau(\epsilon r)^2 r}{L^2 + \pi^2 r^2}  \dd r 
= \int_{1}^\infty \frac{\tilde \tau(s)^2 s}{L^2 \epsilon^2 + \pi^2 s^2}  \dd s
 \rightarrow
 \int_{1}^\infty \frac{\tilde \tau(s)^2}{\pi^2 s}  \dd s. 
\end{align*}
Now 
\begin{align*}
\frac{\tilde \tau(\epsilon r)^2 r}{L^2 + \pi^2 r^2} - \frac{1}{\pi^2 r} =
\frac{(\tilde \tau(\epsilon r)^2 -1) r}{L^2 + \pi^2 r^2} - \frac{L^2}{\pi^2 r (L^2 + \pi^2 r^2)} , 
\end{align*}
for which the term $\frac{L^2}{\pi^2 r (L^2 + \pi^2 r^2)}$ is integrable over $[1,\infty)$. 
Then 
\begin{align*}
\int_1^{\frac{1}{\epsilon}} \frac{(\tilde \tau(\epsilon r)^2 -1) r}{L^2 + \pi^2 r^2} \dd r 
= \int_\epsilon^1 \frac{(\tilde \tau(s)^2 -1) s}{L^2 \epsilon^2 + \pi^2 s^2} \dd s
 \rightarrow \int_0^1 \frac{(\tilde \tau(s) -1)(\tilde \tau(s) +1) }{ \pi^2 s} \dd s.
\end{align*}
As $\tilde \tau'(0)=0$, the function $s\mapsto \frac{\tilde \tau(s)-1}{s}$ is bounded on $[0,1]$, so that the integral on the right-hand side is finite. 
\end{calc}

\begin{theorem}
\label{theorem:fourier_mollifiers_for_smooth_tau_converge_in_B_22_space}
Let $\tilde \tau \in C^1(\R^2,[0,1])$ satisfy $\tilde \tau(0)=1$.
Let $\gamma \in \R$. 
There exists a $C>0$ such that for all $L>0$ and $h\in H_\fn^\gamma(Q_L)$ we have $\| h - \tilde \tau(\epsilon \rD) h\|_{H_\fn^\gamma} \xrightarrow{\epsilon \downarrow 0} 0$ and for $\beta < \gamma$, $\epsilon\in (0,1)$,
\begin{align*}
\| h - \tilde \tau(\epsilon \rD) h\|_{H_\fn^\beta} \le C \epsilon^{(1-\delta) \wedge (\gamma-\beta) \delta} \|h\|_{H_\fn^\gamma}. 
\end{align*}
\end{theorem}
\begin{proof}
By \cite[Theorem 4.14]{ChvZ21} 
\begin{align*}
\| h - \tilde \tau(\epsilon \rD) h\|_{H_\fn^\beta} ^2
& \lesssim 
 \sum_{k\in \N_0^2} (1+|\tfrac{k}{L}|^2)^{\beta} (1- \tilde \tau(\tfrac{\epsilon}{L} k) )^2 \langle h, \fn_k \rangle^2 \\
&  \lesssim \Big( \sup_{k\in \N_0^2} (1+|\tfrac{k}{L}|^2)^{\beta-\gamma} (1- \tilde \tau(\tfrac{\epsilon}{L} k) )^2\Big)  \|h\|_{H_\fn^\gamma}^2.
\end{align*}
Observe that if $ |k| \ge \tfrac{La}{\epsilon^\delta}$, then $(1+|\tfrac{k}{L}|^2)^{\beta-\gamma}
\lesssim \epsilon^{2(\gamma-\beta)\delta}$. 
On the other hand, if $|k| < \frac{La}{\epsilon^\delta}$, then $|1- \tilde \tau(\tfrac{\epsilon}{L} k)| \lesssim \tfrac{\epsilon}{L} |k| \lesssim a \epsilon^{1-\delta}$ (as $\tilde \tau'$ is bounded on the ball of radius $a$). 
Therefore, as $(1- \tilde \tau(\tfrac{\epsilon}{L} k) )\lesssim 1$ and $(1+|\tfrac{k}{L}|^2)^{\beta-\gamma} \lesssim 1$
\begin{align*}
\Big( \sup_{k\in \N_0^2} (1+|\tfrac{k}{L}|^2)^{\beta-\gamma} (1- \tilde \tau(\tfrac{\epsilon}{L} k) )^2\Big) 
 \lesssim \epsilon^{2(1-\delta)\wedge 2 (\gamma-\beta) \delta} .
\end{align*}
\end{proof}

\begin{lemma}
\label{lemma:expectation_reso_xi_eps_C_0_plus_reso_constant}
Let $\tilde \tau$ and $\tilde \xi_{L,\epsilon}$ be as in Lemma \ref{lemma:renormalisation_constant}. 
Then 
$x\mapsto \E [\tilde \xi_{L,\epsilon} \reso \sigma(\rD) \tilde \xi_{L,\epsilon} (x) ] - \frac14 \E [\tilde \xi_{L,\epsilon} \reso \sigma(\rD) \tilde \xi_{L,\epsilon} (0)] $ converges in $\cC_\fn^{-\gamma}$  to a limit that is independent of $\tilde \tau$ as $\epsilon \downarrow 0$ for all $\gamma>0$. Moreover, there exists a $M>0$ such that for all $L\ge 1$ and $\epsilon >0$
\begin{align}
\label{eqn:bound_difference_expectations}
\|\E [\tilde \xi_{L,\epsilon} \reso \sigma(\rD) \tilde \xi_{L,\epsilon} (\cdot) ] -  \tfrac14 \E [\tilde \xi_{L,\epsilon} \reso \sigma(\rD) \tilde \xi_{L,\epsilon} (0)]\|_{\cC_\fn^{-\gamma}} \le  M. 
\end{align}
\end{lemma}
\begin{proof}
The proof of the convergence follows along the same lines as \cite[Theorem 6.15]{ChvZ21} by using Theorem~\ref{theorem:fourier_mollifiers_for_smooth_tau_converge_in_B_22_space} instead of \cite[Theorem 6.14]{ChvZ21}. 
The bound \eqref{eqn:bound_difference_expectations} is not proven in \cite[Theorem 6.15]{ChvZ21}, but can be derived from the decomposition done in the proof and by the bound 
\begin{align*}
\frac{1}{L} \sum_{m \in \frac{1}{L} \N_0} \frac{1}{(1+m)^{1+\gamma}} \lesssim 1,
\end{align*}
which for example can be concluded by \cite[Lemma 11.7]{ChvZ21}. 
\end{proof}

\begin{lemma}
\label{lemma:specific_rho_satisfies_conditions}
Let $\rho$ be as in \eqref{eqn:our_tilde_tau}. 
Then $\rho \in C^\infty(\R,[0,1])$ and $\rho$ is an even function such that $\rho'$ is bounded, $\rho(0)=1$, $\rho'(0)=0$ and $|\rho(x)| \lesssim (1+|x|)^{-1}$. 
\end{lemma}
\begin{proof}
By Leibniz integral rule $\rho$ is $C^\infty$. 
That  $\rho'$ is bounded follows by the identity  $\rho'(x) = - 2 \int_0^\pi z [\varphi * \check \varphi](z) \sin(xz) \dd z$ and that $z\mapsto z [\varphi * \check \varphi](z)$ is integrable. 
The bound $\rho(x) \lesssim (1+|x|)^{-1}$ follows by the following identity that holds for $x\ne 0$
\begin{align*}
\rho(x) 
&= \frac1x \int_0^\pi [\varphi * \check \varphi](z) \frac{\dd}{\dd z} \sin(xz) \dd z 
= - \frac{1}{x} \int_0^\pi [\varphi * \check \varphi]'(z) \sin(xz) \dd z. 
\end{align*}
\end{proof}

\begin{lemma}
\label{lemma:bounds_on_F_and_G}
Let $\tilde \tau$ be as in Theorem~\ref{theorem:two_different_enhanced_convergences}. 
For $X_{k,L}^\epsilon = \langle \xi_\epsilon , \fn_{k,L} \rangle$ and $F_{L,\epsilon}(k,l) = \E[  X_{k,L}^\epsilon X_{l,L}^\epsilon ]$ the following holds 
\begin{align}
\label{eqn:bound_F_better}
& \exists C>0\ \forall L>0 \ \forall k,l \in \N_0^2, k\ne l \ \forall \epsilon >0 ; \quad
 |F_{L,\epsilon}(k,l)|  \le C \prod_{i=1}^2 \Big(\epsilon \wedge \frac{1}{(k_i+l_i)\vee 1}\Big) ,  \\
\label{eqn:bound_G_better}
& \exists C>0 \ \forall L>0 \ \forall k \in \N_0^2  \ \forall \epsilon >0 ; \quad  
|F_{L,\epsilon}(k,k)- \tilde \tau(\tfrac{\epsilon}{L}k)^2 |  \le C (\epsilon \wedge \tfrac{1}{k_1} + \epsilon \wedge \tfrac{1}{k_2}). 
\end{align}
\eqref{eqn:bound_F_better} also holds for $F_{L,\epsilon}(k,l) = \E[ \langle \xi_\epsilon , \fn_{k,L} \rangle \langle \xi, \fn_{l,L} \rangle ]$. 
\end{lemma}
\begin{proof}
Let us first rewrite $F_{L,\epsilon}(k,l)$ by using that 
$\langle \xi_\epsilon , \fn_{k,L} \rangle = 
\langle \psi_\epsilon * \xi , \fn_{k,L} \rangle
= \langle \xi , \check \psi_\epsilon * \fn_{k,L} \rangle$ where $\check \psi_\epsilon(x) = \psi_\epsilon(-x)$ (see also \cite[Theorem 11.5]{DuKo10}), so that 
\begin{align*}
F_{L,\epsilon}(k,l) 
& =  \E[  \langle \xi_\epsilon , \fn_{k,L} \rangle  \langle \xi_\epsilon , \fn_{l,L} \rangle ]
= \langle \check \psi_\epsilon * \fn_{k,L}, \check \psi_\epsilon * \fn_{l,L} \rangle_{L^2(\R^2)} \\
& = \langle \fn_{k,L},  \psi_\epsilon * \check \psi_\epsilon * \fn_{l,L} \rangle_{L^2(\R^2)}
=  \prod_{i=1}^2 \langle \fn_{k_i,\pi},  K_\epsilon * \fn_{l_i,\pi} \rangle_{L^2(\R)},
\end{align*}
where $K_\epsilon(x) = \frac{L}{\pi} \varphi_\epsilon * \check \varphi_\epsilon(\tfrac{L}{\pi}x )$ and $\fn_{m,\pi}$ for $m\in \N_0$ is the element in the Neumann basis for $L^2([-\frac{\pi}{2},\frac{\pi}{2}])$. 
Observe that $K:=K_1$ is a smooth even  function with $\supp K \subset [-\pi,\pi]$ (as $L\ge 1$) that integrates to one. 

Because $\langle \fn_{k_i,\pi},  K_\epsilon * \fn_{l_i,\pi} \rangle_{L^2(\R)}= 
\langle \cT_y \fn_{k_i,\pi},  K_\epsilon * \cT_y \fn_{l_i,\pi} \rangle_{L^2(\R)}$
for all $y\in \R$, and for $y= - (\frac{\pi}{2}, \frac{\pi}{2})$ and all $k\in \N_0$ and $x\in [0,\pi]$, 
\begin{align*}
\tfrac{1}{\nu_k} \cT_{y}  \fn_{k,\pi}(x) = \sqrt{\tfrac{2}{\pi}} \cos( k x), 
\end{align*}
we have for $k,l \in \N_0$ we have by a substitution and Fubini's theorem, 
\begin{calc}
with 
\begin{align*}
A & := \{ (y,z) \in \R^2 : 0 \le y \le \pi, -y \le z  \le  \pi -y \} \\
&  = \{ (y,z) \in \R^2 : -\pi \le  -y \le z \le \pi -y \le \pi \} \\
& = \{ (y,z) \in \R^2 : -\pi \le z \le \pi,   (-z)\vee 0  \le y \le (\pi -z)\wedge \pi  \},
\end{align*}
that
\end{calc}
\begin{align*}
A_{k,l}^\epsilon 
:= \frac{\pi}{2  \nu_k \nu_l }
 \langle \fn_{k,\pi},  K_\epsilon * \fn_{l,\pi} 
  \rangle_{L^2(\R)}
& = \int_0^{\pi}  \int_0^{\pi}  K_\epsilon (x - y) \cos (k x) \cos
  (l y) \dd x \dd y \\
  & = 
\int_0^{\pi}  \int_{-y}^{\pi-y}  K_\epsilon (z) \cos (k (z+y)) \cos
  (l y) \dd z \dd y \\
\cand \begin{calc}
    = 
\int_\R \int_\R \1_A(y,z) K_\epsilon (z) \cos (k (z+y)) \cos
  (l y) \dd z \dd y 
\end{calc} \cnewline 
& = 
\int_{-\pi}^{\pi}    K_\epsilon (z) 
\int_{(-z)\vee 0}^{(\pi-z)\wedge \pi} \cos (k (z+y)) \cos   (l y) \dd y \dd z. 
\end{align*}
We write 
\begin{align*}
f_{k,l}(z) = \int_{(-z)\vee 0}^{(\pi-z)\wedge \pi} \cos (k (z+y)) \cos   (l y) \dd y.
\end{align*}
Note that $( \pi - z) \wedge \pi = \pi - (z \vee 0)= \pi - z^+$ and $(-z) \vee 0 = z^-$. 
Using that the product of cosine functions can be written as the sum of two cosines and that $\sin( \pi m + x) = (-1)^m \sin(x)$ for $m\in \Z$, we obtain for $k \ne l $
\begin{align*}
 f_{k,l}(z)  
& = \tfrac12 \int_{z^-}^{\pi - z^+} \cos (k z+ (k+l)y)) + \cos (k z+ (k-l)y))  \dd y \\
\cand \begin{calc}
 = \tfrac12 \frac{(-1)^{k+l} \sin (k z - (k+l)z^+ ) - \sin (k z+ (k+l)z^-)}{k+l} 
 \end{calc} \cnewline
\cand \begin{calc} \qquad +  \tfrac12\frac{(-1)^{k-l} \sin (k z - (k-l)z^+ ) - \sin (k z+ (k-l)z^-)}{k-l} 
\end{calc} \cnewline
\cand \begin{calc} = 
\begin{cases}
\tfrac12\frac{(-1)^{k+l} \sin (lz) - \sin (k z)}{k+l} 
 + \tfrac12 \frac{(-1)^{k-l} \sin (l z ) - \sin (k z)}{k-l} \cand z \ge 0, \cnewline
\tfrac12 \frac{(-1)^{k+l} \sin (k z  ) - \sin (l z)}{k+l} 
 + \tfrac12 \frac{(-1)^{k-l} \sin (k z ) - \sin (l z)}{k-l} \cand z \le 0 ,
\end{cases} 
\end{calc} \cnewline
& = 
\begin{cases}
\left[ (-1)^{k+l} \sin (lz) - \sin (k z) \right] \frac{k}{k^2 -l^2}  & z \ge 0 , \\
\left[ (-1)^{k+l} \sin (kz) - \sin(l z) \right] \frac{k}{k^2 -l^2}  & z \le 0 .
\end{cases}
\end{align*}
In case $k+l$ is odd, $f_{k,l}$ is an odd function so $\int_{-\pi}^{\pi}  K_\epsilon (z) f_{k,l}(z) \dd z =0$.
If $k+l$ is even, $f_{k,l}$ is an even function. 
As $K$ is even one has $A_{k,l}^\epsilon = A_{l,k}^\epsilon$ and so for $k, l$ such that  $k+l$ is even
\begin{align*}
A_{k,l}^\epsilon 
&= \int_{-\pi}^\pi K_\epsilon(z) \tfrac12[ f_{k,l}(z) +f_{l,z}(z)] \dd z \\
\cand \begin{calc}
=  \int_{0}^\pi K_\epsilon(z) [ 
 \sin (lz) - \sin (k z) ] (\frac{k}{k^2 -l^2} + \frac{l}{l^2 -k^2}) ] \dd z
 \end{calc} \cnewline
& = \frac{1}{k+l} \int_{0}^\pi K_\epsilon(z) [ 
 \sin (lz) - \sin (k z) ]  \dd z.
\end{align*}
As $K$ is a positive function with support in $[-\pi,\pi]$ that integrates to $1$, we obtain $|A_{k,l}^\epsilon| \le  \frac{1}{k+l}$. On the other hand 
\begin{align*}
|A_{k,l}^\epsilon|  
= \Big| \frac{1}{k+l} \int_{0}^\pi K(z) [ 
 \sin (\epsilon l z) - \sin (\epsilon k z) ]  \dd z \Big|
\begin{calc} \lesssim \epsilon \frac{|k-l|}{k+l} \end{calc}
\lesssim \epsilon.
\end{align*}
This proves \eqref{eqn:bound_F_better}. 
As also $ \E[ \langle \xi_\epsilon , \fn_{k,L} \rangle \langle \xi, \fn_{l,L} \rangle ]$ is of the form $\prod_{i=1}^2 \langle \fn_{k_i,\pi},  K_\epsilon * \fn_{l_i,\pi} \rangle_{L^2(\R)}$ but with  $K_\epsilon(x) = \frac{L}{\pi} \psi_\epsilon(\frac{L}{\pi}x)$ for which still $K_1$ is a smooth even  function with $\supp K_1 \subset [-\pi,\pi]$ (as $L\ge 1$) that integrates to one, the above also proves \eqref{eqn:bound_F_better} for $F_{L,\epsilon}(k,l) = \E[ \langle \xi_\epsilon , \fn_{k,L} \rangle \langle \xi, \fn_{l,L} \rangle ]$. 

Now for $k=l \ne 0$, we have 
\begin{align*}
 f_{k,k}(z) 
\cand \begin{calc}
 =   \int_{z^-}^{\pi-z^+} \cos (k (z+y)) \cos (k y) \dd y 
\end{calc} \cnewline
& =  \frac12 \int_{z^-}^{\pi-z^+}  \cos (k z+ 2k y)) + \cos (k z)  \dd y \\
& =  \frac12 \frac{\sin (k z - 2k -z^+) - \sin (k z + 2kz^-)}{2k}  + \frac12 [\pi -|z| ] \cos(kz) \\
\cand\begin{calc}
= \frac12 [\pi -|z| ] \cos(kz) + 
\begin{cases}
 \frac12 \frac{\sin (k z - 2k z) - \sin (k z )}{2k} \cand z \ge 0, \cnewline
 \frac12 \frac{\sin (k z) - \sin (- kz )}{2k} \cand z \le 0 , 
\end{cases}
\end{calc}\cnewline 
& =  \frac{- \sin (k |z| )}{2k}  + \frac12 [\pi -|z| ] \cos(kz),
\end{align*}
and for $k=l=0$, $f_{k,l}(z) = \pi - |z| =  [\pi - |z|]$. 
Observing that $\nu_0^2 = 2$ (for the one dimensional $0$). 
Hence for $k\in \N_0^2$ 
\begin{align*}
F_{\epsilon}(k,k) 
= \prod_{i=1}^2 \frac{2\nu_k^2}{\pi} A_{k_i,k_i}^\epsilon
= \prod_{i=1}^2 2\int_0^\pi K_\epsilon(z) \Big( \frac{- \sin (k_i z )}{\pi k_i} \1_\N(k_i) + [1 -\tfrac{1}{\pi} |z| ]\cos(k_i z) \Big) \dd z.
\end{align*}
Observe that as $K_\epsilon(z) = \frac{L}{\epsilon \pi} [\varphi * \check \varphi](\frac{L}{\epsilon \pi } z)$,
\begin{align*}
\prod_{i=1}^2 2 \int_0^\pi K_\epsilon(z) \cos(k_i z)  \dd z 
= \prod_{i=1}^2  2 \int_0^1 [\varphi * \check \varphi]( z) \cos(\tfrac{\epsilon \pi}{L} k_i z)  \dd z 
= \tilde \tau( \tfrac{\epsilon}{L} k)^2.
\end{align*}
As $F_{L,\epsilon}(k,k) - \tilde \tau(\frac{\epsilon}{L}k)$ is of the form $a_{k_1} a_{k_2} - b_{k_1} b_{k_2}$ and $|a_l|$ and $|b_l|$ are bounded by one, we have $| a_{k_1} a_{k_2} - b_{k_1} b_{k_2}| \le |a_{k_1} - b_{k_1} | + | a_{k_2} -  b_{k_2}|$. Therefore it is sufficient to show the following for all $k\in \N$
\begin{align}
\label{eqn:one_dimensional_difference_bound}
\left| -\frac{1}{k} \int_{0}^{\pi} K (z)  \sin(\epsilon kz)  \dd z
    + \epsilon \int_{0}^{\pi}    K (z)   z \cos(\epsilon kz)  \dd z \right|
    \lesssim \epsilon \wedge \frac{1}{k}.
\end{align}
By partial integration we have 
\begin{align*}
\frac{1}{k} \int_{0}^{\pi} K (z)  \sin(\epsilon kz)  \dd z
& = - \epsilon  \int_0^\pi \int_{-\infty}^z K(x) \dd x \cos(\epsilon k z) \dd z
+ \frac{1}{k} \int_{-\infty}^\pi K(z) \dd z \sin(\epsilon k \pi), \\
\epsilon \int_{0}^{\pi}    K (z)   z \cos(\epsilon kz)  \dd z
& = \frac{1}{ k } \int_{0}^{\pi}    K (z)   z \frac{\dd}{\dd z} \sin(\epsilon kz)  \dd z\\                                                                                                                                                                                                                                                                                                                                                                                                                                                                                                                                                                                                                                                                                                                                                                                                                                                                                                                                                                                                                                                                                                                                                                                                                                                                                                                                                                                                                                                                                                                                                                                                                                                                                                                                                                                                                                                                                                                                                                                                                                                                                                                                                                                                                                                                                                                                                                                                                                                                                                                                                                                                                                                                                                                                                                                                                                                                                                                                                            
& = - \frac{1}{ k } \int_{0}^{\pi}   (z K' (z) +  K(z))  \sin(\epsilon kz)  \dd z.
\end{align*}
By using that $\sin(\epsilon k \pi) \le \epsilon k \pi$ and that 
$\int_{0}^{\pi}   |z K' (z)| + | K(z)|   \dd z $ is bounded by a constant that does not depend on $L$, we obtain \eqref{eqn:one_dimensional_difference_bound}.
\begin{calc}
Indeed, when writing $f= \varphi * \check \varphi$ we have $K_\epsilon(z) = \frac{L}{\epsilon \pi} f( \frac{L}{\epsilon \pi} z)$ and so 
\begin{align*}
\int_0^\pi |K_\epsilon (z)| \dd z 
& = \int_0^{\frac{L}{\epsilon}}  f(x) \dd x 
= \int_0^{1}  f(x) \dd x  = 1, \\
\int_0^\pi |z K_\epsilon' (z)| \dd z 
& = \int_0^\pi \tfrac{L}{\epsilon \pi} | \tfrac{L}{\epsilon \pi} f'  (\tfrac{L}{\epsilon \pi} z)| \dd z 
 = \int_0^{\frac{L}{\epsilon}} |x f'(x)| \dd x 
  = \int_0^{1} |x f'(x)| \dd x <\infty. 
\end{align*}
\end{calc}
\end{proof}

\subsection{Bounds on the noise terms}
\label{section:bounds_on_noise_terms}

In this section we prove bounds on the size of the norms of the noise term $\xi_{L,\epsilon}^\fn$ and its enhancement with respect to the size of the box, $L$, in Lemma~\ref{lemma:growth_of_noise_in_terms_of_box_size}. 
We derive these bounds from $L^p$ bounds, provided in Lemma~\ref{lemma:L_p_bounds_on_xi_L_eps_and_enhancement}. 
First some auxiliary lemmas. 

\begin{lemma}
\label{lemma:delta_i_bounds}
For all $\gamma \in (0,1)$ there exists a $\fC>1$ such that for all $L \ge 1$, 
$i \in \N_{-1}$, $\epsilon >0$, $x\in Q_r$
\begin{align}
\label{eqn:delta_i_bounds}
\E[ | \Delta_i \xi_{L,\epsilon}^\fn |(x)^2]
 \le \fC 
2^{(2+\gamma )i}, 
\qquad 
\E [ |\Delta_i  \Xi_{L,\epsilon}^\fn |(x)^2] 
\le \fC  2^{\gamma i }.
\end{align}
\end{lemma}
\begin{proof}
The proof is similar to the proof of \cite[Theorem 11.1]{ChvZ21}, based on \cite[Theorem 11.5 and 11.12]{ChvZ21}, in which on the right-hand side the factor $L^{2\gamma}$ appears. 
As mentioned before, the bound \eqref{eqn:bound_F_better} we have obtained in Lemma \ref{lemma:bounds_on_F_and_G} is stronger than \eqref{eqn:bound_F}. It also assures that we do not obtain anything that depends on $L$ on the right-hand side of \eqref{eqn:delta_i_bounds}. Indeed, similar to the proof of \cite[Theorem 11.5]{ChvZ21}
\begin{align*}
2^{-(2+2\gamma)i} \E[ | \Delta_i \xi_{L,\epsilon}^\fn (x) |_{L^\infty}^2]
\cand \begin{calc}
\notag
\lesssim L^{-2} \sum_{  k,l \in  \N_0^2} 
\frac{1}{(1+| \frac{k}{L} |)^{\frac{2+2\gamma}{2}}}
\frac{1}{(1+| \frac{l}{L} |)^{\frac{2+2\gamma}{2}}}  
|\E[ X_{k,r}^\epsilon X_{l,r}^\epsilon]| 
\end{calc}\cnewline
& \lesssim 
\Big(
 \sum_{  k,l \in \frac{1}{L} \N_0} 
 L^{-2}
\frac{1}{(1+ k )^{\frac{1+\gamma}{2}}}
\frac{1}{(1+ l )^{\frac{1+\gamma}{2}}}
 \frac{1}{ k+l } \Big)^2 \\
& \lesssim 
\Big(
 \sum_{  k \in \frac{1}{L} \N_0} 
 L^{-1}
\frac{1}{(1+ k )^{1+ \frac{\gamma}{2}}} \Big)^4 \le \fC, 
\end{align*}
for some $\fC$ that does not depend on $L$, by \cite[Lemma 11.7]{ChvZ21}. 
Similarly, one can prove the bound on $\Xi_{L,\epsilon}^\fn$ by following the proof of \cite[Theorem 11.12]{ChvZ21}. 
\end{proof}

We present the following consequence of the Gaussian hypercontractivity, which basically generalises the fact that the $p$-th moment of a Gaussian can be bounded by a multiple of its second moment to the power $\frac{p}{2}$. 

\begin{lemma}\textnormal{\cite[Theorem 1.4.1, see also the end of p.62]{Nu09}}
\label{lemma:hyper_contractivity}
Let $V$ be in the $k$-th Wiener chaos for some $k\in \N$. Then $E[V^p] \le (p-1)^{\frac{k p}{2}} E[V^2]^{\frac{p}{2}}$. For $p>1$
\end{lemma}
\begin{calc}
Indeed, by \cite{Nu09} at the end of page 62: 
For $t>0$ such that $p= 1+ e^{2t} $ it holds that $E[V^p]^\frac{1}{p} \le e^{nt} E[V^2]^{\frac12}$ and $e^{nt} = (e^{2t})^{\frac{n}{2}} = (p-1)^{\frac{n}{2}}$. 
\end{calc}

\begin{lemma}
\label{lemma:L_p_bounds_on_xi_L_eps_and_enhancement}
For all $\gamma\in (0,1)$ there exist $C,M>0$ such that for all $L \ge 1$, 
 $\epsilon \in (0,1)$, $p>1$
\begin{align}
\label{eqn:bound_cC_norm_xi_eps}
E[ \| \xi_{L,\epsilon}^\fn  \|_{\cC_\fn^{-1-\gamma -\frac{2}{p}}}^p ] 
& \le C L^{2}  p^{\frac{p}{2}}
 M^\frac{p}{2}, \\
 \label{eqn:bound_cC_norm_xi_eps_reso_etc}
E[ \| \xi_{L,\epsilon}^\fn \reso \sigma(\rD) \xi_{L,\epsilon}^\fn - c_\epsilon  \|_{\cC_\fn^{-\gamma -\frac{2}{p}}}^p ] 
& \le C L^{2}  p^{p}
 M^\frac{p}{2}.
\end{align}

\end{lemma}
\begin{proof}
Using Lemma~\ref{lemma:hyper_contractivity} along the same lines as \cite[Lemma 6.10]{ChvZ21}, there exists a $C>0$ independent of $\xi_{L,\epsilon}^\fn$ and $\xi_{L,\epsilon}^\fn \reso \sigma(\rD) \xi_{L,\epsilon}^\fn - c_\epsilon$ such that with $\fC$ as in Lemma \ref{lemma:delta_i_bounds} for all $L\ge 1$, $\epsilon>0$ and $p>1$ we have \eqref{eqn:bound_cC_norm_xi_eps} and 
 \begin{align*}
 E[ \| \Xi_{L,\epsilon}^\fn  \|_{\cC_\fn^{-\gamma -\frac{2}{p}}}^p ] 
& \le C L^2  p^{p}
 \fC^\frac{p}{2}.
 \end{align*}
Therefore it is sufficient to show that there exists an $M>0$ such that 
 \begin{align*}
 \| E[\xi_{L,\epsilon}^\fn \reso \sigma(\rD) \xi_{L,\epsilon}^\fn] - c_\epsilon \|_{\cC_\fn^{-\gamma -\frac{2}{p}}}^p 
& \le C L^2  p^{p}
 M^\frac{p}{2}.
 \end{align*}
This follows from Lemma~\ref{lemma:expectation_reso_xi_eps_C_0_plus_reso_constant}
and Lemma~\ref{lemma:renormalisation_constant} 
as 
$E[\xi_{L,\epsilon}^\fn \reso \sigma(\rD) \xi_{L,\epsilon}^\fn] - c_\epsilon$ 
is the sum of $ E[\xi_{L,\epsilon}^\fn \reso \sigma(\rD) \xi_{L,\epsilon}^\fn] 
-\frac14 E[\xi_{L,\epsilon}^\fn \reso \sigma(\rD) \xi_{L,\epsilon}^\fn(0)]$ and $\frac14 E[\xi_{L,\epsilon}^\fn \reso \sigma(\rD) \xi_{L,\epsilon}^\fn(0)] 
- c_\epsilon$. 
\end{proof}

\begin{lemma}
\label{lemma:growth_of_noise_in_terms_of_box_size}
For $L>1$ let $\bxi_L=(\xi_L,\Xi_L) $ be as in  \eqref{eqn:limit_of_enhanced_xi_L_eps_fn}.
For all $\gamma>0$ there exist a $h_0>0$ such that for all $h \in [0,h_0]$ 
\begin{align}
\label{eqn:norm_noise_in_exponent}
& \sup_{\epsilon>0, L > 1} L^{-4} E \Big[ e^{h \|\xi_{L,\epsilon}^\fn \|^2_{\cC_\fn^{-1-\gamma}} } \Big]
+ 
L^{-2} E \Big[ e^{h \|\xi_{L,\epsilon}^\fn \reso \sigma(\rD) \xi_{L,\epsilon}^\fn - c_\epsilon \|_{\cC_\fn^{-\gamma}} } \Big] <\infty, \\
\label{eqn:norm_noise_in_exponent_limit_eps}
& \sup_{ L > 1} L^{-4} E \Big[ e^{h \|\xi_{L} \|^2_{\cC_\fn^{-1-\gamma}} } \Big]
+ 
L^{-2} E \Big[ e^{h \| \Xi_L \|_{\cC_\fn^{-\gamma}} } \Big] <\infty. 
\end{align}
As a consequence, 
for all $\epsilon>0$ almost surely 
\begin{align}
\label{eqn:sup_noise_terms_divided_by_log_L_epsilon}
A_{\gamma,\epsilon}
& := \sup_{L\in \N, L > e}  \frac{\|\xi_{L,\epsilon}^\fn\|^2_{\cC_\fn^{-1-\gamma}}  
+ \|\xi_{L,\epsilon}^\fn \reso \sigma(\rD) \xi_{L,\epsilon}^\fn - c_\epsilon \|_{\cC_\fn^{-\gamma}} }{\log L} <\infty, \\
\label{eqn:sup_noise_terms_divided_by_log_L}
A_{\gamma} 
& := \sup_{L \in \N, L> e}  \frac{\|\xi_{L}\|^2_{\cC_\fn^{-1-\gamma}}  
+ \|\Xi_L  \|_{\cC_\fn^{-\gamma}} }{\log L} <\infty,
\end{align}
and moreover $E[e^{h A_{\gamma}}]<\infty$ and $\sup_{\epsilon>0} E[ e^{h A_{\gamma,\epsilon}}]<\infty$ for all $h \in [0,h_0]$.
\end{lemma}
\begin{proof}
Let $\gamma>0$ and $m\in \N$ be such that $\frac{2}{m} < \gamma$. 
Let us write $S_{L,\epsilon,1} = \| \xi_{L,\epsilon}^\fn\|_{\cC_\fn^{-1-2\gamma}}$ and $S_{L,\epsilon,2} = \| \xi_{L,\epsilon}^\fn \reso \sigma(\rD) \xi_{L,\epsilon}^\fn - c_\epsilon  \|_{\cC_\fn^{-2\gamma}}$. 
For $k\in \{1,2\}$, by using subsequently Jensen's inequality, \eqref{eqn:bound_cC_norm_xi_eps} and \eqref{eqn:bound_cC_norm_xi_eps_reso_etc}, i.e., $E[S_{L,\epsilon,k}^{\frac{2p}{k}}]\le C L^2 p^p \fC^p $ (we may and do assume $\fC > 1$)
 and $\frac{1}{n!} \le (\frac{e}{n})^{n}$, we obtain for $h \ge 0$ and $\epsilon>0$
\begin{align*}
E[ \exp ( h S_{L,\epsilon,k}^{\frac{2}{k}})] 
& = \sum_{n=0}^\infty \frac{h^n}{n!} E[ S_{L,\epsilon,k}^{\frac{2n}{k}}]
\le \sum_{n=0}^{m-1} \frac{h^n}{n!} E[ S_{L,\epsilon,k}^{m}]^{\frac{2n}{km}} 
+ \sum_{n=m}^\infty \frac{h^n}{n!} E[ S_{L,\epsilon,k}^{\frac{2n}{k}}] \\
& \le  \sum_{n=0}^{m-1} \frac{h^n}{n!} ( C L^2 m^{\frac{km}{2}} \fC^{\frac{m}{2}})^{\frac{2n}{km}} 
+ \sum_{n=m}^\infty \frac{h^n}{n!} C L^2 n^n \fC^{\frac{n}{2}} \\
& \le L^{\frac{4}{k}} \exp( h ( C m^{\frac{2m}{k}} \fC^{\frac{m}{2}})^{\frac{2}{km}})) 
+ C L^2 \sum_{n=m}^\infty (h  e \fC^{\frac{1}{2}})^{n}.
\end{align*}
For $h_0 = (2e \fC^{\frac{1}{2}})^{-1}$ the latter series is finite, i.e., there exists a $M>0$ such that for all $h \in [0,h_0]$, $L>1$, $\epsilon>0$ and $k\in \{1,2\}$ we have $E[ \exp ( h S_{L,\epsilon,k}^{\frac{2}{k}})]  \le M L^{\frac{4}{k}} $.  This implies 
\eqref{eqn:norm_noise_in_exponent}.
Let us write $S_{L,0,k} = \lim_{\epsilon\downarrow 0} S_{L,\epsilon,k}$ for $k\in \{1,2\}$. By Fatou's lemma we obtain $E[ \exp ( h S_{L,0,k}^{\frac{2}{k}})]  \le M L^{\frac{4}{k}}$ for $k\in \{1,2\}$ and thus also obtain \eqref{eqn:norm_noise_in_exponent_limit_eps}. 
Moreover, for  $a> \frac{4}{k}+1$ 
\begin{align}
\label{eqn:sum_L_of_L_power_expectation_moment_gen_fu}
 \sum_{L\in \N} L^{-a} \sup_{\epsilon\in [0,\infty)} E[ \exp ( h_0 S_{L,\epsilon,k}^{\frac{2}{k}})] <\infty. 
\end{align}
Therefore $S_\epsilon := \sum_{L\in \N} L^{-a} \exp ( h_0 S_{L,\epsilon,k}^{\frac{2}{k}})$ is almost surely finite for all $\epsilon \in [0,\infty)$. 
\begin{calc}
This implies $L^{-a} \exp ( h_0 S_{L,\epsilon,k}^{\frac{2}{k}}) \le S_\epsilon$ and thus 
$  h_0 S_{L,\epsilon,k}^{\frac{2}{k}} \le a \log L + \log  S_\epsilon$.
\end{calc}
Consequently, $S_{L,\epsilon,k}^{\frac{2}{k}} \le \tfrac{1}{h_0} (a\log L + \log S_\epsilon)$ and thus for $\epsilon\in [0,\infty)$ and $L\in \N$ with $L> e$
\begin{align*}
\frac{S_{L,\epsilon,k}^{\frac{2}{k}} }{\log L} 
\le \tfrac{1}{h_0} (a + \log S_\epsilon). 
\end{align*}
This proves \eqref{eqn:sup_noise_terms_divided_by_log_L_epsilon} and \eqref{eqn:sup_noise_terms_divided_by_log_L}. 
By Jensen's inequality
\begin{align*}
\sup_{\epsilon \in [0,\infty)} E[ \exp( h \tfrac{1}{h_0} (a + \log S_\epsilon)) ] \le 
\exp( \tfrac{ha}{h_0}) \sup_{\epsilon \in [0,\infty)} E[   S_\epsilon  ]^{\frac{h}{h_0}}, 
\end{align*}
hence by \eqref{eqn:sum_L_of_L_power_expectation_moment_gen_fu} the latter is finite for $h\in [0,h_0]$, this concludes the proof. 
\end{proof}

\appendix

\section{Regularity of even extensions}

\begin{lemma}
\label{lemma:equivalence_full_space_and_neumann}
Let $\alpha \in \R$. 
There exists a $C>1$ such that for all $\zeta \in \cC_\fn^\alpha(Q_L)$
\begin{align*}
\frac{1}{C} \| \overline \zeta\|_{\cC^\alpha} 
\le \| \zeta\|_{\cC_\fn^\alpha} 
\le  C \| \overline \zeta\|_{\cC^\alpha},
\end{align*}
where $\overline \zeta = \sum_{k\in \N_0^2} \langle \zeta, \fn_{k,L} \rangle \overline \fn_{k,L}$ (notation as in \ref{obs:notation_neumann}). 
\end{lemma}
\begin{proof}
By definition the norm of $\zeta$ is equivalent to the norm of the even extension of $\zeta$ on the torus $\T_{2L}^2$ (see \cite{ChvZ21}). 
On the other hand, the $\cC^\alpha(\R^2)$-norm of $\overline \zeta$ is equivalent to the $\cC^\alpha(\T_{2L}^2)$-norm of the ``restriction'' of $\overline \zeta$ to the torus, see for example \cite[Theorem 3.8.1, p.195]{ScTr87} (observe that there is the condition $\chi \times \infty>d$, which for $\chi=0$ can be thought of as still valid). 
\end{proof}

\section{The resolvent equation}

Here we show the existence of and a bound for the solution $v$ of 
\begin{align*} 
	\left( \eta - \tfrac{1}{2} \Delta \right) v = f + \nabla v \cdot g \qquad \Leftrightarrow \qquad v= \sigma_{\eta} (\rD) (f + \nabla v \cdot g),
\end{align*}
for suitable distributions $f,g$, where $\sigma_{\eta} (z) := (\eta + 2 \pi^2 | z |^2)^{-1} $ and $\eta$ is a scalar.

Lemma~\ref{lemma:bound_on_sigma_lambda_operator} displays the regularizing effect of the Fourier multiplier $\sigma_\eta$, it follows by an application of \cite[Lemma 2.2]{BaChDa11}.  
\begin{calc}
It may help to observe the identity $\sigma_\lambda(z) = \lambda^{-1} \sigma_1 (\lambda^{-\frac12} z)$. Here the consequence of \cite[Lemma 2.2]{BaChDa11}. 

\newtheorem*{theorem*}{Theorem}

\begin{theorem*} 
Let $d\in \N$. 
There exists a $C>0$ such that the following holds. 
Let $\gamma,m \in \R$. 
If for all $\alpha \in \N_0^d$ with $|\alpha| \le 2 \lfloor 1+ \frac{d}{2} \rfloor$ there exist a $C_\alpha>0$ such that  $|\partial^\alpha \sigma(x)| \le C_\alpha |x|^{m-|\alpha|}$ for all $x \ne 0$, then 
\begin{align}
\label{eqn:fourier_multiplier_bound_on_C}
\|\sigma(\rD)w  \|_{\cC^{\gamma +m }(\R^d)} 
& \le \Big( C \sum_{\alpha \in \N_0^d: |\alpha|\le 2 \lfloor 1+ \frac{d}{2} \rfloor} C_\alpha \Big) \|w \|_{\cC^\gamma(\R^d)}.
\end{align}
\end{theorem*}
\end{calc}

\begin{lemma} 
\label{lemma:bound_on_sigma_lambda_operator}
There exists a $C>0$ only depending on $d$, such that for all $\alpha \in \R$, $\delta \in [0,2]$ and $\eta \ge 1$ 
  \begin{align*} 
  \| \sigma_{\eta} (\rD) v \|_{\cC^{\alpha + \delta}(\R^d)} \le C
     \eta^{- (1 - \frac{\delta}{2})} \| v \|_{\cC^\alpha(\R^d)} . 
     \end{align*}
\end{lemma}

\begin{proposition}
\label{prop:solution_v_lambda}
Let $\alpha \in (-\frac43,-1)$ and $\beta \in (-\alpha, 2\alpha + 4)$ and $f \in \cC^{2\alpha+2}(\R^d)$, $g \in
  \cC^{\alpha+1}(\R^d)$. 
Then there exists a $C > 0$ such that with 
\begin{align*} 
M \ge \max \{ \| f \|_{\cC^{2\alpha+2}}, \| g \|_{\cC^{\alpha+1}} \} , 
\end{align*}  
for all $\eta \ge C (1 + M)^{\frac{2}{2\alpha + 4 - \beta}}$ the equation
  \begin{align*} 
  v = \sigma_{\eta} (\rD) (f + \nabla v \cdot  g) 
     \end{align*}
  has a unique solution $v \in \cC^{\beta}(\R^d)$, and $ \| v \|_{\cC^{\beta}} \le M$.  Moreover, $v$ depends  continuously on $f$ and $g$ in the following sense. 
Suppose additionally that $f_\epsilon \in \cC^{2\alpha+2}(\R^d), g_\epsilon  \in \cC^{\alpha+1}(\R^d)$, $\|f_\epsilon\|_{\cC^{2\alpha+2}},\|g_\epsilon\|_{\cC^{\alpha+1}}\le M$ for $\epsilon >0$ and $f_\epsilon \rightarrow f$, $g_\epsilon \rightarrow g$ in $\cC^{2\alpha+2}(\R^d)$ respectively in $\cC^{\alpha+1}(\R^d)$.
For $\eta \ge C (1 + M)^{\frac{2}{2\alpha + 4 - \beta}}$ and with $v_{\epsilon}$ the solution to $v_{\epsilon} = \sigma_{\eta} (\rD) (f_\epsilon + \nabla v_{\epsilon} \cdot  g_\epsilon)$ we have $v_{\epsilon} \rightarrow v$ in $\cC^{\beta}(\R^d)$.
\end{proposition}

\begin{proof}
We use a Picard iteration: Let $\Phi \colon \cC^{\beta} \rightarrow \cC^{\beta}$, $\Phi (u) = \sigma_{\eta} (\rD) (f + \nabla u \cdot g)$, for $\eta > 0$. Let $C_1 >1$ be as $C$ is in Lemma~\ref{lemma:bound_on_sigma_lambda_operator} and $C_2>1$ be such that $\|\nabla u \cdot g\|_{\cC^{\alpha+1}} \le C_2 \| u\|_{\cC^{\beta}} \|g\|_{\cC^{\alpha+1}}$. 
Let $C= C_1 \cdot C_2$. A short computation shows that if $\| u \|_{\cC^{\beta}} \le M$, then 
  \begin{align*}
    \| \Phi (u) \|_{\cC^{\beta}} 
    \le    C M (1 + M) \eta^{- \frac{2\alpha+4-\beta}{2}}  .
  \end{align*} 
\begin{calc}
Indeed, with $\delta = \beta - (2\alpha + 2) $ (as $2\alpha + 2 \le \alpha +1$) 
 \begin{align*}
    \| \Phi (u) \|_{\cC^{\beta}} 
    & \le  C_1 \eta^{- (1 - \frac{\beta - (2\alpha + 2) }{2})} \| f + \nabla  u \cdot g \|_{\cC^{2\alpha + 2}} \\
  &   \le C_1 \eta^{- \frac{2\alpha + 4 - \beta }{2} } (\| f \|_{\cC^{2\alpha +2 }} +
    C_2 \| u \|_{\cC^{\beta}} \| g \|_{\cC^{\alpha+1}})\\
   &  \le C_1 \eta^{-  \frac{2\alpha + 4 - \beta }{2} } (M + C_2 M^2) \le
    C M (1 + M) \eta^{-  \frac{2\alpha + 4 - \beta }{2} }  .
  \end{align*}
  Observe that as $\beta \in (- \alpha, 2\alpha +4)$, $2 \alpha + 4 - \beta \in (0,3\alpha +4)$. 
\end{calc}
  So if $\eta \ge C (1 + M)^{\frac{2}{2\alpha + 4 - \beta}}$, then $\| \Phi (u) \|_{\cC^{\beta}} \le M$ and thus $\Phi$ leaves
  the ball in $\cC^{\beta}$ with radius $M$ and center $0$
  invariant. Moreover,  
 \begin{align*}
    \| \Phi (u) - \Phi (\tilde{u}) \|_{\cC^{\beta}} 
    & = \|     \sigma_{\eta} (\rD) ((\rD u - \rD \tilde{u}) \cdot g) \|_{\cC^{\beta}} \\
    \cand \begin{calc} 
    \le C_1 \eta^{- \frac{2\alpha + 4 - \beta }{2}} \| (\rD u - \rD \tilde{u}) \cdot g \|_{\cC^{2\alpha + 2 }} 
    \end{calc}\cnewline
    \cand \begin{calc}
    \le C \eta^{- \frac{2\alpha + 4 - \beta }{2}} \| u - \tilde{u}   \|_{\cC^{\beta}} \| g \|_{\cC^{\alpha+1}} 
	\end{calc} \cnewline    
    & \le C M \eta^{- \frac{2\alpha + 4 - \beta }{2}}  \| u - \tilde{u} \|_{\cC^{\beta}} .
  \end{align*}
    So $\Phi$ is a contraction and it has a unique fixed point. The continuity is shown as in \cite[Proposition 2.4]{PevZ}. 
\end{proof}

\textbf{Acknowledgements.} 
The authors are grateful to T. Matsuda for feedback on a previous draft.  
The authors are also grateful to the anonymous referee for their valuable feedback, suggestions and careful reading. 
This work was supported by the German Science Foundation (DFG) via the Forschergruppe FOR2402 ``Rough paths, stochastic partial differential equations and related topics''. WK and WvZ were supported by the DFG through SPP1590 ``Probabilistic Structures in Evolution''. 
NP thanks the DFG for financial support through the Heisenberg programme.
The main part of the work was done while NP was employed at Humboldt-Universit\"at zu Berlin and Max Planck Institute for Mathematics in the Sciences, Leipzig.

\bibliographystyle{abbrv}
\bibliography{references}

\begin{thebibliography}{10}

\bibitem{AlCh15}
R.~Allez and K.~Chouk.
\newblock The continuous anderson hamiltonian in dimension two.
\newblock Preprint available at \url{https://arxiv.org/abs/1511.02718}.

\bibitem{As16}
A.~Astrauskas.
\newblock From extreme values of i.i.d. random fields to extreme eigenvalues of
  finite-volume {A}nderson {H}amiltonian.
\newblock {\em Probab. Surv.}, 13:156--244, 2016.

\bibitem{BaChDa11}
H.~Bahouri, J.-Y. Chemin, and R.~Danchin.
\newblock {\em Fourier analysis and nonlinear partial differential equations},
  volume 343 of {\em Grundlehren der Mathematischen Wissenschaften [Fundamental
  Principles of Mathematical Sciences]}.
\newblock Springer, Heidelberg, 2011.

\bibitem{Bi68}
P.~Billingsley.
\newblock {\em Convergence of probability measures}.
\newblock Wiley Series in Probability and Statistics: Probability and
  Statistics. John Wiley \& Sons, Inc., New York, second edition, 1999.
\newblock A Wiley-Interscience Publication.

\bibitem{BiKodS18}
M.~Biskup, W.~K{\"{o}}nig, and R.~S. dos Santos.
\newblock Mass concentration and aging in the parabolic {A}nderson model with
  doubly-exponential tails.
\newblock {\em Probab. Theory Related Fields}, 171(1-2):251--331, 2018.

\bibitem{CaCh18}
G.~Cannizzaro and K.~Chouk.
\newblock Multidimensional {SDE}s with singular drift and universal
  construction of the polymer measure with white noise potential.
\newblock {\em Ann. Probab.}, 46(3):1710--1763, 2018.

\bibitem{Ch14}
X.~Chen.
\newblock Quenched asymptotics for {B}rownian motion in generalized {G}aussian
  potential.
\newblock {\em Ann. Probab.}, 42(2):576--622, 2014.

\bibitem{ChvZ21}
K.~Chouk and W.~van Zuijlen.
\newblock Asymptotics of the eigenvalues of the anderson hamiltonian with white
  noise potential in two dimensions.
\newblock {\em Ann. Probab.}, 49(4):1917--1964, 2021.

\bibitem{DeDi16}
F.~Delarue and R.~Diel.
\newblock Rough paths and 1d {SDE} with a time dependent distributional drift:
  application to polymers.
\newblock {\em Probab. Theory Related Fields}, 165(1-2):1--63, 2016.

\bibitem{DuKo10}
J.~J. Duistermaat and J.~A.~C. Kolk.
\newblock {\em Distributions}.
\newblock Cornerstones. Birkh\"auser Boston, Inc., Boston, MA, 2010.
\newblock Theory and applications, Translated from the Dutch by J. P. van Braam
  Houckgeest.

\bibitem{DuLa17}
L.~Dumaz and C.~Labb\'{e}.
\newblock Localization of the continuous {A}nderson {H}amiltonian in 1-{D}.
\newblock {\em Probab. Theory Related Fields}, 176(1-2):353--419, 2020.

\bibitem{Fr75}
A.~Friedman.
\newblock {\em Stochastic differential equations and applications. {V}ol. 1}.
\newblock Academic Press [Harcourt Brace Jovanovich, Publishers], New
  York-London, 1975.
\newblock Probability and Mathematical Statistics, Vol. 28.

\bibitem{GaKoMo00}
J.~G\"{a}rtner, W.~K\"{o}nig, and S.~A. Molchanov.
\newblock Almost sure asymptotics for the continuous parabolic {A}nderson
  model.
\newblock {\em Probab. Theory Related Fields}, 118(4):547--573, 2000.

\bibitem{GuImPe15}
M.~Gubinelli, P.~Imkeller, and N.~Perkowski.
\newblock Paracontrolled distributions and singular {PDE}s.
\newblock {\em Forum Math. Pi}, 3:e6, 75, 2015.

\bibitem{GuPe17}
M.~Gubinelli and N.~Perkowski.
\newblock K{PZ} reloaded.
\newblock {\em Comm. Math. Phys.}, 349(1):165--269, 2017.

\bibitem{GuUgZa20}
M.~Gubinelli, B.~Ugurcan, and I.~Zachhuber.
\newblock Semilinear evolution equations for the {A}nderson {H}amiltonian in
  two and three dimensions.
\newblock {\em Stoch. Partial Differ. Equ. Anal. Comput.}, 8(1):82--149, 2020.

\bibitem{Ha14}
M.~Hairer.
\newblock A theory of regularity structures.
\newblock {\em Invent. Math.}, 198(2):269--504, 2014.

\bibitem{HaLa15}
M.~Hairer and C.~Labb\'e.
\newblock A simple construction of the continuum parabolic {A}nderson model on
  {${\bf R}^2$}.
\newblock {\em Electron. Commun. Probab.}, 20:no. 43, 11, 2015.

\bibitem{HyvNVeWe16}
T.~Hyt\"{o}nen, J.~van Neerven, M.~Veraar, and L.~Weis.
\newblock {\em Analysis in {B}anach spaces. {V}ol. {I}. {M}artingales and
  {L}ittlewood-{P}aley theory}, volume~63 of {\em Ergebnisse der Mathematik und
  ihrer Grenzgebiete. 3. Folge. A Series of Modern Surveys in Mathematics
  [Results in Mathematics and Related Areas. 3rd Series. A Series of Modern
  Surveys in Mathematics]}.
\newblock Springer, Cham, 2016.

\bibitem{KaSh91}
I.~Karatzas and S.~E. Shreve.
\newblock {\em Brownian motion and stochastic calculus}, volume 113 of {\em
  Graduate Texts in Mathematics}.
\newblock Springer-Verlag, New York, second edition, 1991.

\bibitem{Ko16}
W.~K{\"{o}}nig.
\newblock {\em The parabolic {A}nderson model. Random walk in random
  potential.}
\newblock Pathways in Mathematics. Birkh\"auser/Springer, [Cham], 2016.

\bibitem{La19}
C.~Labb\'{e}.
\newblock The continuous {A}nderson {H}amiltonian in {$d\leq 3$}.
\newblock {\em J. Funct. Anal.}, 277(9):3187--3235, 2019.

\bibitem{GL}
P.~Y.~G. Lamarre.
\newblock Phase transitions in asymptotically singular anderson hamiltonian and
  parabolic model.
\newblock Preprint available at \url{https://arxiv.org/abs/2008.08116}.

\bibitem{LG16}
J.-F. Le~Gall.
\newblock {\em Brownian motion, martingales, and stochastic calculus}, volume
  274 of {\em Graduate Texts in Mathematics}.
\newblock Springer, [Cham], french edition, 2016.

\bibitem{MaPe19}
J.~Martin and N.~Perkowski.
\newblock Paracontrolled distributions on {B}ravais lattices and weak
  universality of the 2d parabolic {A}nderson model.
\newblock {\em Ann. Inst. Henri Poincar\'{e} Probab. Stat.}, 55(4):2058--2110,
  2019.

\bibitem{Nu09}
D.~Nualart.
\newblock {\em Malliavin calculus and its applications}, volume 110 of {\em
  CBMS Regional Conference Series in Mathematics}.
\newblock Published for the Conference Board of the Mathematical Sciences,
  Washington, DC; by the American Mathematical Society, Providence, RI, 2009.

\bibitem{Pa83}
A.~Pazy.
\newblock {\em Semigroups of linear operators and applications to partial
  differential equations}, volume~44 of {\em Applied Mathematical Sciences}.
\newblock Springer-Verlag, New York, 1983.

\bibitem{PevZ}
N.~Perkowski and W.~B. van Zuijlen.
\newblock Quantative heat kernel estimates for diffusions with distributional
  drift.
\newblock Preprint available at \url{http://arxiv.org/abs/2009.10786}.

\bibitem{ScTr87}
H.-J. Schmeisser and H.~Triebel.
\newblock {\em Topics in {F}ourier analysis and function spaces}.
\newblock A Wiley-Interscience Publication. John Wiley \& Sons, Ltd.,
  Chichester, 1987.

\bibitem{Sz98}
A.-S. Sznitman.
\newblock {\em Brownian motion, obstacles and random media}.
\newblock Springer Monographs in Mathematics. Springer-Verlag, Berlin, 1998.

\end{thebibliography}

\end{document}